\documentclass{amsart}

\usepackage{amsmath}
\usepackage{amssymb}
\usepackage{mathrsfs}
\usepackage{stmaryrd}
\usepackage{mathtools}
\usepackage{amsfonts}

\usepackage{hyperref}
\usepackage{enumitem}

\usepackage{comment}
\usepackage{graphicx}
\usepackage[table]{xcolor}
\usepackage[all]{xy}
\usepackage{tikz}
\usepackage{tikz-cd}
\usepackage{hhline}

\usepackage[alphabetic,nobysame]{amsrefs}
\usepackage{float}

\newtheorem{theorem}{Theorem}[section]
\newtheorem{lemma}[theorem]{Lemma}
\newtheorem{proposition}[theorem]{Proposition}
\newtheorem{cor}[theorem]{Corollary}
\newtheorem{claim}[theorem]{Claim}

\theoremstyle{definition}
\newtheorem{definition}[theorem]{Definition}
\newtheorem{remark}[theorem]{Remark}
\newtheorem{example}[theorem]{Example}

\numberwithin{equation}{section}

\renewcommand{\setminus}{\smallsetminus}

\newcommand{\al}{\alpha}
\newcommand{\epsi}{\varepsilon}

\newcommand{\inv}{^{-1}}
\newcommand{\into}{\hookrightarrow} 
\newcommand{\onto}{\twoheadrightarrow} 
\newcommand{\GL}{\mathrm{GL}} 
\DeclareMathOperator{\coker}{coker}
\DeclareMathOperator{\rank}{rank}

\newcommand{\Gm}{\mathbb{G}_\mathrm{m}}

\DeclareMathOperator\Aut{Aut}

\DeclareMathOperator{\Hom}{Hom} 
\newcommand{\Ext}{\mathrm{Ext}}
\newcommand{\cExt}{\mathcal{E}xt}

\DeclareMathOperator{\Spec}{Spec} 
\DeclareMathOperator{\Pic}{Pic}

\newcommand{\NE}{\overline{\mathrm{NE}}}
\DeclareMathOperator{\Nef}{Nef}

\def\pow#1{ \llbracket  #1 \rrbracket }

\newcommand{\stack}[1]{ \mathcal{M}^\mathrm{Kss}_{#1} } 
\newcommand{\modspace}[1]{M^\mathrm{Kps}_{#1}} 
\newcommand{\stackstab}[1]{ \mathcal{M}^\mathrm{Ks}_{#1} } 
\newcommand{\spacestab}[1]{M^\mathrm{Ks}_{#1} } 

\newcommand{\artC}{(\mathrm{Art}_{\mathbb{C}})}

\newcommand{\set}{(\mathrm{Set})}

\DeclareMathOperator{\Spf}{Spf}
\newcommand{\Def}[1]{\mathrm{Def}_{#1}}

\newcommand{\DefqG}[1]{\mathrm{Def}^{\mathrm{qG}}_{#1}}
\newcommand{\TTqG}[2]{\mathbb{T}^{\mathrm{qG}, #1}_{#2}}
\newcommand{\cTqG}[2]{\mathscr{T}^{\mathrm{qG}, #1}_{#2}}

\def\spanR#1{\mathrm{span}_\mathbb{R} \left( #1 \right)} 
\newcommand{\cone}{\mathrm{cone}}
\def\cone#1{\mathrm{cone} \left\{ #1 \right\}}
\newcommand{\TV}[2]{\operatorname{TV}(#1,#2)}

\newcommand\cA{\mathcal{A}}

\newcommand\cM{\mathcal{M}}
\newcommand\cO{\mathcal{O}}

\newcommand\cU{\mathcal{U}}
\newcommand\cV{\mathcal{V}}
\newcommand\cX{\mathcal{X}}
\newcommand\cY{\mathcal{Y}}

\newcommand\cT{\mathscr{T}}

\renewcommand\AA{\mathbb{A}}
\newcommand\CC{\mathbb{C}}
\newcommand\FF{\mathbb{F}}

\newcommand{\LL}{\mathbb{L}}

\newcommand\PP{\mathbb{P}}
\newcommand\QQ{\mathbb{Q}}
\newcommand\RR{\mathbb{R}}
\newcommand\TT{\mathbb{T}}
\newcommand\ZZ{\mathbb{Z}}

\newcommand\rH{\mathrm{H}}

\newcommand\rN{\mathrm{N}}

\newcommand{\rT}{\mathrm{T}}

\newcommand\rme{\mathrm{e}}
\newcommand\rmi{\mathrm{i}}

\newcommand{\frakm}{\mathfrak{m}}
\newcommand{\frakp}{\mathfrak{p}}
\newcommand{\frakq}{\mathfrak{q}}

\title[K-stable Fano varieties can have obstructed deformations]{K-stable Fano varieties can have \\ obstructed deformations}

\author{Andrea Petracci}
\address{Dipartimento di Matematica, Università di Bologna, Piazza di Porta San Donato 5, 40126 Bologna, Italy}
\email{a.petracci@unibo.it}

\author{Alice Villa}
\address{Dipartimento di Matematica, Università di Bologna, Piazza di Porta San Donato 5, 40126 Bologna, Italy}
\email{alice.villa6@unibo.it}

\begin{document}
	
	\begin{abstract}
		We exhibit examples of K-stable Fano varieties with obstructed qG-deformations; more precisely, their miniversal qG-deformation spaces can be non-reduced. As a consequence, we show that both the K-moduli stack and the K-moduli space can be non-reduced, and hence non-smooth, at points corresponding to K-stable Fano varieties.
	\end{abstract}
	
	\maketitle
	
	\section{Introduction}
	
	We work over an algebraically closed field of characteristic zero, denoted by $\CC$.
	A Fano variety is a normal projective variety whose anticanonical divisor is $\QQ$-Cartier and ample; in particular, it is allowed to be singular.
	All Fano varieties considered in this article will have klt singularities, although we will not use the term `$\QQ$-Fano' variety.
	
	Fano varieties occupy a prominent role in algebraic geometry for many reasons.
	A major source of recent progress has been K-stability \cite{xu_book}, which provides the algebro-geometric counterpart to the existence problem for Kähler--Einstein metrics on Fano varieties through the Yau--Tian--Donaldson correspondence.
	For example, K-stability has been spectacularly applied to construct moduli spaces for Fano varieties (see \cite{xu_book, BX19, openness_of_K-semistab, LXZ} and references therein):
	more precisely,
	for every $n \in \ZZ_{\geq 1}$ and for every $v \in \QQ_{>0}$,
	qG-families of K-semistable Fano varieties
	of dimension $n$ and anticanonical volume $v$
	form an algebraic stack,
	denoted by $\stack{n,v}$, of finite type over $\CC$;
	moreover, the stack $\stack{n,v}$ admits a good moduli space, denoted by $\modspace{n,v}$, which is projective over $\CC$ and whose $\CC$-points
	are in natural bijection with isomorphism classes of K-polystable Fano varieties over $\CC$ of dimension $n$ and anticanonical volume $v$.
	Let $\stackstab{n,v}$ denote the open substack of $\stack{n,v}$ parametrising K-stable Fano varieties of dimension $n$ and anticanonical volume $v$; $\stackstab{n,v}$ is a Deligne--Mumford stack whose coarse moduli space $\spacestab{n,v}$ is an open subscheme of $\modspace{n,v}$.
	We have a cartesian diagram
	\[
	\xymatrix{
		\stackstab{n,v} \ar@{^{(}->}[r] \ar[d] & \stack{n,v} \ar[d] \\
		\spacestab{n,v} \ar@{^{(}->}[r] & \modspace{n,v}}
	\]
	where the horizontal maps are open embeddings,
	the left vertical map is a coarse moduli space and
	the right vertical map is a good moduli space.
	Considerable work has been devoted to understanding the geometry of these stacks and spaces (e.g.\ \cite{LiuXu_cubic3fold, Liu_cubic4fold, LiuZhao, 1-dim_moduli_of_Fano, cheltsov2026degenerationsellipticquarticskmoduli, zhao219, klpz}).
	Here we are concerned with their local geometry, and in particular with smoothness.
	
	By Kodaira--Akizuki--Nakano vanishing it is well known that smooth Fano varieties have unobstructed deformations.
	This implies the following two statements:
	\begin{enumerate}
		\item[(i)] the stack $\stack{n,v}$ is smooth at any point corresponding to a smooth K-semistable Fano variety;
		\item[(ii)] the scheme $\modspace{n,v}$ is normal with rational singularities at any point corresponding to a smooth K-polystable Fano variety, because its automorphism group is reductive by \cite{reductivity}
		and reductive group quotients of smooth schemes have rational singularities by \cite{boutot}.
	\end{enumerate}
	The results (i) and (ii) hold if we remove the assumption on the smoothness of the Fano variety and we impose one of the following conditions which ensure that the qG-deformations of a Fano variety $X$ are unobstructed:
	\begin{itemize}
		\item $\dim X =2$ and $X$ is klt by \cite[\S3]{hacking_prokhorov} or \cite[Lemma~6]{procams};
		\item $\dim X = 3$ and $X$ has terminal singularities by \cite[Theorem~2.11]{Sano} (see also \cite{Namikawa, Minagawa});
		\item $\dim X \geq 3$ and $X$ has lci $1$-DuBois singularities by \cite[Theorem~4.5]{FL25};
		\item $\dim X \geq 3$ and $X$ has isolated lci singularities by \cite[Theorem~3.1]{tenie_preprint} (see also \cite[Theorem~3.17]{friedman_preprint}).
	\end{itemize}
	
	In \cite{KP} Kaloghiros and the first author exhibited examples of K-polystable Fano varieties with obstructed qG-deformations, and consequently they proved the existence of singular points on the K-moduli stack of K-semistable Fano varieties and the existence of non-reduced or non-unibranch points on the K-moduli space of K-polystable Fano varieties.
	Since these varieties are toric, their automorphism groups contain the acting torus and are therefore positive-dimensional; in particular, they are not K-stable.
	Similar toric constructions have been used in \cite{Pet22a, petracci_agplus}
	to deduce further pathologies about the smoothness of K-moduli at some non-smooth toric Fano varieties.
	
	It is therefore natural to ask whether obstructedness and singularities on K-moduli are phenomena confined to the non-K-stable locus, or whether they can occur for K-stable Fano varieties as well.
	Our main result gives a negative answer to the former possibility.
	
	\begin{theorem}\label{thm:main}
		Fix $n \in \ZZ$ such that $n=3$ or $n \geq 5$.
		Then there exists a K-stable Fano variety $Y$ of dimension $n$
		such that the automorphism group of $Y$ is trivial and
		the base of the miniversal qG-deformation of $Y$ is non-reduced; in particular, $Y$ has obstructed qG-deformations.
	\end{theorem}
	
	A consequence is the fact that singularities of K-moduli of Fano varieties do not occur only at the `boundary' $\stack{n,v} \setminus \stackstab{n,v}$, i.e.\ at non-K-stable points: they can already appear inside the open K-stable locus:
	
	\begin{cor} \label{cor:main}
		Fix $n \in \ZZ$ such that $n=3$ or $n \geq 5$.
		Then there exists $v \in \QQ_{>0}$ such that the K-moduli stack $\stackstab{n,v}$ and the K-moduli space $\spacestab{n,v}$ of K-stable Fano $n$-folds with anticanonical volume $v$ are not reduced.
	\end{cor}
	
	\begin{remark} \label{rmk:n=4}
		It is very likely that Theorem~\ref{thm:main} and Corollary~\ref{cor:main} are valid also for $n=4$, but this cannot be immediately proved with the techniques described in this paper.
	\end{remark}
	
	\subsection*{Idea of the proof of Theorem~\ref{thm:main}}
	We begin by analysing the deformation theory of a specific toric K-polystable Fano $3$-fold $X$ (see Theorem~\ref{thm:tor_var_X}): the base of the miniversal qG-deformation of $X$ has four irreducible components, three of which are non-reduced and one is generically reduced.
	Since $X$ is toric, its automorphism group is infinite and thus $X$ is not K-stable.
	
	Heuristically, one would like to choose $Y$ corresponding to a nearby point on one of the non-reduced components of the analytic qG-deformation space of $X$.
	Openness of versality then suggests that the miniversal qG-deformation space of $Y$ should still be non-reduced.
	
	Rather than making this heuristic precise in the complex-analytic category, we work entirely algebraically.
	Since $(-K_X)^3 = 11/3$, the toric Fano $3$-fold $X$ determines a point $[X]$ of the algebraic stack $\stack{3, 11/3}$.
	Set $G := \Aut(X)$.
	The idea is that we consider an algebraisation $S$ of the formal qG-deformation space of $X$; more precisely, we consider an affine scheme $S$ of finite type over $\CC$, a point $s \in S(\CC)$, an action of $G$ on $S$ such that
	$s$ is fixed by $G$
	and the completed local ring of $S$ at $s$ is $G$-equivariantly isomorphic to the hull of the qG-deformation functor of $X$.
	There is no a priori reason for this algebraisation to carry a morphism $S \to \stack{3, 11/3}$ which maps $s$ to $[X]$,
	but by Artin's algebraisation (or more precisely by its generalisation to stacks \cite{AHR20, AHR26}) the two stacks $\stack{3, 11/3}$ and $[S/G]$ are étale-locally equivalent at the points $[X]$ and $s$, respectively.
	Our stack-theoretic result (Theorem~\ref{thm:artin_stack}), which may be of independent interest, allows us to transfer suitable properties of closed points of $S$ which are close to $s$ to
	points in the stack $\stack{3, 11/3}$ which are in an arbitrarily small neighbourhood of $[X]$.
	
	In this way, we obtain a Fano $3$-fold $Y$ with trivial automorphism group whose miniversal qG-deformation space is non-reduced and which is K-polystable (see Theorem~\ref{thm:result_n=3}).
	More specifically, the base of the miniversal qG-deformation of $Y$ is the formal spectrum of the non-reduced ring
	\[
	\CC \pow{w_1,\dots,w_9}/(w_1^2, w_2^2).
	\]
	Since $\Aut(Y)$ is trivial and $Y$ is K-polystable, we have that $Y$ is K-stable. This proves Theorem~\ref{thm:main} in dimension $n=3$.
	
	For $n \geq 5$, it is enough to consider the product of $Y$ with a general hypersurface of degree $n-2$ in $\PP^{n-2}$.
	Unfortunately, for $n=4$ this construction does not work because the automorphism group of $\PP^1$ is infinite (see Remark~\ref{rmk:n=4}).
	
	\begin{remark}
		The obstructed toric K-polystable Fano $3$-fold $X$ we start with
		did not appear in \cite{KP} and is new.
		This is because the two obstructed toric K-polystable Fano $3$-folds exhibited in \cite{KP} are not suitable for the construction above: the base of the miniversal qG-deformation of the $3$-fold in \cite[Theorem~1.2]{KP} is a fat point, so it does not deform to anything,
		whereas every small deformation of the $3$-fold in \cite[Theorem~1.1]{KP} does not satisfy all the three conditions of having obstructed deformations, having finite automorphism group and being K-polystable.
	\end{remark}
	
	\begin{remark} \label{rem:tecnica_differente}
		It is likely that Theorem~\ref{thm:main} (even if the dimension $n$ is $4$) can be proved with a technique different from what we sketched above.
		Instead of deforming an obstructed K-polystable toric Fano variety, one could take a finite quotient of it to obtain a K-stable Fano variety with obstructed qG-deformations.
		Nevertheless, we did not pursue this construction here.
	\end{remark}
	
	\subsection*{Outline}
	The paper is organised as follows.
	
	In \S\ref{sec:prel} we recall the preliminaries on deformation theory, toric varieties, qG-deformations, and K-stability that will be used throughout the paper.
	
	In \S\ref{sec:affine_toric_var} we study an affine toric variety associated with the cone over a specific lattice octagon.
	
	In \S\ref{sec:toric_var_X} we construct a toric K-polystable Fano 3-fold $X$ and study the hull $R$ of its qG-deformation functor $\DefqG{X}$, its automorphism group $\Aut(X)$ and how it acts on $R$.
	
	In \S\ref{sec:result_n=3} we use the previous study of the qG-deformation theory of $X$ to construct a K-stable Fano 3-fold $Y$ with trivial automorphism group and obstructed qG-deformations. The construction relies on a result about algebraic stacks proved in the appendix.
	
	In \S\ref{sec:final_proofs} we prove our main results, namely Theorem~\ref{thm:main} and Corollary~\ref{cor:main}. We first establish a result on automorphism groups of products of Fano varieties, which allows us to extend the $3$-dimensional construction to every dimension $n\geq 5$.
	
	Finally, in the appendix \ref{sec:stack} we prove a result on the existence of points in an algebraic stack with some specific properties. This result is a stack-theoretic tool needed in the construction of $Y$.

	\subsection*{Notation and conventions}
	We work over an algebraically closed field of characteristic zero, denoted by $\CC$.
	Every toric variety considered in this article is normal.
	An (algebraic) torus is an algebraic group isomorphic to $(\Gm)^n$. We will consider several representations of a torus and we will interchangeably use the words `degree' and `weight' for its characters, i.e.\ its $1$-dimensional representations.
	
	\subsection*{Acknowledgements}
	The first author is very grateful to Jarod Alper, Luca Battistella, Davide Gori and Angelo Vistoli for fruitful conversations about the topics treated in the appendix.

	Both authors acknowledge partial financial support from INdAM GNSAGA ``Gruppo Nazionale per le Strutture Algebriche, Geometriche e le loro Applicazioni’’ and from the European Union - NextGenerationEU under the National Recovery and Resilience Plan (PNRR) - Mission 4 Education and research - Component 2 From research to business - Investment 1.1, Prin 2022, ``Geometry of algebraic structures: moduli, invariants, deformations'', DD N.~104, 2/2/2022, proposal code 2022BTA242 - CUP J53D23003720006.
	
	Moreover, A.V.\ acknowledges partial financial support from  the European Union - NextGenerationEU under the National Recovery and Resilience Plan (PNRR) - Mission 4 Education and research -
	Component 2 From research to business - Investment 1.1 Notice Prin 2022 - DD N.
	104 del 2/2/2022, “Symplectic varieties: their interplay with Fano manifolds and
	derived categories”, proposal code 2022PEKYBJ – CUP J53D23003840006.
	
	\subsection*{AI usage}
	Generative AI was used to assist with typesetting, to obtain some hints for the appendix, and to develop Remark~\ref{rem:tecnica_differente}, which might suggest a possible alternative approach to proving Theorem~\ref{thm:main}.

	\section{Preliminaries} \label{sec:prel}
	
	In this section we present a collection of notions and well-known results that will be used in the subsequent parts of the paper.
	
	\subsection{Deformation theory} \label{sec:preliminaries_deformation_theory}
	
	We recall some basic notions and concepts from deformation theory, which will be used 
	later. For further details, we refer the reader to \cite{HarDef, Sernesi_def, TV}.
	
	For a scheme $X$ of finite type over $\CC$, we denote by
	\[
	\Def{X} \colon \artC \longrightarrow \set
	\] 
	its \emph{deformation functor}, where $\artC$ is the category of  local finite $\CC$-algebras with residue field $\CC$. Its \emph{tangent space}, which parametrises first-order deformations of $X$, is
	\[
	\rT\Def{X} := \Def{X}\bigl(\CC[t]/(t^2)\bigr).
	\]
	
	A complete noetherian local $\CC$-algebra $R$ with residue field $\CC$ is called a \emph{hull} of $\Def{X}$ if there exists a smooth natural transformation $h_R\to \Def{X}$ inducing an isomorphism on tangent spaces, where
	\[
	h_R(\cdot):=\Hom (R, \cdot).
	\]
	If a hull $R$ exists, it is unique up to in general non-unique isomorphism and the formal spectrum $\Spf(R)$ is called the \emph{base of the miniversal deformation} (or \emph{Kuranishi space}) of $X$.
	
	Using the cotangent complex $\LL_X$ of $X$ over $\CC$, one can define the $\CC$-vector spaces
	\[
	\TT_{X}^i:= \Ext^i(\LL_X,\cO_X)
	\]
	and the coherent sheaves
	\[
	\cT^i_X:= \cExt^i(\LL_X,\cO_X).
	\]
	We have $\TT_X^1 \simeq \rT\Def{X}$.
	Let $\Omega_X^1$ be the $\cO_X$-module of Kähler differentials of $X$ over $\CC$. If $X$ is reduced, then $\TT_X^1 \simeq \Ext^1(\Omega_X^1,\cO_X)$. If $X$ is normal, then $\TT_X^2 \simeq \Ext^2(\Omega_X^1,\cO_X)$.
	The analogous identifications hold for the sheaves $\cT_X^1$ and $\cT_X^2$.
	
	One can define the notion of obstruction space of a deformation functor. For example, $\TT_{X}^2$ is an obstruction space for $\Def{X}$.
	Moreover, if $k:= \dim_\CC \TT^{1}_X < \infty$ and $\TT^2_X=0$, then the hull of $\Def{X}$ is isomorphic to $\CC \pow{t_1, \dots, t_k}$ and we say that $X$ is \emph{unobstructed}.

	\subsection{Toric varieties} \label{sec:preliminaries_toric_geo}
	
	Toric varieties are algebraic varieties which can be constructed in a combinatorial way using cones and polytopes. The combinatorial properties of the latter translate into geometric properties of the associated varieties.
	For a complete and detailed overview of this topic, see \cite{CLS, Fulton}.
	
	A lattice is a free abelian group of finite rank.
	Let $N$ be a lattice and let $M = \Hom_{\ZZ}(N, \ZZ)$ be its dual lattice.  
	We denote by 
	\[ 
	\langle \cdot, \cdot \rangle \colon M \times N \longrightarrow \ZZ 
	\] 
	both the duality pairing and its extension $\langle \cdot, \cdot \rangle \colon M_{\RR} \times N_{\RR} \to \RR$.
	
	Let $\sigma \subseteq N_\RR$ be a strongly convex rational polyhedral cone.  
	The associated affine toric variety is
	\[
	\TV{\sigma}{N} = \Spec \CC[\sigma^\vee \cap M],
	\]
	where $\sigma^\vee \subseteq M_\RR$ is the dual cone. For each $m \in \sigma^\vee \cap M$, we denote with $\chi^m$ the corresponding element in the monoid algebra $\CC[\sigma^\vee \cap M]$.
	
	A \emph{polytope} in $N_\RR$ is the convex hull of a finite subset of $N_\RR$.
	
	A \emph{Fano polytope} in $N$ is a polytope $P \subseteq N_\RR$ such that $\dim P = \rank N$, the origin is in the interior of $P$ and every vertex of $P$ is a primitive lattice point.
	Whenever we have a full-dimensional polytope $P \subseteq N_\RR $ such that $0 \in \operatorname{Int}(P)$ and all its vertices are in $N_\QQ $, then we can consider the \emph{face fan} of $P$: this is the complete fan whose cones are those constructed over the faces of $P$. Under these hypotheses we can associate a toric variety to $P$, which is Fano whenever $P$ is a Fano polytope.
	
	Given $P \subset N_\RR $ a Fano polytope, the \textit{polar} of $P$ is
	\[
	P^\circ := \{ m \in M_\RR \mid \forall p \in P, \ \langle m, p \rangle \ge -1 \} \subset M_\RR.
	\]
	
	The deformation theory of affine toric varieties has been a subject of extensive research, with fundamental contributions provided by Altmann. A particularly interesting class of singularities consists of those arising from lattice polygons.
	
	\begin{theorem}[{\cite{Alt97}}] \label{prop:Alt_hull}
		Let $F$ be a lattice polygon in $\ZZ^2$ such that every lattice point on the boundary of $F$ is a vertex, and let $V$ denote the affine toric variety over $\CC$ associated with the cone $\RR_{\geq 0}(F \times \{1\})$ in the lattice $\ZZ^3$. Let $R$ be the hull of the deformation functor $\Def{V}$. 
		The following statements hold.
		\begin{enumerate}
			\item $V$ is $3$-dimensional, Gorenstein and its singular locus is either empty or consists of a single point.
			\item $R$ is the completion of the graded algebra $A_F$, introduced in \cite{Alt97} (see also \cite[Construction~4.1]{Pet22a}), at its maximal homogeneous ideal. 
			\item The embedding dimension of $R$, which is the dimension of the $\CC$-vector space $\TT_V^1$, is equal to $\ell - 3$, where $\ell$ is the number of vertices of $F$.
			\item There is a canonical one-to-one correspondence between the minimal prime ideals $\frakp$ of $R$ and the maximal Minkowski decompositions of the polygon $F$. 
			Furthermore, for each integer $s \geq 0$, there exists a canonical one-to-one correspondence between the irreducible components of $\Spec R$ of dimension $s$ and the maximal Minkowski decompositions of $F$ with $s+1$ summands.
		\end{enumerate}
	\end{theorem}
	Every affine toric $3$-fold with an isolated Gorenstein singularity arises via the construction above from a lattice polygon $F$ such that every lattice point on the boundary of $F$ is a vertex.
	
	\medskip
	
	In general, if $X$ is a toric variety, the algebraic torus $T_N = N \otimes_{\ZZ} \CC^*$ acts on $X$, and also acts on $\TT^{1}_{X}$ and $\TT^2_X$, which thus become $T_N$-representations.  
	In particular, we have
	\[
	\TT^1_X = \bigoplus_{m \in M} \TT^1_X(m),
	\]
	where $\TT^1_X(m)$ is the eigenspace where the torus $T_N$ acts with weight $m\in M$.
	
	Deformations of toric varieties are extensively studied (see \cite{Alt97, AltMin, ilten_vollmert, mavlyutov, petracci_homog, acf, ilten_turo, ilten_robins, filip_alg, filip_adv} and references therein).

	\subsection{Index 1 cover and qG-deformations} \label{sec:preliminaries_qG_index1}
	
	Let $X$ be a normal $\QQ$-Gorenstein variety over $\CC$. We introduce the construction of the index 1 cover, which is a local construction (see \cite[3.6-7]{Reid}, \cite[4.5]{Matsuki}). 
	For any closed point $x \in X$, there exists an open neighbourhood $U$ of $x$ such that $rK_U$ is linearly trivial, where $r$ is the local Gorenstein index of $X$ in $x$, i.e.\ the smallest $r\in\ZZ, r\geq1$ such that $rK_X$ is Cartier near $x$. This choice of trivialisation induces an isomorphism of $\cO_U$-modules 
	$\cO_{U}(rK_U) \simeq \cO_U.$
	Using this isomorphism, we define a sheaf of $\cO_U$-algebras
	\[
	\cA = \bigoplus_{i=0}^{r-1} \cO_U(iK_U).
	\]
	The local index 1 cover at $x$ is then defined as the finite $\mu_r$-cover:
	\[
	\pi \colon Y = \Spec_U(\cA) \longrightarrow U.
	\]
	Since $U$ is isomorphic to the schematic quotient $Y/\mu_r$, this construction realises $X$ in a neighbourhood of any point as a quotient of a Gorenstein variety by a cyclic group. If we consider the various quotient stacks $[Y/\mu_r]$ when the point in $X$ varies, they fit together into a Deligne--Mumford stack $\cX$, called the \emph{index $1$ cover stack}, with a coarse moduli morphism $\cX \to X$ which is an isomorphism over the
	Gorenstein locus of $X$.
	
	An \emph{infinitesimal qG-deformation} of $X$ is an infinitesimal deformation of $X$ which is locally induced by an equivariant deformation of the index 1 cover of $X$. This condition defines a subfunctor $\DefqG{X} \subseteq \Def{X}$. Equivalently, qG-deformations are those which satisfy Kollar's condition (see \cite{Kollar, AbraHass}).
	Analogously to the classical case, using the cotangent complex of the index $1$ cover stack $\cX$, it is possible to give the definition of the $\CC$-vector space $\TTqG{i}{X}$ and of the coherent sheaf $\cTqG{i}{X}$ on $X$, for each integer $i \geq 0$.
	We have that $\TTqG{1}{X}$ is the tangent space of
	$\DefqG{X}$, while $\TTqG{2}{X}$ is an obstruction space.
	If $X$ is Gorenstein, then $\DefqG{X}=\Def{X}$ and $\TTqG{i}{X}= \TT^i_X$ for every $i$.

	\subsection{Index 1 cover in the toric case}
	
	Let $N$ be a lattice and let $M = \Hom_{\ZZ}(N, \ZZ)$ be its dual lattice.
	Let $\sigma \subseteq N_\RR$ be a strongly convex rational polyhedral cone.  
	The affine toric variety $X_\sigma = \TV{\sigma}{N}$ is \emph{$\QQ$-Gorenstein} if there exist $m_0 \in M$ and $r \in \ZZ_{> 0}$ such that
	\[
	\langle m_0, v \rangle = r, 
	\]
	whenever $v$ is the primitive generator of each ray of $\sigma$.
	In other words, all ray generators of the cone $\sigma$ lie on the affine hyperplane $H_{m_0,r} = \{ n \in N_{\RR} \mid \langle m_0, n \rangle = r \}$.
	In this case, the \emph{$\QQ$-Gorenstein index} of $X_\sigma$ is the smallest positive integer $r$ for which such an element $m_0$ exists.
	From now on, $r$ will be assumed minimal. If $r = 1$, the variety $X_\sigma$ is Gorenstein.  If $r > 1$, the variety $X_\sigma$ is $\QQ$-Gorenstein of Gorenstein index $r$.
	
	Observe that $m_0$ is primitive by minimality of $r$.
	Indeed, suppose that $m_0=a m_0'$ for some $m_0'\in M$, $a \in \ZZ$ such that $a \geq 2$. Then, for every ray generator $v$ of $\sigma$, we have
	\[
	r=\langle m_0,v\rangle
	=a\langle m_0',v\rangle.
	\]
	It follows that $a$ divides $r$. Hence $r/a$ is a positive integer, and
	\[
	\langle m_0',v\rangle=\frac{r}{a}
	\]
	for every ray generator $v$ of $\sigma$.
	This contradicts the minimality of $r$.
	
	The element $m_0$ is unique when $\sigma$ is full-dimensional, i.e.\ the dimension of $\sigma$ equals the rank of $N$, whereas it is not unique otherwise.
	
	If we fix $m_0$, to construct the \emph{index 1 cover} of $X_\sigma$, we consider the sublattice $N' \subseteq N$ defined by:
	\begin{equation} \label{eq:sublattice_index1}
		N' = \{ n \in N \mid \langle m_0, n \rangle \in r\ZZ \}.
	\end{equation}
	Equivalently, $N'$ is the kernel of the homomorphism
	\[
	N \longrightarrow \mathbb Z/r\mathbb Z,
	\qquad
	n \longmapsto \langle m_0,n\rangle \bmod r.
	\]
	Since $m_0$ is primitive, this homomorphism is surjective, and hence we have an exact sequence
	\[
	0 \longrightarrow N' \longrightarrow N \longrightarrow G \longrightarrow 0,
	\]
	where $G \simeq \ZZ / r \ZZ$.
	Let $M' = \Hom_\ZZ(N', \ZZ)$ be the dual lattice of $N'$.
	The inclusion $N'\hookrightarrow N$ induces an inclusion $M\hookrightarrow M'$.
	
	The group $G=N/N'$ acts on the group algebra $\CC[M']$ as follows.
	For $n\in N$ and $u\in M'$, we set
	\[
	(n+N')\cdot \chi^u
	=
	\rme^{2\pi \rmi \langle u,n\rangle}\chi^u.
	\]
	This action is well-defined, since $\langle u,n'\rangle\in\ZZ$
	for every $n'\in N'$.
	
	Moreover, a monomial $\chi^u$, with $u\in M'$, is
	$G$-invariant if and only if
	\[
	\rme^{2\pi \rmi \langle u,n\rangle}=1
	\qquad\text{for every } n\in N,
	\]
	which is equivalent to $u\in M$. Therefore,
	\begin{equation} \label{eq:G-inv}
		\CC[M']^G=\CC[M],
		\qquad
		\CC[\sigma^\vee\cap M']^G
		=
		\CC[\sigma^\vee\cap M].
	\end{equation}
	
	We can view $\sigma$ as a cone in $N'_{\RR}$. The variety $X'= \TV{\sigma}{N'}$ is the index 1 cover of $X_\sigma$.  Indeed, if we define $u_0 = \frac{1}{r} m_0$, then $u_0 \in M'$ because for any $n \in N'$ the pairing $\langle u_0, n \rangle = \frac{1}{r} \langle m_0, n \rangle$ is an integer by definition of $N'$. Moreover, $v\in N'$ because $v$ is a generator of $\sigma$ in $N'$. Furthermore, for every ray generator $v$ of $\sigma$, we have:
	\[
	\langle u_0, v \rangle = \left\langle \frac{1}{r} m_0, v \right\rangle = \frac{1}{r} \cdot r = 1.
	\]
	So $X'$ is a Gorenstein toric variety. The inclusion $N' \hookrightarrow N$ induces a finite surjective morphism of toric varieties $\pi \colon X' \to X_\sigma$ which is a cyclic cover of degree $r$:
	by \eqref{eq:G-inv} we have
	\[
	X'/G = \Spec \left( \CC [M' \cap \sigma^\vee] \right)^G = \Spec \CC[M \cap \sigma^\vee] = X_\sigma.
	\]
	
	The qG-deformation theory of $X_\sigma$ is closely related to that of its index 1 cover $X'$. Specifically, $\TTqG{i}{X_\sigma}$ can be identified with the $G$-invariant part of $\TT^i_{X'}$. In terms of the grading by the character lattice, it holds that:
	\[
	\TTqG{1}{X_\sigma} = \bigoplus_{m \in M} \TT^1_{X'}(m).
	\]
	
	Although $m_0$ is not uniquely determined when $\sigma$ is not full-dimensional, so that the resulting sublattice $N'$ may depend on this choice, the resulting set of degrees of $\TTqG{1}{X_\sigma}$ is independent of this choice and hence well-defined. This is clarified by the following:
	
	\begin{proposition} 
		Let $N$ be a lattice of rank $n$ and $M$ be its dual lattice. Let $\sigma\subset N_\RR$ be a strongly convex rational polyhedral cone of dimension $d < n$, with primitive ray generators $v_1, \dots, v_k \in N$, and let $X_\sigma = \TV{\sigma}{N}$ be the associated affine toric variety.
		Assume that $X_\sigma$ is $\QQ$-Gorenstein of index $r$, and let $m_0, m_1 \in M$ be two primitive elements such that $\langle m_0, v_j \rangle = \langle m_1, v_j \rangle = r$ for every $j=1, \dots ,k$. Let $N'_0$ and $N'_1$ be the respective sublattices of $N$ defined as in \eqref{eq:sublattice_index1}, with associated index 1 covers $X'_0$ and $X'_1$. 
		
		Then, there exists a lattice isomorphism $\phi \colon N'_0 \xrightarrow{\sim} N'_1$ that preserves the cone $\sigma$. Consequently the following two assertions hold.
		\begin{enumerate}
			\item The isomorphism $\phi$ induces an isomorphism of toric varieties $X'_0 \simeq X'_1$.
			\item The $M$-grading of $\TTqG{1}{X_\sigma}$ is independent of the choice of the element $m_0$ used to construct the index 1 cover.
		\end{enumerate}
	\end{proposition}
	
	\begin{proof}
		Using the same notation as in \cite[\S2]{Fulton}, let $N_\sigma := \spanR{\sigma} \cap N$ be the saturated sublattice of $N$ of rank $d$ obtained by intersecting the real vector subspace generated by $\sigma$ with $N$.
		Since $N_\sigma$ is saturated in $N$, the quotient group
		$N/N_\sigma$ is torsion-free and thus a free abelian group of rank $n-d$.
		
		We first observe that $N_\sigma\cap N'_0=N_\sigma\cap N'_1$. 
		Indeed, since $\langle m_0, v_j \rangle = \langle m_1, v_j \rangle = r$ for all $j=1,\dots,k$, the linear functionals $m_0$ and $m_1$ coincide on the entire subspace $\spanR{\sigma}= (N_\sigma)_\RR$. Thus, $\langle m_0, n \rangle = \langle m_1, n \rangle$ for every $n \in N_\sigma$, which implies $N_\sigma\cap N'_0=N_\sigma\cap N'_1$. We denote this common sublattice by $L$. Moreover, since $N_\sigma$ is saturated in $N$, we have that $L$ is saturated in each $N'_i$.
		Hence for each $i=0,1$, we can choose a complement $K_i\subseteq N'_i$ such that
		\begin{equation} \label{eq:dec_sottoreticoli}
			N'_i=L \oplus K_i,
			\qquad
			K_i\simeq N'_i/L.
		\end{equation}
		Let $\{u_1, \dots, u_d\}$ be a basis of $L$. For each $i \in \{0, 1\}$, let $\{w^i_1,\ldots,w^i_{n-d}\}$ be a basis of $K_i$.
		We consider the lattice isomorphism
		\begin{equation*}
			\phi \colon N'_0 \longrightarrow N'_1 ,
		\end{equation*}
		which sends each basis element $u_s$ to $u_s$ for $s=1, \dots, d$, and each $w^0_t$ to $w^1_t$ for $t=1, \dots, n-d$. By construction, $\phi$ is an isomorphism of lattices whose restriction to $L$ is the identity. Since $\sigma \subset \spanR{L}$, the extension $\phi_\RR \colon (N'_0)_\RR \to (N'_1)_\RR$ of $\phi$ satisfies $\phi_\RR(\sigma) = \sigma$. This induces a toric isomorphism $X'_0 \simeq X'_1$ and proves (1).
		
		Let
		\[
		L^*:=\Hom_\ZZ (L,\ZZ),
		\qquad
		K^*_i:=\Hom_\ZZ (K_i,\ZZ)
		\]
		be the dual lattices. Dualising the decomposition \eqref{eq:dec_sottoreticoli}, we
		obtain
		\[
		M'_i= \Hom_\ZZ(N'_i,\ZZ) \simeq L^* \oplus K^*_i.
		\]
		Since $\sigma$ is contained in $\spanR{L}$, we have
		\[
		\sigma^\vee\cap M'_i
		=
		(\sigma^\vee\cap L^*)\oplus K^*_i.
		\]
		Setting $Y:=\TV{\sigma}{L}$ and
		$T_{K_i}:=\Spec \CC [K^*_i] \simeq (\CC^*)^{n-d}$, we obtain
		\[
		X'_i\simeq Y\times T_{K_i}.
		\]
		Although the decomposition \eqref{eq:dec_sottoreticoli} depends on the choice of $K_i$, the variety $Y$ is the same for both index-one covers since $L$ does not change.
		This product decomposition induces an isomorphism
		\begin{equation} \label{eq:dec_T1}
			\TT^1_{X'_i}
			\simeq
			\TT^1_Y\otimes_\CC \CC[K^*_i],
		\end{equation}
		where in particular
		\[
		\mathbb C[K^*_i]=\bigoplus_{k\in K^*_i}\mathbb C\chi^k
		\]
		is graded by $K^*_i$ and
		\[
		\TT^1_Y= \bigoplus_{\ell \in L^*} \TT^1_Y (\ell)
		\]
		is graded by $L^*$.
		Now let $m\in M$. Under the inclusion $M\hookrightarrow M'_i$,
		the character $m$ corresponds to
		\[
		\bigl(m|_{L},m|_{K_i}\bigr)
		\in L^* \oplus K^*_i.
		\]
		Therefore, by taking the eigenspace corresponding to $m$ in \eqref{eq:dec_T1}, we obtain
		\[
		\TT^1_{X'_i}(m)
		\simeq
		\TT^1_Y(m|_{L})
		\otimes_\CC \CC \chi^{m|_{K_i}}.
		\]
		Since the second factor is 1-dimensional, it follows that
		\[
		\dim_\CC \TT^1_{X'_i}(m)
		=
		\dim_\CC \TT^1_Y(m|_{L}).
		\]
		which is independent of $i$, since both
		$L$ and $Y$ are independent of the choice of $m_0$.
		Thus
		\[
		\dim_\CC\TT^1_{X'_0}(m)= \dim_\CC \TT^1_{X'_1}(m)
		\]
		for every $m\in M$.
		Finally, since $\TTqG{1}{X_\sigma}$ is the $\mu_r$-invariant part of $\TT^1_{X'_i}$, we have
		\[
		\TTqG{1}{X_\sigma}
		\simeq
		\bigoplus_{m\in M}\TT^1_{X'_i}(m).
		\]
		Consequently, the dimensions of
		$\TTqG{1}{X_\sigma}(m)$ are independent of the choice of
		$m_0$. This proves (2).
	\end{proof}
	
	\begin{remark}
		Although $m_0\in M$ is not unique when
		$\sigma$ is not full-dimensional, its restriction to
		$N_\sigma :=\spanR{\sigma} \cap N$ is uniquely determined.
		
		Indeed, let $m_0,m_1\in M$ be two elements satisfying
		\[
		\langle m_0,v_j\rangle
		=
		\langle m_1,v_j\rangle
		=
		r
		\]
		for every primitive ray generator $v_j$ of $\sigma$.
		It follows that the linear functionals $m_0$ and $m_1$ agree on this entire subspace $\spanR{\sigma}$
		In particular,
		\[
		m_0|_{N_\sigma}=m_1|_{N_\sigma}.
		\]
	\end{remark}

	\subsection{K-stability}
	
	As briefly discussed in the introduction, we will use the theory of K-moduli of Fano varieties. For a general overview of this subject, we refer the reader to \cite{Xu_survey} and \cite{xu_book}, and references therein. 
	
	There are three main concepts, each strictly weaker than the previous one: K-stability, K-polystability, and K-semistability.
	Moreover, a K-semistable Fano variety has at worst klt singularities. 
	Furthermore, a Fano variety is K-stable if and only if it is K-polystable and has a finite automorphism group.
	
	In general, deciding whether a given Fano variety is K-polystable is hard.
	In the smooth case, in dimension two we refer to \cite{Tian} and in dimension three we refer to \cite{CalabiProb}.
	
	In the case of toric Fano varieties, this problem becomes easy. In fact, thanks to \cite[Corollary~7.17]{BJ20}, a toric Fano variety is K-polystable if and only if the barycentre of its polar polytope is the origin. We will use this explicit combinatorial criterion to study the K-polystability of the toric variety $X$ constructed in \S\ref{sec:toric_var_X}.

	\section{An affine toric variety} \label{sec:affine_toric_var}
	
	We now study an affine toric variety associated with a cone over an octagon. This variety will serve as one of the affine charts of the toric variety defined in \S\ref{sec:toric_var_X}.
	
	\begin{figure}
		\centering
		\includegraphics[width=0.5\textwidth]{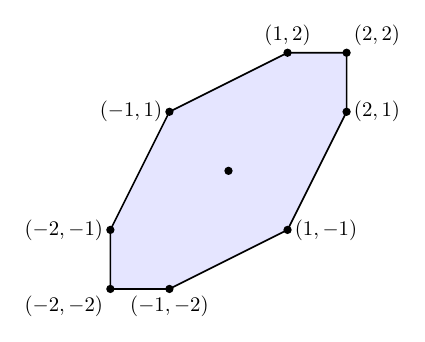}
		\caption{The lattice polygon $F$ appearing in Proposition~\ref{prop:torvar_cone_over_octagon}}
		\label{fig:octagon}
	\end{figure}
	
\begin{figure}
		\begin{center}
			\raisebox{-4mm}{\includegraphics[height=28mm]{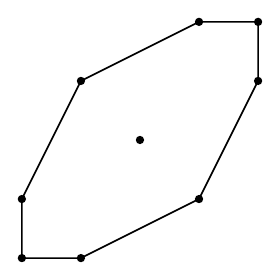}} \quad
			\raisebox{5mm}{=} \quad
			\raisebox{0mm}{\includegraphics[height=14mm]{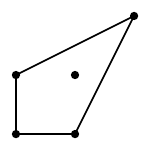}} \quad
			\raisebox{5mm}{+} \quad
			\raisebox{0mm}{\includegraphics[height=14mm]{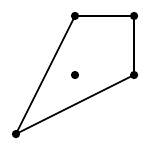}} \quad
			\raisebox{5mm}{=} \quad \\ \raisebox{2mm}{=} \quad
			\raisebox{0mm}{\includegraphics[height=7mm]{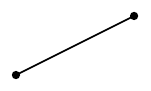}} \quad
			\raisebox{2mm}{+} \quad
			\raisebox{1,5mm}{\includegraphics[height=3mm]{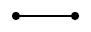}} \quad
			\raisebox{2mm}{+} \quad
			\raisebox{-0,5mm}{\includegraphics[height=7mm]{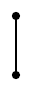}} \quad
			\raisebox{2mm}{+} \quad
			\raisebox{-3mm}{\includegraphics[height=14mm]{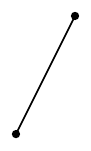}}
		\end{center}
		
		\caption{ The Minkowski decompositions of the lattice polygon $F$ defined in Proposition \ref{prop:torvar_cone_over_octagon}.}
		\label{fig:decMinkF}
	\end{figure}

	\begin{proposition} \label{prop:torvar_cone_over_octagon}
		Let $F$ be the lattice polygon with the following vertices in $\ZZ^2$ (see Figure~\ref{fig:octagon}):
		\[
		\begin{aligned}
			v_{1}=\begin{pmatrix}2\\2 \end{pmatrix}, \quad
			v_{2}=\begin{pmatrix}1\\2 \end{pmatrix}, \quad
			v_{3}=\begin{pmatrix}-1\\1 \end{pmatrix}, \quad
			v_{4}=\begin{pmatrix}-2\\-1 \end{pmatrix}, \\[8pt]
			v_{5}=\begin{pmatrix}-2\\-2 \end{pmatrix}, \quad
			v_{6}=\begin{pmatrix}-1\\-2 \end{pmatrix}, \quad
			v_{7}=\begin{pmatrix}1\\-1 \end{pmatrix}, \quad
			v_{8}=\begin{pmatrix}2\\1 \end{pmatrix}.
		\end{aligned}
		\]
		Let $V$ be the affine toric variety over $\CC$ associated to the cone $\RR_{\geq0}(F \times \{1\})$ in the lattice $\ZZ^3 = \ZZ^2 \oplus \ZZ$.
		
		Then the following statements hold.
		\begin{enumerate}[label=(\roman*)]
			\item The hull of $\Def{V}$ is
			\[ R_F\coloneqq \CC \pow{x_1, \dots, x_5 }/I_F,
			\]
			where $I_F$ is the ideal generated by the polynomials
			\begin{align*}
				f_k   &= p^k - x_1^k - 2x_2^k - x_3^k + 2x_5^k, \\
				g_k   &= q^k - x_2^k - 2x_3^k - x_4^k + x_5^k,
			\end{align*}
			for $k=2, \dots, 6$, with
			\begin{align*}
				p   &= x_1 + 2x_2 + x_3 - 2x_5, \\
				q   &= x_2 + 2x_3 + x_4 - x_5.
			\end{align*}
			\item The embedding dimension of $R_F$ is 5.
			\item The base of the miniversal deformation of $V$ consists of two irreducible components: 
			\begin{itemize}
				\item one has dimension 1, is everywhere non-reduced and corresponds to the minimal prime
				$\frakp_1=(x_1 - x_5, x_2, x_3 - x_5, x_4)$ and to the Minkowski decomposition of $F$ into two quadrilaterals as shown in Figure~\ref{fig:decMinkF};
				\item the other has dimension 3, is generically reduced and corresponds to the minimal prime
				$\frakp_2=(x_2 - x_5, x_3)$  and to the Minkowski decomposition of $F$ into four segments as shown in Figure~\ref{fig:decMinkF}.
			\end{itemize}
			\item If the base field $\CC$ is the field of complex numbers, the fibre over a general point of the 1-dimensional component of the base of the analytic miniversal deformation of $V$ has two isolated singularities, which are analytically isomorphic to the anticanonical affine cone over the first Hirzebruch surface. The fibre over a general point of the 3-dimensional component of the base of the analytic miniversal deformation of $V$ is smooth. 
			\item With respect to the grading induced by the torus action, $\TT_V^1$ is concentrated in degree $(0,0,-1)$.
		\end{enumerate}
	\end{proposition}
	
	\begin{proof}
		Following the construction described in \cite{Alt97} (see also \cite[\S4.1]{Pet22a}), we construct the algebra $A_F$ associated to the polygon $F$. 
		Fix the counterclockwise orientation of $F$, starting with the vertex $v_8$, and order the edges accordingly. The vectors associated to the first seven edges are given by the columns of the following matrix:
		\[
		\begin{array}{ccccccc}
			q & x_1 & x_2 & x_3 & x_4 & p & x_5 \\ \hline
			0 & -1 & -2 & -1 & 0 & 1 & 2 \\
			1 & 0 & -1 &  -2 & -1 & 0 & 1
		\end{array}
		\]
		We thus obtain that the algebra $A_F$ is isomorphic to $\CC[x_1,x_2, x_3, x_4, x_5]/J_F$, where 
		\[
		p=x_1 +2x_2+x_3-2x_5, \qquad
		q=x_2+2x_3+x_4-x_5,
		\]
		and $J_F$ is the ideal generated by
		\[
		f_k= p^k -x_1^k -2x_2^k -x_3^k +2x_5^k, \qquad
		g_k= q^k -x_2^k -2x_3^k -x_4^k +x_5^k,
		\]
		for $k=2, \dots ,6$.
		
		The items (i) and (ii) of the proposition are a direct consequence of Proposition~\ref{prop:Alt_hull}.
		
		To complete the proof, we need to examine the primary decomposition of the ideal $J_F$ and to study the minimal prime ideals of the algebra $A_F$, that are in one-to-one correspondence with the minimal prime ideals of $R_F$ by \cite[Lemma~3.2]{Pet22a}. Using the Magma algebra system \cite{Magma}, we get:
		\begin{equation*}
			J_F= \bigcap_{i=1}^7 \frakq_i ,
		\end{equation*}
		where
		\begin{align*}
			\frakq_1 & =
			\Big(
			x_1 + x_3 - 2x_5,\;
			x_2 + x_4,\;
			x_3^2 - 2x_3x_5 + x_5^2,\;
			x_4^2
			\Big),
		\end{align*}
		\begin{align*}
			\frakq_2 & = 
			\Big( x_2 - x_5,\;
			x_3
			\Big), 
		\end{align*}
		\begin{align*}
			\frakq_3 & =
			\Big(
			x_1^2 - 2x_1x_5 + x_5^2,\;
			x_1x_2 - x_1x_5 + \tfrac{1}{2}x_2^2 - 2x_2x_5 + \tfrac{3}{2}x_5^2, \\
			& \phantom{{}=\frakq_3} \, x_2^3 - 3x_2^2x_5 + 3x_2x_5^2 - x_5^3,\;
			x_3,\;
			x_4 - x_5
			\Big),
		\end{align*}
		\begin{align*}
			\frakq_4 & =
			\Big(
			x_1^2,\;
			x_1x_3 - \tfrac{4}{3}x_3x_5,\;
			x_2 + 2x_3 - x_5,\;
			x_3^2 - \tfrac{2}{3}x_3x_5,\;
			x_5^2
			\Big),
		\end{align*}
		\begin{align*}
			\frakq_5 & =
			\Big(
			x_1,\;
			x_2 - x_5,\;
			x_3^2 + 2x_3x_4,\;
			x_4^2
			\Big),
		\end{align*}
		\begin{align*}
			\frakq_6 & =
			\Big(
			x_2 + \tfrac{1}{2}x_3 - x_5,\;
			x_3^2 - \tfrac{4}{3}x_3x_5,\;
			x_3x_4 + \tfrac{1}{3}x_3x_5,\;
			x_4^4,\;
			x_5^2
			\Big),
		\end{align*}
		
		\begin{align*}
			\frakq_7 & =
			\Big( x_1^4, x_1x_2 + \tfrac{1}{2}x_1x_3 - x_1x_5
			+ \tfrac{1}{2}x_2^2 - \tfrac{1}{2}x_2x_4 - \tfrac{3}{2}x_2x_5
			- \tfrac{1}{2}x_3^2 - x_3x_4 + \tfrac{1}{2}x_4x_5, \\
			&\, x_1x_3^2 - \tfrac{4}{3}x_1x_3x_5
			+ \tfrac{2}{3}x_2^3 - \tfrac{1}{2}x_2^2x_4 - \tfrac{3}{2}x_2^2x_5
			+ \tfrac{4}{3}x_2x_4x_5
			- \tfrac{4}{3}x_3^3 - \tfrac{3}{2}x_3^2x_4 + \tfrac{3}{2}x_3^2x_5 + \notag\\[0.3em]
			&+ 3x_3x_4x_5,
			\, x_1x_3x_4 + \tfrac{1}{3}x_1x_3x_5
			+ \tfrac{2}{3}x_3^3 + x_3^2x_4 - x_3^2x_5 - 2x_3x_4x_5,
			\, x_2^4 - \tfrac{8}{3}x_2^3x_5, \notag\\[0.3em]
			&\, x_2^3x_4 + \tfrac{1}{3}x_2^3x_5 - 3x_2^2x_4x_5
			- \tfrac{4}{3}x_3^3x_5 + 3x_3^2x_4^2- 2x_3^2x_4x_5 - 2x_3x_4^2x_5, \notag\\[0.3em]
			&\, x_2x_3 + \tfrac{1}{2}x_2x_4 - \tfrac{1}{2}x_2x_5
			+ \tfrac{1}{2}x_3^2 + x_3x_4 - x_3x_5 - \tfrac{1}{2}x_4x_5, \; x_3^4 - \tfrac{4}{3}x_3^3x_5, \notag\\[0.3em]
			&\, x_2x_4^2 - \tfrac{2}{3}x_2x_4x_5
			+ \tfrac{2}{3}x_3^3 + x_3^2x_4 - x_3^2x_5
			+ 2x_3x_4^2 - 2x_3x_4x_5 - x_4^2x_5, \notag\\[0.3em]
			&\, x_3^3x_4 - \tfrac{1}{3}x_3^3x_5
			+ \tfrac{3}{2}x_3^2x_4^2 - 2x_3^2x_4x_5 - x_3x_4^2x_5,\; x_3^2x_4^3 - x_3^2x_4^2x_5 - \tfrac{2}{3}x_3x_4^3x_5,\; x_4^4, \; x_5^2
			\Big),
		\end{align*}
		with associated prime ideals $\frakp_i= \sqrt{\frakq_i}$.
		
		There are only two minimal primes: 
		\begin{align*}
			\frakp_1=(x_1 - x_5, x_2, x_3 - x_5, x_4), && \frakp_2=(x_2 - x_5, x_3).
		\end{align*}
		The ideal $\frakp_1$ corresponds to the $1$-dimensional irreducible component, while $\frakp_2$ corresponds to the $3$-dimensional irreducible component. In particular, they correspond to the two maximal Minkowski decompositions of $F$ shown in Figure \ref{fig:decMinkF} (see Proposition~\ref{prop:Alt_hull}).
		
		We examine the structure of $\Def{V}$ at its generic points. We shall analyse the localisations of the algebra $A_F$ at its two minimal prime ideals.
		
		Let $S:= \CC[x_1,x_2, x_3, x_4, x_5]$, then
		\begin{equation*}
			(A_F)_{\frakp_1} = (S/J_F)_{\frakp_1}= S_{\frakp_1} / J_F S_{\frakp_1}.
		\end{equation*}
		The only ideals $\frakq_i$ whose localisation at $\frakp_1$ is not the whole ring $S_{\frakp_1}$ are those for which $\frakp_i \subseteq \frakp_1$. As $\frakp_1$ is minimal, it follows that:
		\begin{equation*}
			J_F S_{\frakp_1} = \frakq_1 S_{\frakp_1},
		\end{equation*}
		where $\frakq_1= (x_1+x_3-2x_5, x_2+x_4, (x_3-x_5)^2, x_4^2)$.
		
		Let us now make a change of coordinates: 
		\[
		\begin{cases}
			u_1 = x_1 - x_5, \\
			u_2 = x_2, \\
			u_3 = x_3 - x_5, \\
			u_4 = x_4, \\
			u_5 = x_5.
		\end{cases}
		\]
		In the new coordinate system the ideals may be expressed as:
		\begin{align*}
			\frakp_1= (u_1, u_2, u_3, u_4), &&
			\frakq_1= (u_1+u_3, u_2+u_4, u_3^2, u_4^2).
		\end{align*}
		It follows that
		\begin{align*}
			(A_F)_{\frakp_1}
			&= \frac{S_{\frakp_1}}{J_F S_{\frakp_1}}
			= \frac{S_{\frakp_1}}{\frakq_1 S_{\frakp_1}} \\[4pt]
			&= \frac{\CC(u_5)[u_1,u_2,u_3,u_4]_{(u_1,u_2,u_3,u_4)}}
			{(u_1+u_3,\; u_2+u_4,\; u_3^2,\; u_4^2)} \\[6pt]
			&\simeq \frac{\CC(u_5)[u_3,u_4]_{(u_3,u_4)}}
			{(u_3^2,\; u_4^2)} \\[6pt]
			&\simeq \frac{\CC(u_5)[u_3,u_4]}
			{(u_3^2,\; u_4^2)} .
		\end{align*}
		The last isomorphism is justified by the fact that the localisation at $(u_3,u_4)$ does not change the quotient.
		Indeed, set
		\[
		A=\CC(u_5)[u_3,u_4], \qquad
		H=(u_3^2,u_4^2), \qquad
		M=A/H,
		\]
		and let $\frakm=(u_3,u_4)\in \mathrm{Specm}(A)$.
		We apply the following general fact: if $U\subseteq T$ are multiplicative
		subsets of $A$ such that for every $\frakp \in \Spec(A)$ one has
		\[
		\frakp\cap U=\varnothing,\quad
		\frakp\supseteq \mathrm{ann}_A(M)
		\ \Longrightarrow\
		\frakp\cap T=\varnothing,
		\]
		then the natural localisation morphism $U^{-1}M\to T^{-1}M$ is bijective.
		In our situation we choose $U=\{1\}$ and $T=A\setminus \frakm$.
		Since $\sqrt{\mathrm{ann}_A(M)}=\sqrt{H}=(u_3,u_4)= \frakm$, the only prime ideal
		containing $\mathrm{ann}_A(M)$ is $\frakm$ itself, and clearly
		$\frakm \cap (A\setminus \frakm)=\varnothing$.
		Therefore the hypotheses of the general fact are satisfied, and the localisation map
		$M \longrightarrow M_\frakm$
		is an isomorphism. 
		
		Similarly, for the prime ideal $\frakp_2= \frakq_2=(x_2 - x_5, x_3)$, via the following change of coordinates
		\[
		\begin{cases}
			w_1 = x_1, \\
			w_2 = x_2-x_5, \\
			w_3 = x_3, \\
			w_4 = x_4, \\
			w_5 = x_5,
		\end{cases}
		\]
		we obtain:
		\begin{equation*}
			(A_F)_{\frakp_2} = \frac{S_{\frakp_2}}{J_F S_{\frakp_2}}
			= \frac{S_{\frakp_2}}{\frakq_2 S_{\frakp_2}} = \CC(w_1, w_4, w_5).
		\end{equation*}
		
		These calculations imply that the $3$-dimensional component is reduced and smooth at its general point. Conversely, the presence of nilpotent elements in $(A_F)_{\frakp_1}$ shows that the $1$-dimensional component is everywhere non-reduced.
		
		To conclude the proof of item~(iv) about the geometry of the fibres, we apply \cite[Proposition~2.5]{CHP_smoothing} that relates the singularities of the general fibre to the Minkowski decomposition of the polygon:
		\begin{itemize}
			\item the $1$-dimensional component corresponds to the prime ideal $\frakp_1$ and to the Minkowski decomposition $F= Q_1 + Q_2$ into two quadrilaterals (see Figure~\ref{fig:decMinkF}). Each of these, following \cite{AltMin}, corresponds to the anticanonical affine cone $C_a(\FF_1)$ over the first Hirzebruch surface, which is an isolated Gorenstein singularity. Using \cite[Proposition~2.5]{CHP_smoothing}, the singular locus of the general fibre over this component consists of two isolated points analitically isomorphic to $C_a(\FF_1)$.
			
			\item The Minkowski decomposition corresponding to the $3$-dimensional component is given by the sum of four primitive lattice segments (see again Figure~\ref{fig:decMinkF}). Consequently, the affine toric varieties associated with these segments are smooth, and hence the fibre over the generic point of this component is smooth.
		\end{itemize}
		
		Finally, by \cite[\S 6.3]{Alt97}, $\TT_V^1$ is concentrated in degree $(0,0,-1)$.
		
		This concludes the proof of Proposition~\ref{prop:torvar_cone_over_octagon}.
	\end{proof}

	\section{A toric Fano variety} \label{sec:toric_var_X}
	
	In this section we construct and study a toric Fano variety called $X$, that will later be used in the proof of Theorem \ref{thm:result_n=3}.
	
	\begin{figure}
		\centering
		\includegraphics[width=0.5\textwidth]{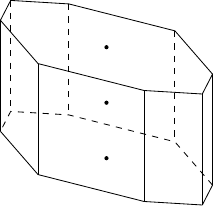}
		\caption{The octagonal prism defined in Theorem~\ref{thm:tor_var_X}.}
		\label{fig:prism}
	\end{figure}
	
	\begin{theorem} \label{thm:tor_var_X}
		Let $P$ be the polytope in $N=\ZZ^3$ whose sixteen vertices are:
		\[
		\begin{aligned}
			v_{1, \pm}&=\begin{pmatrix}2\\2\\ \pm1\end{pmatrix}, \quad
			v_{2, \pm}=\begin{pmatrix}1\\2\\ \pm1\end{pmatrix}, \quad
			v_{3, \pm}=\begin{pmatrix}-1\\1\\ \pm1\end{pmatrix}, \quad
			v_{4, \pm}=\begin{pmatrix}-2\\-1\\ \pm1\end{pmatrix}, \\[8pt]
			v_{5, \pm}&=\begin{pmatrix}-2\\-2\\ \pm1\end{pmatrix}, \quad
			v_{6, \pm}=\begin{pmatrix}-1\\-2\\ \pm1\end{pmatrix}, \quad
			v_{7, \pm}=\begin{pmatrix}1\\-1\\ \pm1\end{pmatrix}, \quad
			v_{8, \pm}=\begin{pmatrix}2\\1\\ \pm1\end{pmatrix}.
		\end{aligned}
		\]
		It is the octagonal prism shown in Figure~\ref{fig:prism}. 
		Let $X$ be the toric 3-fold associated to the face fan of $P$. Then the following statements hold.
		\begin{enumerate}[label=(\roman*)]
			\item $X$ is a K-polystable Fano 3-fold with anticanonical volume 11/3.
			
			\item The automorphism group $\Aut(P)$ of the polytope $P$ is isomorphic to $(\ZZ/2\ZZ)^3$, generated by the following three involutions: 
			\[M_1=
			\begin{pmatrix}
				1 & 0 & 0 \\
				0 & 1 & 0 \\
				0 & 0 & -1
			\end{pmatrix}
			, \quad M_2=
			\begin{pmatrix}
				-1 & 0 & 0 \\
				0 & -1 & 0 \\
				0 & 0 & 1
			\end{pmatrix}
			, \quad M_3=
			\begin{pmatrix}
				0 & -1 & 0 \\
				-1 & 0 & 0 \\
				0 & 0 & 1
			\end{pmatrix}.
			\]
			
			\item The automorphism group $\Aut(X)$ of $X$ is isomorphic to $T_N\rtimes \Aut(P)$, where $T_N = \Gm^3 = (\CC^*)^3$ is the 3-dimensional algebraic torus acting on the toric variety $X$.
			
			\item The base space of the miniversal qG-deformation of $X$ is the formal spectrum of the local noetherian $\CC$-algebra
			\[ R:= \CC \pow{y_1, \dots, y_6, x_{1,+}, \dots x_{5,+}, x_{1,-}, \dots, x_{5,-}} /I,
			\]
			where $I$ is the ideal generated by the following polynomials:
			\begin{align*}
				f_{\varepsilon,k}   &= p_\varepsilon^k - x_{1,\varepsilon}^k - 2x_{2,\varepsilon}^k - x_{3,\varepsilon}^k + 2x_{5,\varepsilon}^k, \\
				g_{\varepsilon,k}   &= q_\varepsilon^k - x_{2,\varepsilon}^k - 2x_{3,\varepsilon}^k - x_{4,\varepsilon}^k + x_{5,\varepsilon}^k, 
			\end{align*}
			for $k=2, \dots ,6$ and $\varepsilon \in \{+,-\}$, with
			\begin{align*}
				p_\varepsilon   &= x_{1,\varepsilon} + 2x_{2,\varepsilon} + x_{3,\varepsilon} - 2x_{5,\varepsilon}, \\
				q_\varepsilon   &= x_{2,\varepsilon} + 2x_{3,\varepsilon} + x_{4,\varepsilon} - x_{5,\varepsilon}.
			\end{align*}
			
			\item 
			The base space of the miniversal qG-deformation of $X$ has four irreducible components: three are everywhere non-reduced (one of dimension 8 and two of dimension 10), while the fourth is generically reduced (of dimension 12).
			These components correspond to the following minimal prime ideals of $R$:
			\[
			\mathfrak{P}_{1,1}
			= (x_{1,+} - x_{5,+},\, x_{2,+},\, x_{3,+} - x_{5,+},\, x_{4,+},\,
			x_{1,-} - x_{5,-},\, x_{2,-},\, x_{3,-} - x_{5,-},\, x_{4,-}),
			\]
			\[
			\mathfrak{P}_{1,2}
			= (x_{1,+} - x_{5,+},\, x_{2,+},\, x_{3,+} - x_{5,+},\, x_{4,+},\,
			x_{2,-} - x_{5,-},\, x_{3,-}),
			\]
			\[
			\mathfrak{P}_{2,1}
			= (x_{2,+} - x_{5,+},\, x_{3,+},\,
			x_{1,-} - x_{5,-},\, x_{2,-},\, x_{3,-} - x_{5,-},\, x_{4,-}),
			\]
			\[
			\mathfrak{P}_{2,2}
			= (x_{2,+} - x_{5,+},\, x_{3,+},\, x_{2,-} - x_{5,-},\, x_{3,-}).
			\]
			\begin{align*}
				\dim(R/\mathfrak{P}_{1,1})=8,\qquad
				\dim(R/\mathfrak{P}_{1,2})=10,\\
				\dim(R/\mathfrak{P}_{2,1})=10,\qquad
				\dim(R/\mathfrak{P}_{2,2})=12.
			\end{align*}
			
			\item The automorphism group $\Aut(X)=T_N\rtimes \Aut(P)$ acts on $\DefqG{X} \simeq \Spf R$ as follows.
			\begin{itemize}
				\item The torus $T_N$ acts linearly and diagonally with weights 
				\[
				(0,-1,0)\textrm{ on } y_1, \quad 
				(0,1,0) \textrm{ on } y_2, 
				\]
				\[
				(-1,0,0)\textrm{ on } y_3, \quad 
				(1,0,0)\textrm{ on } y_4,
				\]
				\[
				(-1,1,0)\textrm{ on } y_5, \quad 
				(1,-1,0)\textrm{ on } y_6,
				\]
				\[
				(0,0,-1) \textrm{ on } x_{1,-}, \ldots, x_{5,-} , \quad
				(0,0,1) \textrm{ on } x_{1,+}, \ldots, x_{5,+}.
				\]
				
				\item The involution $M_1$ fixes the variables $y_j$ for $j=1,\dots,6$, and swaps the variables $x_{i,+}$ and $x_{i,-}$ for each $i=1,\dots,5$.
				
				\item The involution $M_2$ acts on the variables $y_1,\ldots,y_6$ as follows:
				\[
				M_2 \cdot y_1 = y_2, \quad M_2 \cdot y_3 = y_4, \quad M_2 \cdot y_5 = y_6,
				\]
				\item The involution $M_3$ acts on the variables $y_1,\ldots,y_6$ as follows:
				\[
				M_3 \cdot y_1 = y_4, \quad M_3 \cdot y_2 = y_3, \quad M_3 \cdot y_5 = y_5, \quad M_3 \cdot y_6 = y_6.
				\]
			\end{itemize}
		\end{enumerate}
	\end{theorem}
	
	\begin{remark} \label{rem:action}
		We do not explicitly describe the action of $M_2$ and $M_3$ on the variables $x_{i,\pm}$, since this will not be needed. Nevertheless, we can deduce that these involutions preserve the two sets of variables $\{x_{i,+}\}_{i=1, \dots, 5}$ and $\{x_{i,-}\}_{i=1, \dots, 5}$. Indeed, the torus $T_N$ acts on these two sets of variables with weights $(0,0,-1)$ and $(0,0,1)$, respectively, and both $M_2$ and $M_3$ preserve the third coordinate of the weights. Hence, they preserve the corresponding weight spaces and do not mix the two sets.
		
		It follows that the action of the whole group $\Aut (X)=T_N\rtimes\Aut(P)$ preserves the two sets of variables $\{ x_{i,\pm}\}_{i=1, \dots, 5}$ and $\{ y_{1}, \dots, y_{6}\}$. 
		If we write the hull $R$ as
		\[
		R\simeq R_\mathbf{y} \widehat{\otimes}_{\CC}R_\mathbf{x},
		\]
		where
		\[
		R_\mathbf{y}:=\CC\pow{y_1,\ldots,y_6}
		\]
		and
		\[
		R_\mathbf{x}:=
		\frac{
			\CC\pow{x_{1,+},\ldots,x_{5,+},x_{1,-},\ldots,x_{5,-}}}{I_\mathbf{x}},
		\]
		where $I_\mathbf{x}$ denotes the ideal generated by the polynomials $f_{\varepsilon,k}$, $g_{\varepsilon,k}$, we have that both $R_\mathbf{y}$ and $R_\mathbf{x}$ inherit an action of $\Aut(X)$.
	\end{remark}
	
	\begin{proof}[Proof of items~(i)-(iii) of Theorem\ref{thm:tor_var_X}]
		$P$ is a Fano 3-dimensional polytope, hence the associated toric variety is a toric Fano 3-fold. 
		The polar polytope $P^\circ$ has the following vertices:
		\[
		\pm
		\begin{pmatrix} 0 \\ \tfrac{1}{2} \\ 0 \end{pmatrix},
		\;
		\pm
		\begin{pmatrix} 0 \\ 0 \\ 1 \end{pmatrix},
		\;
		\pm
		\begin{pmatrix} \tfrac{1}{3} \\ -\tfrac{2}{3} \\ 0 \end{pmatrix},
		\;
		\pm
		\begin{pmatrix} \tfrac{1}{2} \\ 0 \\ 0 \end{pmatrix},
		\;
		\pm
		\begin{pmatrix} \tfrac{2}{3} \\ -\tfrac{1}{3} \\ 0 \end{pmatrix}.
		\]
		Since $P$ is centrally symmetric, also $P^\circ$ is centrally symmetric. Hence the barycentre of $P^\circ$ is the origin and we can conclude that $X$ is  K-polystable by \cite[Corollary~7.17]{BJ20}. Using Magma algebra system \cite{Magma}, we compute the anticanonical volume, which is 11/3. So item~(i) of the theorem is proved. 
		
		Let $\Aut(P) \subset \GL_3(\ZZ)$ be the automorphism group of the polytope $P$.
		It is easy to check that
		$\Aut(P)$ is generated by the involution
		\[ M_1=
		\begin{pmatrix}
			1 & 0 & 0 \\
			0 & 1 & 0 \\
			0 & 0 & -1
		\end{pmatrix},
		\]
		which swaps the top facet of $P$ and the bottom facet of $P$, 
		and by the involutions
		\[M_2=
		\begin{pmatrix}
			-1 & 0 & 0 \\
			0 & -1 & 0 \\
			0 & 0 & 1
		\end{pmatrix}
		, \quad  M_3=
		\begin{pmatrix}
			0 & -1 & 0 \\
			-1 & 0 & 0 \\
			0 & 0 & 1
		\end{pmatrix},
		\]
		that preserve the third coordinate and act on the octagon.
		$M_2$ is the reflection through the origin in the plane of the first two
		coordinates, while $M_3$ is the reflection with respect to the bisector of the
		second and fourth quadrants in the plane of the first two coordinates.
		Hence, $\Aut(P)$ has order 8 and is isomorphic to $(\ZZ/2\ZZ)^3$.
		
		By \cite[Proposition~2.8]{KP} (see also \cite{Nill_2006, Cox95, WG99}), since no facet of $P^\circ$ has interior lattice points, $\Aut(X)$ is the semidirect product of 
		$T_N = N \otimes_{\ZZ} \CC^* = (\CC^*)^3$ with $\Aut(P)$.	
		This proves items~(ii) and (iii) of Theorem~\ref{thm:tor_var_X}.	
	\end{proof}
	
	\medskip
	
	Before continuing with the proof, we first fix the notation for the cones of the face fan of the polytope $P$ of Theorem \ref{thm:tor_var_X} and the associated affine charts of the toric variety $X$.
	
	\subsection{Notation for cones and affine charts}\label{sec:notation_cones_charts}
	Let
	\[
	\sigma_{+}=\cone{v_{1,+},\ldots,v_{8,+}}, 
	\quad
	\sigma_{-}=\cone{v_{1,-},\ldots,v_{8,-}} \subseteq N_{\RR}
	\]
	be the 3-dimensional cones over the two octagonal facets of the polytope $P$ defined in Theorem~\ref{thm:tor_var_X}, and denote by $U_{+}, U_{-}$ the corresponding affine charts of $X$.
	
	For each $i \in \ZZ/8\ZZ$, consider the cone
	\[
	\sigma_{i}=\cone{v_{i, \pm},v_{i+1, \pm}}\subseteq N_{\RR}.
	\]
	It is the cone over a rectangular facet of $P$. Denote by $U_i$ the corresponding affine chart of $X$.
	We have obtained an open affine cover of $X$:
	\[
	X = U_{+} \cup U_{-}  \cup \bigcup_{i \in \ZZ/8\ZZ} U_i.
	\]
	
	Set $U_{i,j}=U_i\cap U_j$ and
	\[
	U=\bigcup_{i \in \ZZ/8\ZZ} U_i.
	\]
	
	Let
	\[
	\tau_i= \cone{v_{i,\pm}} = \sigma_i \cap \sigma_{i-1}
	\]
	be the 2-dimensional cone corresponding to $U_{i,i-1}$.
	
\begin{figure}
		\centering
		\includegraphics[width=0.5\textwidth]{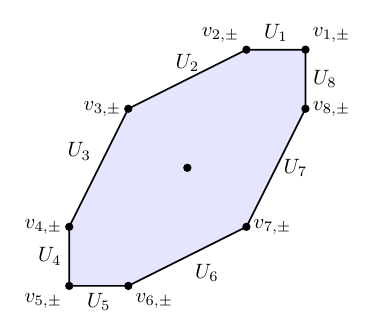}
		\caption{A top view of the charts defined in \S\ref{sec:notation_cones_charts} of the toric variety $X$ of Theorem~\ref{thm:tor_var_X}, together with the vertices of the polytope.}
		\label{figure:charts}
	\end{figure}

	See Figure~\ref{figure:charts} for a top view of the charts of $U$, together with the vertices of the polytope.
	
	Before continuing with the proof of the Theorem~\ref{thm:tor_var_X}, we first study the deformations of $U$.
	
	\subsection{Deformations of $U$.} \label{sec:defU}
	
	In the following lemmas, we will study the various charts $U_i$ that cover $U$, as well as their pairwise intersections $U_{i,j}$.
	
	\begin{lemma}  \label{lem:U18}
		The affine chart $U_{1,8}$ is isomorphic to 
		$
		\Gm \times \frac{1}{4}(1,1)
		$ 
		and is $\QQ$-Gorenstein with Gorenstein index 2. Its index 1 cover is isomorphic to
		\[
		\Gm \times \Spec  \frac{\CC[z_1, z_2, z_3]}{(z_1 z_3 - z_2^2)} .
		\]
		Furthermore,
		\begin{equation} \label{eq:T1_U18}
			\dim \TTqG{1}{U_{1,8}}(m) =
			\begin{cases}
				1 & \text{if } m = (-1-k,k,0) \text{ for } k \in \mathbb{Z}, \\[6pt]
				0 & \text{otherwise.}
			\end{cases}
		\end{equation}
	\end{lemma}
	\begin{proof}
		We want to study the affine chart $U_{1,8}$ corresponding to the 2-dimensional cone $\tau_1= \cone{v_{1,\pm}}$. 
		
		Using the notation introduced in \S\ref{sec:preliminaries_qG_index1}, this affine variety has Gorenstein index $r=2$ with $m_0=(1,0,0)$. 
		
		We can construct the index 1 cover using the sublattice 
		\[
		N'=\{ (n_1, n_2, n_3)^T \in N \text{ such that } n_1\in 2 \ZZ\} \subseteq N.
		\]
		A $\ZZ$-basis of $N'$ is $(2, 0, 0)^T, (0,1,0)^T, (0,0,1)^T$.
		Its dual lattice $M'$ has $\ZZ$-basis $(1/2, 0, 0)$, $(0,1,0)$, $(0,0,1)$.
		
		If we consider the cone $\tau_1$ inside $N'$, the corresponding affine variety is 
		\[
		\TV{\tau_1}{N'} = U_{1,8}' = \Spec \dfrac{\CC[z_1,z_2,z_3,t^{\pm}]}{(z_1z_3-z_{2}^2)},
		\]
		where 
		\begin{align*}
			t=& \chi^{(-1,1,0)}, \\
			z_1=&\chi^{(1/2,0,-1)}, \\
			z_2=&\chi^{(1/2,0,0)}, \\
			z_3=&\chi^{(1/2,0,1)}.
		\end{align*}
		
		Consider the action of the cyclic group $\mu_4=\langle g\rangle$ on the affine plane $\AA^2= \Spec \CC [u,v]$ given by
		\[
		g(u,v)=(\rmi u, \rmi v),
		\qquad \rmi ^2=-1.
		\]
		
		Since $g^2(u,v)=(-u,-v)$, the subgroup $\mu_2=\langle g^2\rangle$ is precisely the group defining the quotient
		\[
		\mathbb C[u,v]^{\mu_2}
		=\mathbb C[u^2,uv,v^2]
		\simeq
		\frac{\mathbb C[z_1,z_2,z_3]}{(z_1z_3-z_2^2)},
		\]
		where $z_1=u^2$, $z_2=uv$, and $z_3=v^2$.
		
		The action of $\mu_4$ on $\mathbb C[u,v]$ induces an action of the quotient group $\mu_4/\mu_2 \simeq \mu_2$ on the invariant ring $\mathbb C[u,v]^{\mu_2}$. Indeed, applying $g$ to the generators $z_1,z_2,z_3$, we obtain
		\begin{align*}
			g(z_1)= g(u^2)= \rmi^2 u^2= -z_1,\\
			g(z_2)= g(uv)=  \rmi^2 uv= -z_2,\\
			g(z_3)=g(v^2)= \rmi^2v^2= -z_3.
		\end{align*}
		The coordinate $t$ of the torus $\Gm$ remains fixed.
		Hence, taking the quotient of $U'_{1,8}$ by this action is equivalent to taking the quotient of $\AA^2$ by the whole group $\mu_4$. Therefore,
		\[
		\begin{aligned}
			U_{1,8}
			&\simeq \Gm \times
			\operatorname{Spec}
			\left(\CC[u,v]^{\mu_2}\right)^{\mu_4/\mu_2}\\
			&\simeq \Gm\times
			\Spec\CC[u,v]^{\mu_4}\\
			&\simeq \Gm\times\frac14(1,1).
		\end{aligned}
		\]
		
		Let
		\[
		S := \CC[z_1,z_2,z_3,t^{\pm}]/(z_1z_3 - z_2^2).
		\]
		We may view $S$ as a tensor product
		\[
		S \simeq A \otimes_\CC B,
		\]
		where
		\[
		A := \CC[z_1,z_2,z_3]/(z_1z_3 - z_2^2)
		\quad\text{and}\quad
		B := \CC[t^{\pm}].
		\]
		Accordingly, the corresponding affine variety decomposes as a product
		\[
		\Spec S
		\simeq
		\Spec A \times \Spec B
		=
		W \times \Gm,
		\]
		where we set
		\[
		W := \Spec A.
		\]
		Using this decomposition and denoting by $p_W$ and $p_{\Gm}$ the projections onto the two factors, we get
		\[
		\Omega^1_{U_{1,8}'} = p^{*}_W \Omega^1_{W} \oplus p^{*}_{\Gm} \Omega^1_{\Gm}.
		\]
		As modules, we have
		\[
		\Omega^1_{S} =  (\Omega^1_{A} \otimes_A S) \oplus (\Omega^1_{B}\otimes_B S) = (\Omega^1_{A} \otimes_A S) \oplus S,
		\]
		where in the last equality we have used that the $B$-module $\Omega^1_B$ is free of rank one.
		Applying the additive functor $ \Ext^1_S(\cdot,S)$, we obtain
		\[ \TT_{U'_{1,8}}=
		\Ext^1_S(\Omega^1_S,S)
		\simeq
		\Ext^1_S(\Omega^1_A \otimes_A S,S)
		\oplus
		\Ext^1_S(S,S) = \Ext^1_S(\Omega^1_A \otimes_A S,S).
		\]
		Since $S = A \otimes_{\CC} \CC[t^{\pm 1}]$ is flat over $A$ and $\Omega^1_A$ is a finitely presented $A$-module, \cite[Theorem~7.11]{Matsumura} yields
		\[
		\Ext^1_S(\Omega^1_A \otimes_A S,S)
		\simeq
		\Ext^1_A(\Omega^1_A,A) \otimes_A S.
		\]
		
		We now need to calculate $\Ext^1_A(\Omega^1_A,A) \simeq \TT^1_{W}$.
		Using the conormal sequence of the closed embedding $W \hookrightarrow \AA^3$ given by the ideal $(z_1z_3- z_2^2)$, we obtain a free resolution of $\Omega^1_{W}$:
		\begin{equation}\label{eq:res}
			0 \longrightarrow \cO_{W}
			\xrightarrow{\begin{pmatrix}
					z_3 \\ -2z_2 \\ z_1
			\end{pmatrix}}
			\cO_{W}^{\oplus 3} \longrightarrow \Omega^1_{W} \longrightarrow 0.
		\end{equation}
		If we now apply the functor $\Hom_{\cO_{W}}(\cdot, \cO_{W})$ to \eqref{eq:res}, we get:
		\[
		\TT^1_{W} = \cO_{W} / (z_1, -2z_2, z_3) 
		= \CC[z_1,z_2,z_3]/(z_1z_3-z_2^2, z_1,z_2,z_3) \simeq \CC.
		\]
		Hence, 
		\[
		\TT^1_{U_{1,8}'} = \CC \otimes_A S = \CC \otimes_A A \otimes_{\CC} \CC[t^{\pm}] = \CC[t^{\pm}].
		\]
		Therefore, the $\CC$-vector space $\TT^1_{U_{1,8}'}$ has a homogeneous basis
		\[
		\{t^k \mid k\in \ZZ \} 
		\]
		with respect to the torus action.
		Therefore, since $U'_{1,8}$ is an affine hypersurface, it is unobstructed and the miniversal (qG-)deformation of $U'_{1,8}$ is
		\[
		z_1z_3-z_2^2+\sum_{k\in \ZZ}c_k t^k=0.
		\]
		Since the degree of the equation $z_1z_3-z_2^2$ is $(1,0,0) \in M$ and the degree of $t$ is $(-1,1,0) \in M$, the degree of the deformation parameter $c_k$ is
		\[
		\deg c_k =(1,0,0)-(-k,k,0)=(1+k,-k,0)
		\]
		for each $k \in \ZZ$.
		
		By duality, the corresponding degree in $\TT^1_{U'_{1,8}}$ is the opposite, namely $(-1-k,k,0)$.
		Then, the eigenspace $\TT_{U'_{1,8}}^1(m)$ with respect to the character $m$ has dimension 1 when $m=(-1-k,k,0)$ for any $k\in \ZZ$.
		
		Since the index of $U_{1,8}$ is $2$, $\TTqG{1}{U_{1,8}}$ is the $\mu_2$-invariant part of $\TT^1_{U'_{1,8}}$ (see \S\ref{sec:preliminaries_qG_index1}). Thus, only those degrees lying in the sublattice $M\subset M'$ are considered and we get \eqref{eq:T1_U18}.
	\end{proof}

	\begin{lemma} \label{lem:U_12}
		The affine chart $U_{1,2}$ is a Gorenstein affine variety isomorphic to
		\[
		\Gm \times \Spec \dfrac{\CC[z_1,z_2,z_3]}{(z_1z_3-z_2^2)}.
		\]
		We have:
		\begin{equation} \label{eq:T1_U12}
			\dim \TTqG{1}{U_{1,2}}(m) =
			\begin{cases}
				1 & \text{if } m = (-2-2k,k,0) \text{ for } k \in \ZZ, \\[6pt]
				0 & \text{otherwise.}
			\end{cases}
		\end{equation}
		The miniversal (qG-)deformation of $U_{1,2}$ is given by the equation
		\[
		z_1z_3-z_2^2+\sum_{k\in \ZZ} c_k t^k=0,
		\]
		where the action of the torus $T_N$ has the following weights:
		\begin{align*}
			\deg(c_k) &= (2+2k,-k,0), \\
			\deg (t) &= (-2,1,0), \\
			\deg(z_1) &= (1,0,-1), \\
			\deg(z_2) &= (1,0,0), \\ 
			\deg(z_3) &= (1,0,1).
		\end{align*}
	\end{lemma}
	\begin{proof}
		
		The chart corresponding to the 2-dimensional cones $\tau_{2}$ is a Gorenstein variety of the form
		\[
		\Spec \dfrac{\CC[z_1,z_2,z_3,t^{\pm}]}{(z_1z_3-z_2^2)} = \Gm \times \Spec \dfrac{\CC[z_1,z_2,z_3]}{(z_1z_3-z_2^2)},
		\]
		where
		\begin{align*}
			t   &= \chi^{(-2,1,0)}, &\qquad 
			z_1 &= \chi^{(1,0,-1)}, \\
			z_2 &= \chi^{(1,0,0)},  &\qquad 
			z_3 &= \chi^{(1,0,1)}.
		\end{align*}
		Using the same methods as in the proof of Lemma~\ref{lem:U18}, we obtain that the $\CC$-vector space $\TTqG{1}{U_{1,2}}$ has basis $\{t^k | k\in \ZZ \}$.
		Therefore, since $U_{1,2}$ is an affine hypersurface, it is unobstructed and the miniversal (qG-)deformation of $U_{1,2}$ is
		\[
		z_1z_3-z_2^2+\sum_{k\in \ZZ}c_k t^k=0.
		\]
		
		We now determine the degrees of the parameters $c_k$. Since 
		\[ 
		\deg(z_1z_3-z_2^2)=(2,0,0) \qquad\text{and}\qquad \deg(t^k)=(-2k,k,0), 
		\] 
		if we want the defining equation of the family to be homogeneous, we have that 
		\[ \deg(c_k)+\deg(t^k)=(2,0,0), 
		\] 
		and hence 
		\[ 
		\deg(c_k)=(2+2k,-k,0). 
		\]
		By duality, the corresponding degree in $\TTqG{1}{U_{1,2}}$ is the opposite of the degree of the deformation parameter. Hence, we get \eqref{eq:T1_U12}.
	\end{proof}
	
	\begin{lemma} \label{lem:U23}
		The affine chart $U_{2,3}$ is a Gorenstein affine variety isomorphic to
		\[
		\Gm \times \Spec \dfrac{\CC[z_1,z_2,z_3]}{(z_1z_3-z_2^2)}.
		\]
		We have:
		\[
		\dim \TTqG{1}{U_{2,3}}(m) =
		\begin{cases}
			1 & \text{if } m = (k,k-2,0) \text{ for } k \in \ZZ, \\[6pt]
			0 & \text{otherwise.}
		\end{cases}
		\]
		The miniversal qG-deformation of $U_{2,3}$ is given by the equation
		\[
		z_1z_3-z_2^2+\sum_{k\in \ZZ} c_k t^k=0,
		\]
		where the action of the torus $T_N$ has the following weights:
		\begin{align*}
			\deg(c_k) &= (-k,2-k,0), \\
			\deg (t) &= (1,1,0), \\
			\deg(z_1) &= (0,1,-1), \\
			\deg(z_2) &= (0,1,0), \\ 
			\deg(z_3) &= (0,1,1).
		\end{align*}
	\end{lemma}
	\begin{proof}
		The chart corresponding to the 2-dimensional cone $\tau_{3}$ is a Gorenstein variety isomorphic to
		\[
		\Spec \dfrac{\CC[z_1,z_2,z_3,t^{\pm}]}{(z_1z_3-z_2^2)} = \Gm \times \Spec \dfrac{\CC[z_1,z_2,z_3]}{(z_1z_3-z_2^2)},
		\]
		where:
		\begin{align*}
			t=& \chi^{(1,1,0)}, &\qquad
			z_1=&\chi^{(0,1,-1)}, \\
			z_2=&\chi^{(0,1,0)}, &\qquad
			z_3=&\chi^{(0,1,1)}.
		\end{align*}
		
		Using the same methods as in the proof of Lemma~\ref{lem:U18} and of Lemma~\ref{lem:U_12}, we obtain that the $\CC$-vector space $\TTqG{1}{U_{2,3}}$ has basis $\{t^k | k\in \ZZ \}$, and the degrees in the decomposition are as described in the statement of the lemma being proved.
	\end{proof}
	
	\begin{lemma} \label{lem:U1}
		The affine variety $U_{1}$ is $\QQ$-Gorenstein with Gorenstein index 2. Its index 1 cover is isomorphic to
		\[
		\Spec \dfrac{\CC[z_1,z_2,t_1,t_2]}{(z_1z_2 - t_1^2 t_2^2)}.
		\]
		We have $\TTqG{2}{U_{1}} = 0$ and
		\begin{equation} \label{eq:T1_U1}
			\dim \TTqG{1}{U_{1}}(m)=
			\begin{cases}
				1 & \text{if } m=(-k,k-1,0), \; k\geq 0,\\
				1 & \text{if } m=(2k,-1-k,0), \; k\geq 1,\\
				0 & \text{otherwise.}
			\end{cases}
		\end{equation}
	\end{lemma}
	\begin{proof}
		The 3-dimensional cone $\sigma_{1}$ gives rise to a variety of Gorenstein index $2$ with $m_0=(0,1,0)$ (see \S\ref{sec:preliminaries_qG_index1}).
		
		Consider the sublattice
		\[
		N'=\{(n_1,n_2,n_3)^T \in N \mid n_2\in 2\ZZ\},
		\]
		which has basis $(1,0,0)^T, (0,2,0)^T, (0,0,1)^T$. 
		
		The index 1 cover of $U_1$ is the hypersurface
		\[
		\TV{\sigma_1}{N'}= U'_{1}= \Spec \dfrac{\CC[z_1,z_2,t_1,t_2]}{(z_1z_2 - t_1^2 t_2^2)},
		\]
		where
		\begin{align*}
			z_1=&\chi^{(0,1/2,-1)}, \\
			z_2=&\chi^{(0,1/2,1)}, \\
			t_1=&\chi^{(-1,1,0)}, \\ 
			t_2=& \chi^{(1,-1/2,0)}.
		\end{align*}
		Since $U_1'$ is affine and l.c.i., $\TT^2_{U'_1}=0$. Hence $\TTqG{2}{U_1}=(\TT^2_{U'_1})^{\mu_2}=0$.
		
		Using the conormal sequence of the closed embedding $U_{1}' \hookrightarrow \AA^4$ given by the ideal $(z_1z_2 - t_1^2 t_2^2)$, we have a free resolution of $\Omega^{1}_{U'_{1}}$:
		\begin{equation*} \label{res2}
			0 \longrightarrow \cO_{U'_{1}}
			\xrightarrow{\begin{pmatrix}
					-2t_1t_2^2 \\
					-2t_1^2t_2 \\
					z_2 \\
					z_1
			\end{pmatrix}}
			\cO_{U'_{1}}^{\oplus 4}
			\longrightarrow \Omega^1_{U'_{1}}
			\longrightarrow 0 .
		\end{equation*}
		After applying the functor $\Hom_{\cO_{U'_{1}}}(\cdot,\cO_{U'_{1}})$ to this sequence, we have that
		\[
		\TT^1_{U'_{1}}
		=
		\dfrac{\cO_{U'_{1}}}{(t_1t_2^2, t_1^2t_2, z_1, z_2)}
		=
		\dfrac{\CC[z_1,z_2,t_1,t_2]}{(t_1t_2^2, t_1^2t_2, z_1, z_2)} = \dfrac{\CC[t_1,t_2]}{(t_1t_2^2, t_1^2t_2)}
		\]
		In this case, the basis of $\TT^1_{U'_{1}}$ is
		\[
		\{1,t_1t_2\} \cup \{t_1^k \mid k \ge 1\} \cup \{t_2^k \mid k \ge 1\}.
		\]
		Since the defining equation $z_1z_2-t_1^2t_2^2$ is homogeneous of degree $(0,1,0)$, while the degrees of $t_1$ and $t_2$ are respectively $(-1,1,0)$ and $(1,-1/2,0)$, the degrees of the deformation parameters are determined by requiring the deformed equation to remain homogeneous of degree $(0,1,0)$. So the degrees of these parameters are:
		\begin{align*}
			(0,1,0), \\
			(0,1,0) - (-1,1,0) - (1,-1/2,0)=& (0,1/2,0),\\
			(0,1,0) - (-k, k, 0)=&(k,1-k,0), \\
			\left(0,1,0\right) - \left(k,-\frac{k}{2},0\right)
			=& \left(-k,1+\frac{k}{2},0\right).
		\end{align*}
		By duality, the corresponding degrees in $\TT^1_{U'_{1}}$ are the opposite and we have 
		\[
		\dim \TT^1_{U'_{1}}(m)=
		\begin{cases}
			1 & \text{if } m=(0,-1/2,0), \\
			1 & \text{if } m=(-k,k-1,0), \; k\geq 0,\\
			1 & \text{if } m=\left(k,-1-k/2,0 \right), \; k\geq 1,\\
			0 & \text{otherwise.}
		\end{cases}
		\]
		Since the Gorenstein index of $U_1$ is $2$, $\TTqG{1}{U_1}$ is the $\mu_2$-invariant part of $\TT^1_{U'_1}$ (see \S\ref{sec:preliminaries_qG_index1}). Thus, only those degrees lying in the sublattice $M\subset M'$ are considered and we get \eqref{eq:T1_U1}.
	\end{proof}
	
	\begin{lemma} \label{lem:U2}
		The affine variety $U_{2}$ is $\QQ$-Gorenstein with Gorenstein index 3. Its index 1 cover is isomorphic to
		\[
		\Spec \dfrac{\CC[z_1,z_2,t_1,t_2]}{(z_1z_2 - t_1^2 t_2^2)}.
		\]
		We have $\TTqG{2}{U_{2}} = 0$ and
		\[
		\dim \TTqG{1}{U_2}(m)=
		\begin{cases}
			1 & \text{if } m=(k+2,k,0), \; k\geq -1,\\
			1 & \text{if } m=(-2-2k,k,0), \; k\geq -1,\\
			0 & \text{otherwise.}
		\end{cases}
		\]
	\end{lemma}
	\begin{proof}
		The variety defined by the 3-dimensional cone $\sigma_{2}$ has Gorenstein index $3$ with $m_0=(-1,2,0)$, using the same notation as in \S\ref{sec:preliminaries_qG_index1}.
		
		Using the sublattice
		\[
		N'=\{(n_1,n_2,n_3)^T \in N \mid -n_1+2n_2\in 3\ZZ\},
		\]
		which has basis $(3,0,0)^T, (2,1,0)^T, (0,0,1)^T$, we can construct the index 1 cover $U'_{2}$ of $U_2$.
		The dual lattice $M'$ has basis $(1/3, -2/3, 0), (0,1,0), (0,0,1)$.
		We have that
		\[
		\TV{\sigma_{2}}{N'}= U'_{2}= \Spec \dfrac{\CC[z_1,z_2,t_1,t_2]}{(z_1z_2 - t_1^2 t_2^2)},
		\]
		with
		\begin{align*}
			z_1 &= \chi ^{(-1/3,2/3,-1)},\\
			z_2 &= \chi ^{(-1/3,2/3,1)},\\
			t_1 &= \chi ^{(-2/3,1/3,0)},\\
			t_2 &= \chi ^{(1/3,1/3,0)}.
		\end{align*}
		Following the same steps and doing similar calculations as in the proof of Lemma~\ref{lem:U1}, we obtain the statement.
	\end{proof}
	
	Thanks to the symmetries of the polytope $P$ (see item~(ii) of Theorem~\ref{thm:tor_var_X}), the eigenspace decomposition of $\TTqG{1}{}$ for all other charts of $U$, defined in \S\ref{sec:notation_cones_charts}, can be completely recovered from the cases above (see Lemmas~\ref{lem:U18}, \ref{lem:U_12}, \ref{lem:U23}, \ref{lem:U1}, \ref{lem:U2}).
	We collect the degrees of the non-zero homogeneous components of $\TTqG{1}{U_{i,j}}$ and $\TTqG{1}{U_i}$ in Table~\ref{tab:degrees_T1_Uij} and Table~\ref{tab:degrees_T1_Ui}, respectively.
	
	\begin{table}
		\centering
		\caption{Degrees of the non-zero homogeneous components of
			$\TTqG{1}{U_i}$.}
		\label{tab:degrees_T1_Ui}
		\renewcommand{\arraystretch}{2.5}
		\rowcolors{2}{gray!12}{white}
		\setlength{\tabcolsep}{9pt}
		\begin{tabular}{c c}
			\hline
			chart & degrees $m$ with $\TTqG{1}{U_i}(m)\neq 0$ \\
			\hhline{--}
			$U_1$ & $\begin{aligned}
				&(-k,k-1,0), && k\geq 0,\\
				&(2k,-1-k,0), && k\geq 1
			\end{aligned}$ \\
			
			$U_2$ & $\begin{aligned}
				&(k+2,k,0), && k\geq -1,\\
				&(-2-2k,k,0), && k\geq -1
			\end{aligned}$ \\
			
			$U_3$ & $\begin{aligned}
				&(-k,-k-2,0), && k\geq -1,\\
				&(-k,2+2k,0), && k\geq -1
			\end{aligned}$ \\
			
			$U_4$ & $\begin{aligned}
				&(1-k,k,0), && k\geq 0,\\
				&(1+k,-2k,0), && k\geq 1
			\end{aligned}$ \\
			
			$U_5$ & $\begin{aligned}
				&(k,1-k,0), && k\geq 0,\\
				&(-2k,1+k,0), && k\geq 1
			\end{aligned}$ \\
			
			$U_6$ & $\begin{aligned}
				&(-2-k,-k,0), && k\geq -1,\\
				&(2+2k,-k,0), && k\geq -1
			\end{aligned}$ \\
			
			$U_7$ & $\begin{aligned}
				&(k,k+2,0), && k\geq -1,\\
				&(k,-2-2k,0), && k\geq -1
			\end{aligned}$ \\
			
			$U_8$ & $\begin{aligned}
				&(k-1,-k,0), && k\geq 0,\\
				&(-1-k,2k,0), && k\geq 1
			\end{aligned}$ \\
			\hline
		\end{tabular}
	\end{table}

	\begin{table}
		\centering
		\caption{Degrees of the non-zero homogeneous components of
			$\TTqG{1}{U_{i,j}}$.}
		\label{tab:degrees_T1_Uij}
		\renewcommand{\arraystretch}{1.5}
		\begin{tabular}{c@{\qquad}c}
			\hline
			chart & degrees $m$ with $\TTqG{1}{U_{i,j}}(m)\neq 0$ \\
			\hline
			$U_{1,8}$ & $(-1-k,k,0)$, \quad $k\in\ZZ$ \\
			$U_{1,2}$ & $(-2-2k,k,0)$, \quad $k\in\ZZ$ \\
			$U_{2,3}$ & $(k,k-2,0)$, \quad $k\in\ZZ$ \\
			$U_{3,4}$ & $(-k,2+2k,0)$, \quad $k\in\ZZ$ \\
			$U_{4,5}$ & $(1+k,-k,0)$, \quad $k\in\ZZ$ \\
			$U_{5,6}$ & $(2+2k,-k,0)$, \quad $k\in\ZZ$ \\
			$U_{6,7}$ & $(-k,2-k,0)$, \quad $k\in\ZZ$ \\
			$U_{7,8}$ & $(k,-2-2k,0)$, \quad $k\in\ZZ$ \\
			\hline
		\end{tabular}
	\end{table}

	\begin{lemma} \label{lem:Hi_U}
		Let $U$ be the non-affine open subscheme of $X$ defined in \S\ref{sec:notation_cones_charts}. Then $\rH^0(U, \cTqG{1}{U}) = \CC^{6} $ and
		$\rH^1(U, \cTqG{1}{U}) = 0.$ More precisely,
		\[
		\dim \rH^0(U, \cTqG{1}{U})(m) =
		\begin{cases}
			1 & \text{if } m = \pm (0,1,0), \\
			1 & \text{if } m = \pm (1, 0,0), \\
			1 & \text{if } m =  \pm (1,-1,0), \\
			0 & \text{otherwise}.
		\end{cases}
		\]
	\end{lemma}
	\begin{proof}
		Whenever $j\notin \{i, i \pm 1\}$, the chart $U_{i,j} = U_i \cap U_j$ is the affine toric variety associated with the zero cone, which is the algebraic torus $T_N$. Since $T_N$ is a smooth affine variety, the sheaf $\cTqG{1}{U}$ vanishes when restricted to $U_{i,j}$.
		As a result, when we compute the \v{C}ech complex of $\cTqG{1}{U}$ with respect to the affine cover $\{U_i \mid i \in \ZZ/8\ZZ \}$, all contributions from intersections of three or more distinct charts disappear, while the only surviving pairwise intersections are those of the form $U_{i, i+1}$, and the complex has nontrivial terms only in degrees 0 and 1:
		\[
		\bigoplus_{i \in \mathbb{Z}/8\mathbb{Z}} \TTqG{1}{U_i} 
		\xrightarrow{\quad d \quad} 
		\bigoplus_{i \in \mathbb{Z}/8\mathbb{Z}} \TTqG{1}{U_{i,i+1}}
		\]
		As this is a homomorphism of $M$-graded $\CC$-vector spaces, we analyse the homogeneous components of this complex for every degree $m \in M$.
		The degrees of the non-zero homogeneous components of
		$\TTqG{1}{U_i}$ and $\TTqG{1}{U_{i,i+1}}$ are collected in
		Tables~\ref{tab:degrees_T1_Ui} and
		\ref{tab:degrees_T1_Uij}, respectively.
		\begin{itemize}
			
			\item If $ m=\pm(0,1,0)$, the complex is
			\[
			\CC^3 \xrightarrow{\pm \begin{pmatrix} 0 & 1 & -1 \\ -1 & 0 & 1\end{pmatrix}} \CC^2.
			\]
			
			\item If $m= \pm(1, 0, 0)$, the complex is
			\[
			\CC^3 \xrightarrow{\pm \begin{pmatrix} -1 & 1 & 0 \\ 0 & -1 & 1\end{pmatrix}} \CC^2.
			\]
			
			\item The complex is 
			\[
			\CC^2 \xrightarrow{\pm \begin{pmatrix} 1 & 0 \\ 0 & -1\end{pmatrix}} \CC^2.
			\]
			in the following cases:
			\[
			m= 
			\begin{cases}
				\pm(-2, 0, 0), \\
				\pm(0, -2, 0), \\
				\pm(3, -4, 0), \\
				\pm(-2, 2,0), \\
				\pm(-4, 3,0).
			\end{cases}
			\]
			
			\item  If $m= \pm(-1, 1,0)$, we get
			\[
			\CC^2 \xrightarrow{\pm \begin{pmatrix} 1 & -1\end{pmatrix}} \CC.
			\]
			
			\item In all the other degrees for which at least one of the graded pieces are non-zero, the complex is 
			\[
			\CC \xrightarrow{\pm 1} \CC.
			\]
		\end{itemize}
		
		Since $\rH^0(U, \cTqG{1}{U}) = \ker d$ we have that
		\[
		\dim \rH^0(U, \cTqG{1}{U})(m) =
		\begin{cases}
			1 & \text{if } m = \pm (0,1,0), \\
			1 & \text{if } m = \pm (1, 0,0), \\
			1 & \text{if } m =  \pm (1,-1,0), \\
			0 & \text{otherwise}.
		\end{cases}
		\]
		
		Moreover, the differential $d$ is surjective, so  
		\[
		\rH^1(U, \cTqG{1}{U}) = \coker d = 0. \qedhere
		\]
	\end{proof}
	
	\begin{remark}\label{rem:degrees_chart}
		If $m\in M$ is one of the six degrees appearing in Lemma~\ref{lem:Hi_U}, then $\TTqG{1}{U_{i,i+1}}(m)$ is non-zero for one value or two value of $i \in \ZZ/8\ZZ$.
		The complete list, together with the corresponding values of $k$ as in Table~\ref{tab:degrees_T1_Uij}, is the following.
		\[
		\begin{array}{c@{\qquad}c}
			\text{degree }m & (\text{chart},k) \\[1.5pt]
			\hline \rule{0pt}{10pt}
			(0,1,0) & (U_{4,5},-1),\ (U_{5,6},-1) \\
			(0,-1,0) & (U_{1,8},-1),\ (U_{1,2},-1) \\
			(1,0,0) & (U_{3,4},-1),\ (U_{4,5},0) \\
			(-1,0,0) & (U_{1,8},0),\ (U_{7,8},-1) \\
			(1,-1,0) & (U_{2,3},1) \\
			(-1,1,0) & (U_{6,7},1)
		\end{array}
		\]
	\end{remark}

	\begin{proof}[Proof of items~(iv)-(vi) of Theorem~\ref{thm:tor_var_X}]
		
		In order to compute the base space of the miniversal qG-deformation of $X$, we consider the local-to-global spectral sequence for $\Ext$:
		\[
		E^{p,q}_2 = \rH^p(X, \cTqG{q}{X}) \;\;\Rightarrow\;\; \TTqG{p+q}{X}.
		\]
		We first determine $E^{p,q}_2$ that will be needed below.
		
		We begin with $E^{1,1}_2 = \rH^1(X,\cTqG{1}{X})$.
		We compute the Čech complex of $\cTqG{1}{X}$ with respect to the affine open cover $\cU=\{U_+, U_-, U_1, \dots, U_8 \}$ of $X$.
		Consider the affine open cover $\cV= \{ U_1, \dots, U_8 \}$ of $U$. Since $\cTqG{1}{X}$ is supported on the singular locus of $X$, and the singularities contained in $U_+$ and $U_-$ are isolated, we have that the sheaf $\cTqG{1}{X}$ vanishes on every intersection of at least two distinct open subsets and involving $U_+$ or $U_-$.
		Hence, for $p\geq 1$, the components of the Čech cochain groups corresponding to such intersections vanish, and therefore
		\[
		\check{C}^p(\cU,\cTqG{1}{X})
		\simeq
		\check{C}^p(\cV,\cTqG{1}{X}),
		\qquad p\geq 1.
		\]
		Hence the forgetful homomorphism of complexes
		\[
		\check{C}^{\boldsymbol{\cdot}}(\cU,\cTqG{1}{X}) \longrightarrow \check{C}^{\boldsymbol{\cdot}}(\cV,\cTqG{1}{X})
		\]
		is bijective in degrees $p\geq 1$ and surjective in degree $p=0$. Hence, $U_+$ and $U_-$ do not contribute to the first cohomology and we have
		\[
		\rH^1(X, \cTqG{1}{X}) \simeq \rH^1(U, \cTqG{1}{X}).
		\]
		By Lemma~\ref{lem:Hi_U}, we conclude that
		\[
		E^{1,1}_2=\rH^1(X, \cTqG{1}{X})=0.
		\]
		
		We now consider the remaining cohomology groups. Since $X$ is a toric Fano variety, $\rH^p(X,\cT^0_X)=0$ for $p>0$ by \cite[\S5]{Totaro} (see also \cite[Lemma~4.4]{Pet22b}).
		Moreover, $\cT^0_X = \cTqG{0}{X}$. Hence
		\[
		E^{p,0}_2= \rH^p(X,\cTqG{0}{X})=0
		\qquad\text{for }p>0.
		\]
		
		If $q>0$, the sheaves $\cTqG{q}{X}$ are supported on the singular locus of $X$, which is $1$-dimensional, so it follows that 
		$\rH^p(X,\cTqG{q}{X})=0$ for $p\geq 2$. In particular,
		\[
		E^{2,1}_2=\rH^2(X,\cTqG{1}{X})=0, \qquad
		E^{2,2}_2=\rH^2(X,\cTqG{2}{X})=0.
		\]
		
		To calculate $\TTqG{1}{X}$ and $\TTqG{2}{X}$, we consider the two diagonals $p+q=1$ and $p+q=2$. Along the diagonal $p+q=1$, the only possible terms are
		\[
		E^{1,0}_2=\rH^1(X,\cTqG{0}{X}),
		\qquad
		E^{0,1}_2=\rH^0(X,\cTqG{1}{X}),
		\]
		while along the diagonal $p+q=2$, they are
		\[
		E^{2,0}_2=\rH^2(X,\cTqG{0}{X}),
		\qquad
		E^{1,1}_2=\rH^1(X,\cTqG{1}{X}),
		\qquad
		E^{0,2}_2=\rH^0(X,\cTqG{2}{X}).
		\]
		By the results above, $E^{1,0}_2, E^{2,0}_2, E^{1,1}_2$ vanish.
		Therefore
		\[
		\TTqG{1}{X} = \rH^0(X,\cTqG{1}{X}), \qquad \TTqG{2}{X}=\rH^0(X,\cTqG{2}{X}).
		\]
		
		Now observe that by the analysis in \S\ref{sec:defU}, we have $\cTqG{2}{X}\vert_{U_i}=0$ for every $i=1,\dots,8$ because the index 1 cover of $U_i$ is l.c.i., so $\cTqG{2}{U}=0$.
		Since the charts $U_+$ and $U_-$ are affine and Gorenstein, we have
		\[
		\TT^{j}_{U_+} = \TTqG{j}{U_+}=\rH^0(U_+, \cTqG{j}{U_+}), \qquad \TT^{j}_{U_-} = \TTqG{j}{U_-}=\rH^0(U_-, \cTqG{j}{U_-}),
		\]
		for $j=1,2$
		
		Let $\{U,U_+,U_-\}$ be the non-affine open cover of $X$.
		Since $\cTqG{1}{X}$ and $\cTqG{2}{X}$ are supported on the singular locus of $X$ and the singularities contained in $U_+$ and $U_-$ are isolated, these sheaves vanish on the pairwise intersections of the open sets in $\{U,U_+,U_-\}$.
		Using this cover, we obtain
		\begin{align*}
			\TTqG{1}{X} &= \rH^0(X, \cTqG{1}{X}) = \rH^0(U, \cTqG{1}{U}) \oplus \rH^0(U_+, \cTqG{1}{U_+}) \oplus \rH^0(U_-, \cTqG{1}{U_-}) \\
			&\phantom{{}= \rH^0(X, \cTqG{1}{X})} = \rH^0(U, \cTqG{1}{U}) \oplus \TT^{1}_{U_{+}} \oplus \TT^{1}_{U_{-}}
		\end{align*}
		and
		\begin{align*}
			\TTqG{2}{X} &= \rH^0(X, \cTqG{2}{X}) = \rH^0(U, \cTqG{2}{U}) \oplus \rH^0(U_+, \cTqG{2}{U_+}) \oplus \rH^0(U_-, \cTqG{2}{U_-}) \\
			&\phantom{{}= \rH^0(X, \cTqG{2}{X})} = \TT^{2}_{U_{+}} \oplus \TT^{2}_{U_{-}}.
		\end{align*}
		Consider the morphism
		\[
		\Phi \colon \DefqG{X} \longrightarrow \DefqG{U_+} \times \DefqG{U_-} = \Def{U_+} \times \Def{U_-}
		\]
		which maps a qG-deformation of $X$ to the pair consisting of its restrictions to $U_+$ and $U_-$.
		From the descriptions above we have that $\Phi$
		induces a bijection on the obstruction space $\TTqG{2}{}$ and a surjection on the tangent space $\TTqG{1}{}$. By the smoothness criterion (see \cite[Corollary~3.6.10]{Manetti_def}), $\Phi$ is smooth with relative dimension $\dim \rH^0(U,\cTqG{1}{U})=6$ (see Lemma~\ref{lem:Hi_U}). Hence, using Proposition~\ref{prop:torvar_cone_over_octagon}, items (iv) and (v) of Theorem~\ref{thm:tor_var_X} are proved.
		
		Now we need to analyse how the automorphism group $\Aut(X)$ acts on the hull $R$ of $\DefqG{X}$. 
		We consider the functor
		\[
		F := \DefqG{U_+} \times \DefqG{U_-} \times \prod_{i \in \ZZ/8\ZZ} \DefqG{U_i} \times \prod_{i \in \ZZ / 8\ZZ} \DefqG{U_{i,i+1}}
		\]
		and the commutative diagram of deformation functors
		\begin{equation} \label{eq:triangolo_funtori_deformazione}
			\xymatrix{
				\DefqG{X} \ar[rd]_\Phi \ar[r]^r & F \ar[d]^\pi \\
				& \DefqG{U_+} \times \DefqG{U_-}
			}
		\end{equation}
		where $r$ maps a qG-deformation of $X$ to the tuple of its restrictions to the open subsets in the affine open cover
		\[
		\cU' = \{ U_+, U_-, U_1, \dots, U_8, U_{1,2}, \dots, U_{8,1} \}
		\]
		of $X$
		and $\pi$ is the projection onto the first two factors.
		Let
		\[
		\sigma \colon \DefqG{U_+} \times \DefqG{U_-} \longrightarrow F
		\]
		be the section of $\pi$ given by considering the trivial deformations of $U_i$ and $U_{i,i+1}$.
		Since $\rH^1(X,\cT^0_X)=0$ by \cite[§5]{Totaro} and by a version of \cite[Lemma~2.1]{Pet22a} which is valid for qG-deformations, $r$ is injective.
		Since $\rH^2(X,\cT^0_X)=0$ by \cite[§5]{Totaro}, by the proof of \cite[Proposition~2.2]{Pet22a} and of \cite[Proposition~2.6]{Pet22a}, the image of $\sigma$ is contained in the image of $r$.
		Hence there exists a section
		\[
		s \colon \DefqG{U_+} \times \DefqG{U_-} \longrightarrow \DefqG{X}
		\]
		of $\Phi$ such that $r \circ s = \sigma$.
		
		Now let us see how the torus $T_N$ acts.
		Since every open subset in $\cU'$ is $T_N$-invariant, we have that the diagram \eqref{eq:triangolo_funtori_deformazione} is $T_N$-equivariant.
		Clearly, also $\sigma$ is $T_N$-equivariant, and consequently also $s$ is.
		This implies that we can choose an isomorphism
		\[
		R \simeq \CC \pow{y_1, \dots, y_6, x_{1,+}, \dots x_{5,+}, x_{1,-}, \dots, x_{5,-}} /I,
		\]
		where the torus $T_N$ acts on the three sets of variables $\{y_1,\dots,y_6\}$, $\{x_{1,+},\dots, x_{5,+}\}$ and $\{x_{1,-},\dots, x_{5,-}\}$ separately.
		
		The unique torus weight of the $5$-dimensional $\CC$-vector space $\TTqG{1}{U_+}$ (resp.\ $\TTqG{1}{U_-}$) is $(0,0,-1)$ (resp.\ $(0,0,1)$) by \cite[\S6.3]{Alt97} (see Proposition~\ref{prop:torvar_cone_over_octagon}(v)). This determines the torus action on the $x_{i,+}$'s (resp.\ $x_{i,-}$'s) at the first order.
		Since the weights are contained in a strongly convex cone in $M$, we deduce that the action of $T_N$ on the $x_{i,+}$'s (resp.\ $x_{i,-}$'s) is linear with weight $(0,0,1)$ (resp.\ $(0,0,-1)$).
		
		Now let us focus on the action of the torus $T_N$ on the variables $y_1,\dots,y_6$.
		These are the variables that correspond to the $6$-dimensional $T_N$-representation $\rH^0(U, \cTqG{1}{U})$ studied in Lemma~\ref{lem:Hi_U}.
		The choice of the variables $y_i$ associated with the 1-dimensional weight spaces $\TTqG{1}{X}(m_i)$ is not canonical, since each variable
		can be replaced by a nonzero scalar multiple.
		We fix the choice of the $y_i$ so that, on each chart $U_{k,l}$ associated with the weight $m_i$, the variable $y_i$ is a specific parameter in the miniversal deformation of that chart.	
		More precisely, for example we choose $y_2$ such that, as an element of the hull $R$ of $\DefqG{X}$, corresponds, via the forgetful morphism $\DefqG{X} \to \DefqG{U_{1,2}}$, to the parameter $c_{-1}$ appearing in the miniversal deformation of $U_{1,2}$ described in Lemma~\ref{lem:U_12}. The other variables $y_i$ are fixed analogously, using the corresponding charts as specified in Remark~\ref{rem:degrees_chart}.
		
		It remains to prove that the torus $T_N$ acts linearly on $y_1,\dots,y_6$.
		Let us just prove it for $y_2$. Its torus weight is $(0,1,0)$ and is the chosen coordinate on $\TTqG{1}{X}(0,-1,0)$. By looking at the miniversal (qG-)deformation of $U_{1,2}$ and at the $T_N$-action on it, we deduce that the torus acts linearly. 
		
		To conclude the proof of item~(vi), we have to describe the action of $\Aut(P)$ on the variables of the miniversal qG-deformation space.
		
		The involution $M_1$ swaps the top and bottom octagonal facets and maps the vertex $v_{i,+}$ to $v_{i,-}$ for each $i=1, \dots, 8$. 
		This geometric symmetry exchanges the maximal cones $\sigma_+$ and $\sigma_-$, which in turn swaps the corresponding affine charts $U_+$ and $U_-$. 
		This induces an isomorphism between their respective formal deformations. Therefore, we have that $M_1$ swaps $x_{i,+}$ and $x_{i,-}$ for each $i=1, \dots, 5$.
		
		To determine the action of the three involutions $M_1$, $M_2$, $M_3$  on the variables $y_1, \dots, y_6$, it is necessary to study how these automorphisms act on the charts $U_{i,j}$, from which the corresponding weights are obtained (see Table~\ref{tab:degrees_T1_Uij}). 
		The action of an automorphism on the toric variety permutes the affine charts; this geometrically induces an isomorphism between their formal deformations and a specific action on the respective deformation parameters. 
		Combining this action on the charts with the weights of the torus action on the variables (as established in Remark~\ref{rem:degrees_chart}), one can uniquely determine how $M_1$, $M_2$ and $M_3$ act on the $y_i$. 
		As a demonstrative example of this procedure, the complete calculation for the variable $y_2$ is carried out in Lemma~\ref{lem:action_matrices} below. We can apply the exact same method to recover the action on all the $y_j$ variables as described in the statement.
		
		This concludes the proof of Theorem~\ref{thm:tor_var_X}.
	\end{proof}

	\begin{lemma} \label{lem:action_matrices}
		With the notation of Theorem~\ref{thm:tor_var_X}, the matrices
		\[M_1=
		\begin{pmatrix}
			1 & 0 & 0 \\
			0 & 1 & 0 \\
			0 & 0 & -1
		\end{pmatrix}
		, \quad M_2=
		\begin{pmatrix}
			-1 & 0 & 0 \\
			0 & -1 & 0 \\
			0 & 0 & 1
		\end{pmatrix}
		, \quad M_3=
		\begin{pmatrix}
			0 & -1 & 0 \\
			-1 & 0 & 0 \\
			0 & 0 & 1
		\end{pmatrix}.
		\]
		act on the coordinate $y_2$ of $R$ via:
		\[
		M_1 \cdot y_2 = y_2, \quad M_2 \cdot y_2 = y_1, \quad M_3 \cdot y_2 = y_3.
		\]
	\end{lemma}
	\begin{proof}
		By Lemma~\ref{lem:U_12}, we have that the miniversal deformation of 
		\[
		U_{1,2}= \Gm \times \Spec \dfrac{\CC[z_1,z_2,z_3]}{(z_1z_3-z_2^2)}.
		\]
		is given by
		\[
		f_{1,2}+\sum_{k\in\mathbb Z}c_k^{(12)}t^k=0, \qquad f_{1,2}=z_1z_3-z_2^2,
		\]
		where
		\[
		t=\chi^{(-2,1,0)}, \qquad z_1=\chi^{(1,0,-1)}, \qquad z_2=\chi^{(1,0,0)}, \qquad z_3=\chi^{(1,0,1)}.
		\]
		and the degree of the deformation parameter $c_k^{(12)}$ is
		\[
		\deg(c_k^{(12)}) = (2+2k,-k,0).
		\]
		By duality, the corresponding degree in $\TTqG{1}{U_{1,2}}$ is the opposite of the degree of the deformation parameter.
		
		Since $y_2$ is the chosen coordinate of the 1-dimensional vector space $\TTqG{1}{U_{1,2}}(m)$, with $m=(0,-1,0)$, then $y_2$ corresponds to the deformation parameter of degree $(0,1,0)$, i.e.\ $c_{-1}^{(12)}$. To determine the action of the three matrices on $y_2$, we first need to understand how they act on the chart $U_{1,2}$, on its formal deformation, and, in particular, on the parameter $c_{-1}^{(12)}$.
		
		We start with the matrix $M_1$. Since this matrix swaps the top and bottom facets of the polytope $P$ in Theorem~\ref{thm:tor_var_X}, the affine chart $U_{1,2}$ is preserved by $M_1$.
		We now compute the action on the characters defining the chart, and hence on its deformation. We have
		\[
		(-2,1,0)
		\begin{pmatrix}
			1 & 0 & 0 \\
			0 & 1 & 0 \\
			0 & 0 & -1
		\end{pmatrix} 
		=(-2,1,0) \\
		\]
		\[
		(1,0,-1)
		\begin{pmatrix}
			1 & 0 & 0 \\
			0 & 1 & 0 \\
			0 & 0 & -1
		\end{pmatrix}
		=(1,0, 1)  \\
		\]
		\[
		(1,0,0) 
		\begin{pmatrix}
			1 & 0 & 0 \\
			0 & 1 & 0 \\
			0 & 0 & -1
		\end{pmatrix} 
		=(1,0,0) \\
		\]
		\[
		(1,0,1)
		\begin{pmatrix}
			1 & 0 & 0 \\
			0 & 1 & 0 \\
			0 & 0 & -1
		\end{pmatrix} 
		=(1,0,-1)
		\]
		Therefore, 
		\[
		M_1 \cdot t=t \qquad M_1 \cdot z_1=z_3 \qquad M_1 \cdot z_2=z_2 \qquad M_1 \cdot z_3 =z_1
		\]
		It follows that $f_{1,2}=z_1z_3-z_2^2$ is fixed. Thus the deformation
		\[
		f_{1,2}+\sum_{k\in\mathbb Z}c_k^{(12)}t^k=0 
		\]
		is sent to itself, and consequently
		\[
		M_1 \cdot c_k^{(12)} =c_k^{(12)}
		\]
		for every $k\in\ZZ$. In particular,
		\[
		M_1 \cdot y_2 = M_1 \cdot c_{-1}^{(12)}=c_{-1}^{(12)}=y_2.
		\]
		
		Next, consider the matrix $M_2$. This automorphism preserves the third coordinate and is a reflection
		through the origin in the plane of the first two coordinates. Hence, $M_2$ maps the chart $U_{1,2}$ to $U_{5,6}$.
		
		By the symmetries of the polytope $P$, the miniversal deformation of $U_{5,6}$ is given by the equation
		\[
		f_{5,6}+\sum_{k\in\mathbb Z}c_k^{(56)}(t')^{k}=0, \qquad f_{5,6}=z'_1z'_3-(z'_2)^{2},
		\]
		where
		\[
		t'=\chi^{(-2,1,0)}, \qquad z'_1=\chi^{(-1,0,-1)}, \qquad z'_2=\chi^{(-1,0,0)}, \qquad z'_3=\chi^{(-1,0,1)}.
		\]
		and the degree of the deformation parameter $c_k^{(56)}$ is
		\[
		\deg(c_k^{(56)}) = (-2+2k,-k,0).
		\]
		Since $y_1$ is the chosen coordinate of the 1-dimensional vector space $\TTqG{1}{U_{5,6}}(m)$ with $m=(0,1,0)$ by Table~\ref{tab:degrees_T1_Uij}, then $y_1$ is the deformation parameter of degree $-m=(0,-1,0)$, i.e.\ $c_{1}^{(56)}$.
		Repeating similar computations as for the matrix $M_1$, we obtain
		\[
		M_2 \cdot t= (t')^{-1} \qquad M_2 \cdot z_1=z'_1 \qquad M_2 \cdot z_2=z'_2 \qquad M_2 \cdot z_3 =z'_3
		\]
		It follows that $f_{1,2}$ maps to $f_{5,6}$. Consequently, the miniversal deformation is transformed into
		\[
		f_{5,6}+\sum_{k\in\mathbb Z}c_k^{(56)}(t')^{-k}=0 
		\]
		Reindexing with $l=-k$, this becomes
		\[
		f_{5,6}+\sum_{l\in\mathbb Z}c_{-l}^{(56)}(t')^{l}=0 
		\]
		It follows that
		\[
		M_2 \cdot c_k^{(12)} = c_{-k}^{(56)}.
		\]
		In particular,
		\[
		M_2 \cdot y_2 = M_2 \cdot c_{-1}^{(12)}=c_{1}^{(56)}=y_1.
		\]
		
		Finally, the matrix $M_3$ acts as a reflection with respect to the bisector of the second and fourth quadrants in the plane of the first two coordinates and maps $U_{1,2}$ to $U_{3,4}$.
		The miniversal deformation of $U_{3,4}$ is described by the equation
		\[
		f_{3,4}+\sum_{k\in\mathbb Z}c_k^{(34)}(t'')^k=0,
		\qquad
		f_{3,4}=z_1''z_3''-(z_2'')^2,
		\]
		where
		\[
		t''=\chi^{(1,-2,0)}, \qquad
		z_1''=\chi^{(0,-1,-1)}, \qquad
		z_2''=\chi^{(0,-1,0)}, \qquad
		z_3''=\chi^{(0,-1, 1)}.
		\]
		As in the previous cases, we compute the action of $M_3$ on these characters and hence on the formal deformation. Since the defining equation $f_{1,2}$ is mapped to $f_{3,4}$, we obtain the induced action on the deformation parameters.
		Moreover, $y_3$ is the chosen coordinate of the 1-dimensional vector space $\TTqG{1}{U_{3,4}}(m)$ with $m=(1,0,0)$ by Table~\ref{tab:degrees_T1_Uij},
		and it is the deformation parameter $c_{1}^{(34)}$ of degree $(-1,0,0)$. Hence,
		\[
		M_3\cdot y_2 = M_3\cdot c_{-1}^{(12)} = c_{1}^{(34)} =y_3.
		\]
		This concludes the proof of Lemma~\ref{lem:action_matrices}.
	\end{proof}
	
	\section{An obstructed K-stable Fano 3-fold}\label{sec:result_n=3}
	
	In this section we aim to prove the following theorem.
	\begin{theorem} \label{thm:result_n=3}
		There exists a K-stable Fano 3-fold such that its automorphism group is trivial and the base of its miniversal qG-deformation is the formal spectrum of
		\[
		\CC \pow{ w_1, \dots, w_9 } / (w_1^2, w_2^2).
		\]
	\end{theorem}
	The strategy of the proof is to study the local structure of the moduli stack $\cM$ of K-semistable Fano 3-folds with a fixed anticanonical volume around the point
	$[X]$ represented by the toric Fano 3-fold $X$ of Theorem~\ref{thm:tor_var_X}. 
	To this end we apply Theorem~\ref{thm:artin_stack}: we consider an affine scheme $S$, obtained by looking at the equations that define the hull of $\DefqG{X}$, together with an action of $G=\Aut(X)$ and a point $s\in S(\CC)$ fixed by $G$, in such a way that the completion of the local ring $\cO_{S,s}$ is $G$-equivariantly isomorphic to the hull of $\DefqG{X}$.
	This implies that the local structure of $\cM$ at $[X]$ is the same as the local structure of the quotient stack $[S/G]$ at the point $\bar{s}$, which is the image of $s$.
	
	By \cite[Theorem~4.19(3)]{AHR20} we obtain an affine scheme $W$ of finite type over $\CC$ with an action of $G$, a point $w\in W(\CC)$ fixed by $G$, and a diagram of \'etale morphisms
	\begin{equation*}
		\xymatrix{
			& [W/G] \ar[ld]_\Phi \ar[rd]^\Psi \\
			\cM & & [S/G]
		}
	\end{equation*}
	such that $\Psi(\bar{w})=\bar{s}$ and $\Phi(\bar{w})=[X]$, where $\bar{w}$ denotes the image of $w$ in $[W/G]$.
	By choosing a point $s'$ in a neighborhood of $s$ in $S$ and climbing the roof above via $[W/G]$, we obtain a point $[Y]$ in neighborhood of $[X]$ in $\cM$. If we choose $s'$ carefully (we need to consider a closed subscheme $S' \subseteq S$ satisfying the hypothesis of Theorem~\ref{thm:artin_stack}), we can control $\Aut(Y)$ and $\DefqG{Y}$ and we can show that the point $[Y]$ is closed in $\cM$, which implies that $Y$ is K-polystable.

	\subsection{Notation and setting} \label{sec:main_notation}
	
	Let $X$ be the toric Fano 3-fold associated to the face fan of the Fano polytope $P$ of Theorem~\ref{thm:tor_var_X}.
	Let $R$ be the hull of $\DefqG{X}$. $R$ is isomorphic to the quotient
	\[ 
	\CC \pow{y_1, \dots, y_6, x_{1,+}, \dots x_{5,+}, x_{1,-}, \dots, x_{5,-}} /I,
	\]
	where $I$ is the ideal generated by the polynomials $f_{\epsi, k}$ and $g_{\epsi, k}$ of Theorem~\ref{thm:tor_var_X}(iv).
	By a slight abuse of notation, we will use $I$ to denote both the ideal in the polynomial ring $\CC [y_1, \dots, y_6, x_{1,+}, \dots x_{5,+}, x_{1,-}, \dots, x_{5,-}]$ and the corresponding ideal in the power series ring $\CC \pow{y_1, \dots, y_6, x_{1,+}, \dots x_{5,+}, x_{1,-}, \dots, x_{5,-}}$ defining the ring $R$.
	
	Let $S$ be the closed subscheme of 
	\[
	\AA^{16}= \Spec \CC [y_1, \dots, y_6, x_{1,+}, \dots x_{5,+}, x_{1,-}, \dots, x_{5,-}]
	\]
	defined by $I$.
	
	Let $G$ be the automorphism group of $X$, which is isomorphic to $\Gm^3 \rtimes (\ZZ/2\ZZ)^3$ by Theorem~\ref{thm:tor_var_X}(iii). Let $G$ act on $S$ in the same way as it acts on $\Spf R$. By Theorem~\ref{thm:tor_var_X}(vi), we have that the three copies of $\ZZ/2\ZZ$ are generated by the matrices 
	\[M_1=
	\begin{pmatrix}
		1 & 0 & 0 \\
		0 & 1 & 0 \\
		0 & 0 & -1
	\end{pmatrix}
	, \quad M_2=
	\begin{pmatrix}
		-1 & 0 & 0 \\
		0 & -1 & 0 \\
		0 & 0 & 1
	\end{pmatrix}
	, \quad M_3=
	\begin{pmatrix}
		0 & -1 & 0 \\
		-1 & 0 & 0 \\
		0 & 0 & 1
	\end{pmatrix}, 
	\]
	while the torus $\Gm^3$ acts linearly with the following weights:
	\[
	(0,-1,0)\textrm{ on } y_1, \quad 
	(0,1,0) \textrm{ on } y_2, 
	\]
	\[
	(-1,0,0)\textrm{ on } y_3, \quad 
	(1,0,0)\textrm{ on } y_4,
	\]
	\[
	(-1,1,0)\textrm{ on } y_5, \quad 
	(1,-1,0)\textrm{ on } y_6,
	\]
	\[
	(0,0,-1) \textrm{ on } x_{1,-}, \ldots, x_{5,-} , \quad
	(0,0,1) \textrm{ on } x_{1,+}, \ldots, x_{5,+}.
	\] 
	
	Before proceeding with the proof of Theorem~\ref{thm:result_n=3}, we want to construct the slice $S'$ of Theorem~\ref{thm:artin_stack} and verify the hypothesis in
	Theorem~\ref{thm:artin_stack}(3), namely, the étaleness of the morphism $G\times^{G'} S' \to S$ at $[(e,s')]$, where $e$ is the identity of $G$ and $s'$ is a general closed point of $S'$.
	For this purpose, we temporarily assume that $G'$ is trivial, so that the morphism above is simply the restricted action morphism $G\times S' \to S$.
	The triviality of $G'$ will be proved later, in \S\ref{sec:proof}.
	
	Setting 
	\[
	\mathbf y=(y_1,\dots,y_6),\qquad
	\mathbf x=(x_{1,+},\dots,x_{5,+},x_{1,-},\dots,x_{5,-}),
	\]
	by Remark~\ref{rem:action}, the action of $G$ preserves the two sets of $\mathbf y$-variables and $\mathbf x$-variables. Hence, both factors in the above decomposition of $S$ inherit an action of $G$.
	Moreover, since the variables $y_1, \dots, y_6$ do not appear in the ideal $I$, we can write $S$ as
	\begin{equation} \label{eq:S_dec}
		S= \Spec \CC[y_1,\dots,y_6] \times \Spec \CC[x_{1,+},\dots,x_{5,+},x_{1,-},\dots,x_{5,-}]/I,
	\end{equation}
	where again by a slight abuse of notation $I$ is the ideal in the polynomial ring $\CC[x_{1,+},\dots,x_{5,+},x_{1,-},\dots,x_{5,-}]$ generated by the polynomials $f_{\epsi, k}$ and $g_{\epsi, k}$ of Theorem~\ref{thm:tor_var_X}(iv).
	
	Since we want to prove that the restricted action morphism $G\times S' \to S$ is étale at $(e,s')$ for a suitable choice of $S'$ and $s'$, we will restrict to the connected component of the identity of $G$, which is isomorphic to $\Gm^3$. Notice that we can decompose the torus as 
	\[
	\Gm^3=\Gm^2\times \Gm,
	\]
	where $\Gm^2$ and $\Gm$ are the subtori of $\Gm^3$ corresponding respectively to the first two coordinates and to the last coordinate.
	In particular, it follows that the first
	factor $\Gm^2$ of the torus acts only on the $\mathbf y$-variables and the second factor $\Gm$ acts only on the $\mathbf x$-variables.
	
	The decomposition \eqref{eq:S_dec} of $S$ allows us to construct the slice by treating
	the $\mathbf y$-variables and $\mathbf x$-variables separately. We carry
	out these two constructions in \S\ref{sec:slice_y}
	and \S\ref{sec:slice_x}, respectively, and then combine them in \S\ref{sec:proof}.
	
	In \S\ref{sec:slice_y} and \S\ref{sec:slice_x}, we will use the notation $G$, $S$, and $S'$ for objects playing analogous roles in the two settings. The meaning of this notation will be specified at the beginning of each section.

	\subsection{Slice for the y-variables}  \label{sec:slice_y} 
	We begin by considering the \textbf{y}-variables.
	
	\begin{lemma}[Slice for the \textbf{y}-variables] \label{lem:slice_y}
		Consider the group $G=\Gm^2$ acting on
		$ S:= \Spec \CC[y_1, \dots, y_6] \simeq \AA^6 $
		diagonally with weights
		\[
		(0,-1)\textrm{ on } y_1, \quad 
		(0,1) \textrm{ on } y_2, 
		\]
		\[
		(-1,0)\textrm{ on } y_3, \quad 
		(1,0)\textrm{ on } y_4,
		\]
		\[
		(-1,1)\textrm{ on } y_5, \quad 
		(1,-1)\textrm{ on } y_6.
		\]
		Let $S' \hookrightarrow S$ be the closed affine subscheme defined by the equations
		\[
		y_1=y_2, \quad y_3=y_4.
		\]
		If $s'\in S'(\CC)$ is a general point, then the morphism $G \times S' \to S$ obtained by restricting the $G$-action on $S$ to $S'$ is \'etale at the point $(e,s')$ and there is an isomorphism of complete local $\CC$-algebras
		\[
		\widehat{\cO_{S',s'}} \simeq \CC \pow{y_1, y_3, y_5, y_6}.
		\]
	\end{lemma}
	\begin{proof}
		We have $s'=(a,a,b,b,c,d) \in S'(\CC) \subseteq S(\CC) \subseteq \CC^6$, with general $a,b,c,d \in \CC$. 
		Notice that $S' \simeq \AA^4= \Spec \CC [y_1, y_3, y_5, y_6]$.
		Therefore, the morphism to be studied is the map
		\[
		\Phi \colon \Spec \CC[t_1^{\pm}, t_2^{\pm}] \times \Spec \CC[y_1, y_3, y_5, y_6] \longrightarrow \Spec \CC[y_1, y_2, y_3, y_4, y_5, y_6].
		\]
		given by
		\[
		((t_1, t_2), (y_1, y_3, y_5, y_6)) \mapsto (t_2^{-1} y_1, t_2y_1, t_1^{-1} y_3, t_1y_3, t_1^{-1} t_2y_5, t_1t_2^{-1} y_6).
		\]
		Computing the determinant of the Jacobian matrix of $\Phi$ we obtain
		\[
		\det(J_\Phi)=4y_1y_3t_1^{-1}t_2^{-1}.
		\]
		In particular, at the point $(e,s')$ we have
		\[
		\det(J_\Phi(e,s'))=4ab.
		\]
		By the generality assumption, we have that this determinant is nonzero. Therefore $\Phi$ is étale at $(e, s')$.
		
		The last statement on the completion of the local ring $\cO_{S',s'}$ follows immediately from the isomorphism $S' \simeq \AA^4$.
	\end{proof}

	\subsection{Slice for the x-variables}  \label{sec:slice_x}
	
	In order to define and analyse the slice for the $\mathbf{x}$-variables and to study the étaleness of the action of the torus, we present the following auxiliary result, which we will subsequently apply.
	
	\begin{lemma} \label{lem:etale}
		Let $S''$ be an affine scheme of finite type over $\CC$ equipped with an algebraic action of the torus $\Gm$. Let $u, v \in \cO_{S^{''}}(S'')$ be two homogeneous regular functions of weights $-1$ and $1$, respectively. Let $S' \hookrightarrow S''$ be the closed subscheme defined by the ideal $(u - v) \subseteq \cO_{S^{''}}(S'')$, and let $U\subseteq S''$ be the principal open subset corresponding to the localisation of $u$. Then the restricted action morphism
		\[
		\Gm \times (S' \cap U) \longrightarrow U, \quad (t, s') \mapsto t \cdot s'
		\]
		is étale.
	\end{lemma}
	\begin{proof}
		Let $A = \cO_{S^{''}}(S'')$ and let $\Gm = \Spec \CC[t^{\pm 1}]$. The algebraic action of $\Gm$ on the affine scheme $S''$ corresponds to the morphism
		\[
		\alpha \colon A \longrightarrow A[t^{\pm1}].
		\]
		This action induces a $\ZZ$-grading on $A$: if $a \in A$ is homogeneous of weight $d$, then the action is given by $\alpha(a)=at^d$. We have:
		\[ 
		\alpha(u) = u t\inv, \quad \alpha(v) = v t. 
		\]
		
		Since the closed subscheme $S' \hookrightarrow S''$ is defined by the ideal $(u - v) \subseteq A$, we have that its coordinate ring is $A' = \cO_{S'}(S') = A / (u - v)$.
		For each $a\in A$ denote its image in $A'=A/(u-v)$ by $\bar{a}$.
		
		Since $U$ is the principal open subset corresponding to the localisation of $u$, we have
		\[
		U=\operatorname{Spec} A[u^{-1}],
		\qquad
		S'\cap U=\operatorname{Spec} A'[\bar u^{-1}].
		\]
		The restriction of the action morphism
		\[
		\Phi\colon \Gm\times(S'\cap U)\longrightarrow U
		\]
		corresponds to the homomorphism
		\[
		\Phi^*\colon A[u^{-1}] \longrightarrow A'[\bar u^{-1},t^{\pm1}],
		\]
		given by $\Phi^*(a)=\bar at^d$, if $a\in A$ has weight $d$.
		
		Set
		\[
		B:=  \frac{A[u\inv, T^{\pm 1}]}{(T^2 - vu^{-1})}.
		\]
		Observe that $\Phi^*$ naturally extends to an $A[u\inv]$-algebra homomorphism
		\[
		\varphi \colon B \longrightarrow A'[\bar{u}\inv, t^{\pm 1}]
		\]
		by setting $\varphi(T) = t$. 
		This map is well defined because $\bar{u}=\bar{v}$ in $A'[\bar{u}\inv]$, and hence
		\[
		\varphi(vu^{-1}) = \Phi^*(v)\Phi^*(u)^{-1} = \frac{\bar{v}t}{\bar{u}t\inv} = t^2 = \varphi(T^2).
		\]
		We can define a map
		\[
		\psi \colon A'[\bar{u}\inv, t^{\pm 1}] \longrightarrow B
		\]
		by setting $\psi(t) = T$ and, for any homogeneous element $a \in A[u\inv]$ of weight $d \in \ZZ$, defining $\psi(\bar{a}) = a T^{-d}$. To see that $\psi$ is well-defined, we just need to check the relation $\bar{u}- \bar{v}$:
		\[
		\psi(\bar{u} - \bar{v}) = uT - vT\inv = T\inv(uT^2 - v) = 0 \quad \text{in } B.
		\]
		It is straightforward to check that $\varphi$ and $\psi$ are mutually inverse. Hence we have the $A[u^{-1}]$-algebra isomorphism
		\[
		A'[\bar{u}\inv, t^{\pm 1}] \simeq B= \frac{A[u\inv, T, T\inv]}{(T^2 - vu^{-1})} \simeq \left( \frac{A[u\inv][T]}{(T^2 - vu^{-1})} \right)_T
		\]
		The polynomial $P(T)=T^2- vu^{-1}$ has derivative $P'(T)=2T,$ which is invertible in $B$.
		By the standard étaleness criterion, it follows that
		\[
		A[u^{-1}]\longrightarrow B
		\]
		is étale, so the corresponding morphism
		\[
		\Phi\colon
		\Gm\times(S'\cap U)\longrightarrow U
		\]
		is étale (see \cite[TAG 02GT]{stacks-project}).
	\end{proof}
	
	\begin{lemma}[Slice for the \textbf{x}-variables] \label{lem:slice_x}
		Let $I$ be the ideal in the polynomial ring $\CC [x_{1,+}, \dots x_{5,+}, x_{1,-}, \dots, x_{5,-}]$ generated by the polynomials $f_{\epsi, k}$ and $g_{\epsi, k}$ of Theorem~\ref{thm:tor_var_X}(iv).
		Let $S$ be the closed subscheme of 
		\[
		\AA^{10}= \Spec \CC [x_{1,+}, \dots x_{5,+}, x_{1,-}, \dots, x_{5,-}]
		\]
		defined by $I$.
		Let the torus $G:=\Gm$ act diagonally with weights $1$ and $-1$ on $\{ x_{1,+}, \dots x_{5,+} \}$ and $\{ x_{1,-}, \dots x_{5,-} \}$, respectively. Let $S' \into S$ be the closed affine subscheme defined by the ideal
		\[
		\Big(
		x_{1,+} + x_{3,+} - 2x_{5,+},\;
		x_{2,+} + x_{4,+},\;
		(x_{3,+}-x_{5,+})^2,\;
		x_{4,+}^2, \;
		x_{2,-}-x_{5,-}, \;
		x_{3,-}, \;
		x_{5,+}-x_{5,-}
		\Big)
		\]
		of $\CC [x_{1,+}, \dots, x_{5,+}, x_{1,-}, \dots, x_{5,-}]$. Let $s'\in S'(\CC)$ be a general point.
		
		Then the morphism $G \times S' \to S$  given by the restriction of the $G$-action is \'etale at the point $(e,s')$ and there is an isomorphism of complete local $\CC$-algebras
		\[
		\widehat{\mathcal{O}_{S', s'}} \simeq \frac{\CC \pow{w_1, \ldots, w_5}}{\big(w_1^2, \; w_2^2\big)}.
		\]
	\end{lemma}
	\begin{proof}
		We have $s'=(u,0,u,0,u,e,u,0,f,u)\in S'(\CC) \subseteq S(\CC) \subseteq \AA^{10}$, with general $u,e,f \in \CC$. 
		We start by proving the étaleness at the point $(e, s')$ of the map $G \times S' \to S$ obtained by restricting the $G$-action on $S$ to $S'$.
		Let $S''$ be the non-reduced closed subscheme of $S$ defined by the ideal
		\[
		\begin{aligned}
			\bigl(x_{1,+}+x_{3,+}-2x_{5,+},\;
			x_{2,+}+x_{4,+},\;
			(x_{3,+}-x_{5,+})^2,\;
			x_{4,+}^2,\;
			x_{2,-}-x_{5,-},\;
			x_{3,-}\bigr).
		\end{aligned}
		\]
		We have that $S''$ is non-reduced and its underlying topological space is the irreducible component of $S$ corresponding to the minimal prime ideal $\mathfrak{P}_{1,2}$ (see Theorem~\ref{thm:tor_var_X}(v) and \cite[Lemma 3.2]{Pet22a}).
		Note that $S'$ is the closed subscheme of $S''$ obtained by imposing the additional equation $x_{5,+} = x_{5,-}$. 
		Since the point $s'$ is away from the intersection of $S''$ with the union of any other irreducible component of $S$, we have that locally around $s'$, the inclusion $S'' \hookrightarrow S$ is an open immersion, which is in particular étale.
		Applying Lemma~\ref{lem:etale} to $S' \hookrightarrow S''$ (with $u=x_{5,-}$ and $v=x_{5,+}$), the restricted action morphism $\Gm \times S' \to S''$ is étale over the open subset where $x_{5,+} \neq 0$. Since the composition of étale morphisms is étale, we conclude that the composed map
		\[
		\Gm \times S' \longrightarrow S'' \hookrightarrow S
		\]
		is étale in a neighborhood of $s'$.
		
		Finally, to compute the local ring $\cO_{S',s'}$, we make the following change of coordinates to move the point $s'$ to the origin: 
		\[
		\begin{aligned} 
			\varepsilon_{1,+} &= x_{1,+} - u, & \quad \varepsilon_{1,-} &= x_{1,-} - e, \\ \varepsilon_{2,+} &= x_{2,+}, & \quad \varepsilon_{2,-} &= x_{2,-} - u, \\ \varepsilon_{3,+} &= x_{3,+} - u, & \quad \varepsilon_{3,-} &= x_{3,-}, \\ \varepsilon_{4,+} &= x_{4,+}, & \quad \varepsilon_{4,-} &= x_{4,-} - f, \\ \varepsilon_{5,+} &= x_{5,+} - u, & \quad \varepsilon_{5,-} &= x_{5,-} - u .
		\end{aligned}
		\]
		In this new setting the equations that define $S'$ in $\AA^{10}$ become
		\begin{align*}
			x_{1,+} + x_{3,+} - 2x_{5,+} &= \varepsilon_{1,+} + \varepsilon_{3,+} - 2\varepsilon_{5,+}, \\
			x_{2,+} + x_{4,+} &= \varepsilon_{2,+} + \varepsilon_{4,+}, \\
			(x_{3,+} - x_{5,+})^2 &=  (\varepsilon_{3,+} - \varepsilon_{5,+})^2 ,\\
			x_{4,+}^2 &=  \varepsilon_{4,+}^2, \\
			x_{2,-} - x_{5,-} &=  \varepsilon_{2,-} - \varepsilon_{5,-}, \\
			x_{3,-} &=  \varepsilon_{3,-},\\
			x_{5,+} - x_{5,-} &=  \varepsilon_{5,+} - \varepsilon_{5,-}.
		\end{align*}
		We compute that
		\[
		\mathcal{O}_{S', s'} \simeq \left( \frac{\CC[\varepsilon_{2,+}, \varepsilon_{3,+}, \varepsilon_{5,+}, \varepsilon_{1,-}, \varepsilon_{4,-}]}{\big(\varepsilon_{2,+}^2, \; (\varepsilon_{3,+} - \varepsilon_{5,+})^2\big)} \right)_{(\varepsilon_{2,+}, \varepsilon_{3,+}, \varepsilon_{5,+}, \varepsilon_{1,-}, \varepsilon_{4,-})}.
		\]
		Setting $\eta= \varepsilon_{3,+} - \varepsilon_{5,+}$, taking the completion, and giving new names to the variables we conclude the proof.
	\end{proof}
	
	\subsection{Combining the slices} \label{sec:proof}
	
	To prove Theorem~\ref{thm:result_n=3} we need to combine the slices of Lemma~\ref{lem:slice_y} and Lemma~\ref{lem:slice_x}.
	
	\begin{proof} [Proof of Theorem~\ref{thm:result_n=3}]
		Let $\stack{3,11/3}$ be the stack of K-semistable Fano varieties of dimension $3$ with anticanonical volume $11/3$, which is an algebraic stack of finite type over $\CC$ with affine diagonal and admitting a good moduli space \cite{Xu_survey}. 
		Let $p = [X] \in \stack{3,11/3}(\CC)$ be the $\CC$-point representing the toric Fano variety $X$ studied in Theorem~\ref{thm:tor_var_X}. 
		Since the variety $X$ is K-polystable by Theorem~\ref{thm:tor_var_X}(i), the point $p$ is closed in $\stack{3, 11/3}$  by \cite[Proof of Theorem 1.3]{openness_of_K-semistab}. 
		
		We denote by $S_\mathbf{x}$ and $S'_\mathbf{x}$ the schemes $S$ and
		$S'$ constructed in Lemma~\ref{lem:slice_x}, and by $S_\mathbf{y}$ and
		$S'_\mathbf{y}$ the schemes $S$ and $S'$ constructed in
		Lemma~\ref{lem:slice_y}. We then set
		\[
		S:=S_\mathbf{y}\times S_\mathbf{x},
		\qquad
		S':=S'_\mathbf{y}\times S'_\mathbf{x} \hookrightarrow S.
		\]
		Let $s'=(s'_\mathbf{y},s'_\mathbf{x})\in S'(\CC)$ be a general point on $S'$. The coordinates are
		\[
		s'_\mathbf{y}=(a,a,b,b,c,d),\qquad
		s'_\mathbf{x}=(u,0,u,0,u,e,u,0,f,u),
		\]
		with general $a,b,c,d,e,f,u\in\CC$.
		
		Consider the automorphism group $G = \Aut(X)$, which is a reductive algebraic group isomorphic to the semidirect product $T_N \rtimes (\ZZ/2\ZZ)^3$ by Theorem~\ref{thm:tor_var_X}(iii).
		Let $G$ act on $S$ in the same way as it acts on $\Spf R$ (see Theorem~\ref{thm:tor_var_X}(vi)), with $R$ the hull of $\DefqG{X}$ (Theorem~\ref{thm:tor_var_X}(iv)).
		
		In particular, the origin $0 \in S$ is fixed by $G$ and we have a $G$-equivariant isomorphism $\widehat{\cO_{S, 0}} \simeq R$.
		Thus we have an isomorphism of stacky germs (see Appendix~\ref{sec:stack})
		\[
		(\stack{3,11/3}, p) \simeq (G \curvearrowright \widehat{\cO_{S,0}}).
		\]
		
		In order to apply Theorem~\ref{thm:artin_stack}, we check its hypotheses.
		
		\medskip
		
		(1) 
		The torus weights corresponding to the non-zero coordinates of $s'$ are
		\[
		\pm(0,1,0),\quad \pm(1,0,0),\quad
		\pm(1,-1,0),\quad \pm(0,0,1).
		\]
		Since the origin lies in the interior of the convex hull of these weights, the torus orbit $T_N\cdot s'$ is closed in $\AA^{16}$ by the Hilbert--Mumford criterion \cite[Theorem 9.2]{Dolgachev}.
		
		Since $G=T_N \rtimes \Aut(P)$ and $\Aut(P)$ is finite, we have that the orbit
		\[G\cdot s'
		=
		\bigcup_{h\in\Aut(P)}h\cdot(T_N\cdot s').
		\]
		is the union of finitely many closed torus orbits in $\AA^{16}$. It follows that the orbit $G\cdot s'$ is closed in $\AA^{16}$ and hence also in $S$.
		
		\medskip
		
		(2)
		Let $G'$ be the stabiliser of $s'$. We want to prove that $G'$ is trivial.
		An element of the group $G$ can be written as 
		\[
		g= t M_{1}^{i}M_{2}^{j}M_{3}^{k}
		\]
		with $t=(t_1, t_2, t_3)\in \Gm^3$ and $i,j, k\in \{0,1\}$.
		To determine the stabiliser $G'$ of $s'$, we compute the elements $g\in G$ satisfying $g \cdot s'=s'$.
		Since $M_1$ and the third torus factor act trivially on the $\mathbf y$-variables, the condition
		\[
		g\cdot s'_\mathbf{y}=s'_\mathbf{y}
		\]
		depends only on $t_1,t_2$ and on $M_2^jM_3^k$.
		
		We carry out the computation explicitly in the case $j=k=1$:
		\begin{align*}
			t M_{2} M_{3} \cdot s'_\mathbf{y} = (t_2^{-1} b, t_{2}b, t_{1}^{-1}a, t_{1}a, t_{1}^{-1}t_{2}d, t_{1}t_{2}^{-1}c).
		\end{align*}
		If $t M_{2} M_{3} \cdot s'_\mathbf{y} = s'_\mathbf{y}$, then comparing the first two coordinates we have
		\[
		t_{2}b=a, \quad t_{2}^{-1}b=a.
		\]
		Hence $t_{2}^{2}=1$ and $t_{2}=\pm 1$. If $t_{2}=1$, then the first equation becomes $a=b$, whereas if $t_{2}=-1$, it becomes $a=-b$. Since $s'$ is a general point of $S'$, neither equality holds. Therefore, there is no element of the form $t M_{2} M_{3}$ fixing $s'_\mathbf{y}$.
		The remaining cases $(j=k=0)$, $(j=1,k=0)$, $(j=0,k=1)$ are analogous and we omit them.
		Thus any element stabilising $s'_\mathbf{y}$ is of the form
		\[
		tM_1^i \textrm{ with } t=(1,1,t_3)
		\]
		It remains to consider the action on the $s'_\mathbf{x}$.
		Firstly, consider the case $i=0$:
		\begin{align*}
			t \cdot s'_\mathbf{x} = (tu,0,tu,0,tu,t^{-1}e,t^{-1}u,0,t^{-1}f,t^{-1}u).
		\end{align*}
		If $t\cdot s'_\mathbf{x} = s'_\mathbf{x}$, then comparing the first coordinate we have $tu=u$, hence $t=1$ since $u\neq 0$ for a general point in $S'$.
		If we now consider $i=1$, we have
		\begin{align*}
			tM_1 \cdot s'_\mathbf{x} =
			(te,tu,0,tf,tu,t^{-1}u,0,t^{-1}u,0,t^{-1}u).
		\end{align*}
		This is equal to $s'_\mathbf{x}$ if $u=0$ (check the third component), which is not true for a general point $s'\in S'$. 
		We can conclude that the stabiliser $G'$ of $s'=(s'_\mathbf{y},s'_\mathbf{x})$ is trivial.
		
		\medskip
		
		(3)
		Since $G'$ is trivial, we have $G\times^{G'}S'=G\times S'$. To verify the étaleness condition (3) in Theorem~\ref{thm:artin_stack}, it is enough to consider the identity component of $G$, which is isomorphic to
		\[ 
		\Gm^3=\Gm^2\times\Gm. 
		\]
		The restricted action morphism $\Phi \colon \Gm^3 \times S' \to S$ can be written as the product $\Phi_\mathbf{y} \times \Phi_\mathbf{x}$, where
		\[
		\Phi_\mathbf{y} \colon \Gm^2 \times S'_\mathbf{y} \longrightarrow S_\mathbf{y},
		\]
		\[
		\Phi_\mathbf{x} \colon \Gm \times S'_\mathbf{x} \longrightarrow S_\mathbf{x}.
		\]
		Since
		$\Phi_\mathbf{y}$ is étale at $(e,s'_\mathbf{y})$ by Lemma~\ref{lem:slice_y}
		and
		$\Phi_\mathbf{x}$ is étale at $(e,s'_\mathbf{x})$ by
		Lemma~\ref{lem:slice_x},
		their product $\Phi = \Phi_\mathbf{y} \times \Phi_\mathbf{x}$
		is étale at $(e,s')$.
		\medskip
		
		(4)
		By Lemma~\ref{lem:slice_y} and Lemma~\ref{lem:slice_x}, it is easy to compute that 
		\[
		\widehat{\mathcal{O}_{S', s'}} \simeq \mathbb{C}\pow{y_1, y_3, y_5, y_6} \widehat{\otimes}_{\mathbb{C}} \frac{\mathbb{C} \pow{w_1, \ldots, w_5}}{\big(w_1^2, \; w_2^2\big)} .
		\]
		Defining this ring to be $R'$, we obtain the desired isomorphism of stacky germs
		$(G' \curvearrowright \widehat{\mathcal{O}_{S',s'}}) \simeq (G' \curvearrowright R')$
		where the trivial group $G'$ acts trivially.
		
		Hence, by Theorem~\ref{thm:artin_stack} there exists a closed $\CC$-point $p' \in \stack{3, 11/3}$ corresponding to a Fano 3-fold $Y$ with miniversal ring $R'$ and trivial automorphism group. Since $p'$ is a closed point of the stack of K-semistable Fano varieties, $Y$ is K-polystable by \cite[Proof of Theorem 1.3]{openness_of_K-semistab}. This fact, together with the finiteness of its automorphism group, implies that $Y$ is K-stable.
		This concludes the proof of Theorem~\ref{thm:result_n=3}.
	\end{proof}
	
	\begin{remark}
		If the base field $\CC$ is the field of complex numbers, the Fano 3-fold $Y$ constructed in the proof of Theorem~\ref{thm:result_n=3} has at least 2 isolated singularities analytically isomorphic to the anticanonical affine cone over the first Hirzebruch surface by Proposition~\ref{prop:torvar_cone_over_octagon}(iv). We did not study the other singularities of $Y$, which must exist because $Y$ is not Gorenstein as its anticanonical volume is not an integer.
	\end{remark}
	
	\section{Final proofs} \label{sec:final_proofs}
	
	The goal of this final section is to give the proofs of Theorem~\ref{thm:main} and Corollary~\ref{cor:main}. To this end, we first introduce a key lemma.
	
	\begin{lemma} \label{lem:automorfismi_prodotto_Fano}
		Let $X_1$ and $X_2$ be two non-isomorphic klt Fano varieties over $\CC$.
		For each $i \in \{1,2\}$, assume that $X_i$ is not isomorphic to the product of varieties of smaller dimension.
		Then the group homomorphism
		\begin{equation} \label{eq:aut_product}
			\Aut(X_1) \times \Aut(X_2) \longrightarrow \Aut(X_1 \times X_2)
		\end{equation}
		given by $(f_1, f_2) \mapsto f_1 \times f_2$
		is an isomorphism.
		The same statement holds for the corresponding
		automorphism group schemes.
	\end{lemma}
	
	\begin{proof}
		If one of $X_1, X_2$ is isomorphic to $\Spec \CC$, the assertion is immediate. We may therefore assume $\dim X_1 >0$ and $\dim X_2 > 0$.
		
		We begin proving the statement for the automorphism groups as abstract groups.
		It is obvious that the homomorphism \eqref{eq:aut_product} is injective.
		Thus we need to prove that it is surjective.
		Set $X=X_1\times X_2$.
		Fix $f \in \Aut(X)$. We have to show that $f$ lies in the image of \eqref{eq:aut_product}.
		
		Denote the projections by
		$\pi_i \colon X\to X_i$ for $i =1,2$.
		Since $\rH^1(X_1,\cO_{X_1})=0$ by Kawamata--Viehweg vanishing and since $X_2$ is connected,
		by \cite[Exercise~III.12.6(b)]{hartshorne} the group homomorphism
		\begin{equation} \label{eq:isom_Pic}
			\Pic(X_1) \times \Pic(X_2) \longrightarrow \Pic(X)
		\end{equation}
		given by $(L_1,L_2) \mapsto \pi_1^* L_1 \otimes_{\cO_X} \pi_2^*L_2$ is an isomorphism.
		By tensoring \eqref{eq:isom_Pic} with $\RR$ and by taking numerical equivalence, we obtain that the natural $\RR$-linear map
		\begin{equation} \label{eq:isom_N^1}
			\rN^1(X_1)_\RR \times \rN^1(X_2)_\RR \longrightarrow \rN^1(X)_\RR
		\end{equation}
		is surjective.
		By intersecting with curves on $X$ of the type $C_1 \times \{ x_2 \}$ or $\{ x_1 \} \times C_2$, where $C_i$ is a curve in $X_i$ and $x_i$ is a closed point of $X_i$, it is easy to see that \eqref{eq:isom_N^1} is injective, and hence an isomorphism.
		We always use this isomorphism and thus we always consider the decomposition
		\begin{equation*}
			\rN^1(X)_{\RR}=
			\rN^1(X_1)_{\RR} \times
			\rN^1(X_2)_{\RR}.
		\end{equation*}
		Using the same types of curves as above, we obtain the equality of nef cones
		\begin{equation*}
			\Nef(X) = \Nef(X_1) \times \Nef(X_2).
		\end{equation*}
		By dualising we obtain the decomposition
		\begin{equation*}
			\rN_1(X)_{\RR}=
			\rN_1(X_1)_{\RR} \times
			\rN_1(X_2)_{\RR}
		\end{equation*}
		and the equality of Mori cones
		\begin{equation} \label{eq:NE_product}
			\NE(X) = \NE(X_1) \times \NE(X_2).
		\end{equation}
		See also \cite[Lemma~6.1]{fanelli_floris_zimmermann}.
		Since $X_1, X_2, X$ are klt Fano varieties,
		by the cone theorem \cite[Theorem~3.7]{kollar_mori}
		the cones
		$\NE(X_1), \NE(X_2), \NE(X)$ are strongly convex (i.e.\ pointed), rational, polyhedral and generated by classes of rational curves.
		
		\begin{claim} \label{claim:existence_f_1}
			There exists $f_1 \in \Aut(X_1)$ such that $\pi_1 \circ f = f_1 \circ \pi_1$.
		\end{claim}
		
		We now prove Claim~\ref{claim:existence_f_1}.
		Consider the contraction $q := \pi_1 \circ f \colon X \to X_1$ and the cone
		\[
		F := \ker q_* \cap \NE(X)
		\]
		which is a face of $\NE(X)$ by \cite[Proposition~1.14(a)]{debarre}.
		By \eqref{eq:NE_product} there exist two faces $F_1, F_2$ of the Mori cones $\NE(X_1), \NE(X_2)$, respectively, such that $F = F_1 \times F_2$.
		For each $i \in \{1,2\}$, let $c_i \colon X_i \to Y_i$ be the contraction associated to $F_i$ by \cite[Theorem~3.7(3)]{kollar_mori}.
		From $F_i = \ker (c_i)_* \cap \NE(X_i)$ we deduce that the contractions $c_1 \times c_2 \colon X \to Y_1 \times Y_2$ and $q \colon X \to X_1$ contract exactly the same curve classes.
		By the uniqueness of the contraction associated to a face of the Mori cone (see \cite[Proposition~1.14(b)]{debarre}), there exists an isomorphism $g \colon Y_1 \times Y_2 \to X_1$ such that $q = g \circ (c_1 \times c_2)$.
		We have a commutative diagram
		\begin{equation*}
			\xymatrix{
				X \ar[r]^f \ar[d]_{c_1 \times c_2} \ar[rd]^q & X \ar[d]^{\pi_1} \\
				Y_1 \times Y_2 \ar[r]_g & X_1
			}
		\end{equation*}
		where the horizontal arrows are isomorphisms and the other arrows are contractions.
		In particular $X_1$ is isomorphic to $Y_1 \times Y_2$.
		Since $X_1$ is not isomorphic to the product of varieties of smaller dimension, we have that at least one among $Y_1$ and $Y_2$ is isomorphic to $\Spec \CC$.
		
		We now show that $Y_1$ cannot be a point.
		For a contradiction, assume that $Y_1$ is a point.
		Hence the projection $p_2 \colon Y_1 \times Y_2 \to Y_2$ is an isomorphism.
		We have a commutative diagram
		\begin{equation*}
			\xymatrix{
				X_2 \ar[d]_{c_2} & X \ar[r]^f \ar[d]_{c_1 \times c_2} \ar[rd]^q \ar[l]_{\pi_2} & X \ar[d]^{\pi_1} \\
				Y_2 & Y_1 \times Y_2 \ar[r]_g \ar[l]^{p_2} & X_1
			}
		\end{equation*}
		where the bottom horizontal arrows are isomorphisms.
		Fix a closed point $y_2 \in Y_2$, set $x_1 := g ( (p_2)^{-1}(y_2)  ) \in X_1$ and observe that $f$ and $g \circ (p_2)^{-1}$ induce an isomorphism between the fibre $\pi_1^{-1}(x_1) = \{ x_1\} \times X_2 \simeq X_2$ and the fibre $(c_2 \circ \pi_2)^{-1}(y_2) = X_1 \times (c_2)^{-1}(y_2)$;
		thus we obtain that $X_2$ is isomorphic to $X_1 \times (c_2)^{-1}(y_2)$.
		Since $X_2$ is integral, we deduce that $(c_2)^{-1}(y_2)$ is integral.
		Since $X_2$ is not isomorphic to the product of varieties of smaller dimension and since $\dim X_1 > 0$, we obtain that $(c_2)^{-1}(y_2)$ is isomorphic to $\Spec \CC$.
		Hence $X_2$ is isomorphic to $X_1$, which is a contradiction.
		This implies that it was absurd to assume that $Y_1$ is a point.
		
		Since $Y_1$ cannot be a point, we have that $Y_2$ must be a point.
		Hence the projection $p_1 \colon Y_1 \times Y_2 \to Y_1$ is an isomorphism.
		We have a commutative diagram
		\begin{equation*}
			\xymatrix{
				X_1 \ar[d]_{c_1} & X \ar[r]^f \ar[d]_{c_1 \times c_2} \ar[rd]^q \ar[l]_{\pi_1} & X \ar[d]^{\pi_1} \\
				Y_1 & Y_1 \times Y_2 \ar[r]_g \ar[l]^{p_1} & X_1
			}
		\end{equation*}
		where the bottom horizontal arrows are isomorphisms.
		For every closed point $y_1 \in Y_1$, set $x_1 := g ( (p_1)^{-1}(y_1)  ) \in X_1$ and observe that $f$ and $g \circ (p_1)^{-1}$ induce an isomorphism between the fibre $\pi_1^{-1}(x_1) = \{ x_1\} \times X_2 \simeq X_2$ and the fibre $(c_1 \circ \pi_1)^{-1}(y_1) =  (c_1)^{-1}(y_1) \times X_2$; thus
		$X_2$ is isomorphic to $(c_1)^{-1}(y_1) \times X_2$.
		Hence $(c_1)^{-1}(y_1)$ is isomorphic to $\Spec \CC$ for every closed point $y_1 \in Y_1$.
		In particular, $c_1$ is quasi-finite and proper, hence it is finite.
		Since $c_1$ is finite and $(c_1)_* \cO_{X_1}= \cO_{Y_1}$ because $c_1$ is a contraction, we deduce that $c_1 \colon X_1 \to Y_1$ is an isomorphism.
		The automorphism $f_1 := g \circ (p_1)^{-1} \circ c_1$ of $X_1$ satisfies $\pi_1 \circ f = f_1 \circ \pi_1$.
		This concludes the proof of Claim~\ref{claim:existence_f_1}.
		
		\medskip
		
		We have proved that there exists $f_1 \in \Aut(X_1)$ such that $\pi_1 \circ f = f_1 \circ \pi_1$.
		By swapping the roles of $X_1$ and $X_2$, we obtain that there exists $f_2 \in \Aut(X_2)$ such that $\pi_2 \circ f = f_2 \circ \pi_2$.
		We deduce that $f = f_1 \times f_2$, hence $f$ lies in the image of \eqref{eq:aut_product}.
		This concludes the proof that the abstract group homomorphism \eqref{eq:aut_product} is an isomorphism.
		
		\medskip
		
		Finally, consider \eqref{eq:aut_product} as a morphism of
		automorphism group schemes. These group schemes are of finite type over $\CC$ because a
		sufficiently large multiple of the anticanonical line bundle is very ample and is preserved
		by every automorphism.
		Since we are in characteristic zero, these group schemes are smooth.
		
		We have the isomorphism of coherent sheaves $\cT^0_X\simeq
		\pi_1^*\cT^0_{X_1}\oplus\pi_2^*\cT^0_{X_2}$. This follows from the corresponding statement for the Kähler differentials, by restricting to the smooth locus, dualising and using the fact that the tangent sheaf of a normal variety is reflexive.
		Therefore we have the isomorphism
		\begin{equation} \label{prodotto_H0}
			\rH^0(X,\cT^0_X)\simeq
			\rH^0(X_1,\cT^0_{X_1})\oplus \rH^0(X_2,\cT^0_{X_2}).
		\end{equation}
		which means that the differential at the identity of the group scheme homomorphism \eqref{eq:aut_product} is bijective.
		This implies that \eqref{eq:aut_product} is étale.
		
		It is clear that \eqref{eq:aut_product} is a monomorphism of schemes. Since it is also étale, it is an open immersion.	
		The surjectivity on closed points proved above forces this open immersion to be
		surjective, and therefore it is an isomorphism of group schemes.
	\end{proof}

	\begin{proof}[Proof of Theorem~\ref{thm:main}]
		For $n=3$, Theorem~\ref{thm:main} immediately follows from Theorem~\ref{thm:result_n=3}. 
		
		Now suppose $n\geq 5$. Let $Z \subset \PP^{n-2}$ be a general hypersurface of degree $n-2$. We have $\dim Z = n-3$.
		By \cite[Corollary~1.5]{Fujita}, $Z$ is K-stable, and by \cite{MatMon}, the group of automorphisms $\Aut(Z)$ is trivial.
		
		Now consider the product $Y \times Z$, where $Y$ is the Fano 3-fold in Theorem~\ref{thm:result_n=3}. We first verify that $Y$ and $Z$ satisfy the hypotheses of Lemma~\ref{lem:automorfismi_prodotto_Fano}. Since $Y$ is singular and $Z$ is smooth, we have that $Y$ and $Z$ are not isomorphic.
		Now we prove that $Y$ is not the product of varieties of smaller dimension. For a contradiction, suppose that $Y\simeq Y_1\times Y_2$, where both factors have positive dimension. Since $Y$ is a Fano 3-fold, both $Y_1$ and $Y_2$ are Fano varieties and, up to exchanging the factors, we have $\dim Y_1=1$ and $\dim Y_2=2$. In particular, $Y_1$ must be isomorphic to $\PP^1$. By \eqref{prodotto_H0}, this would imply $\rH^0(Y,\cT_Y^0)\neq 0$, contradicting the triviality of $\Aut(Y)$, since $\rH^0(Y,\cT_Y^0)$ is its tangent space at the identity. 
		Hence $Y$ cannot be written as a product of varieties of smaller dimension.
		
		Now we prove that $Z$ is not the product of varieties of smaller dimension. If $n\geq 6$, by Lefschetz hyperplane theorem we have that the Picard rank $\rho(Z)$ is $1$. However, if $Z$ were isomorphic to a product of two positive-dimensional varieties, both factors would be Fano and hence have positive Picard number, implying $\rho(Z)\geq 2$, a contradiction.
		If $n=5$, then $Z$ is a smooth cubic surface. A product decomposition would force $Z\simeq\PP^1\times\PP^1$, contradicting the triviality of $\Aut(Z)$. Thus $Z$ is not isomorphic to the product of varieties of smaller dimension.
		
		Since both $\Aut(Y)$ and $\Aut(Z)$ are trivial, by Lemma~\ref{lem:automorfismi_prodotto_Fano} it follows that the automorphism group $\Aut(Y \times Z)$ is trivial. 
		Furthermore, $Y \times Z$ is K-stable by \cite[Theorem~1.1]{Zhuang_product}. 
		Finally, by \cite[Proposition~4.1]{KP}, the base of the miniversal qG-deformation of $Y \times Z$ is the product of the bases of the respective miniversal qG-deformations of $Y$ and $Z$. Since the base of the miniversal qG-deformation of $Y$ is non-reduced, we have that the base for $Y \times Z$ is also non-reduced. This concludes the proof of Theorem~\ref{thm:main}. 
	\end{proof}

	\begin{proof}[Proof of Corollary~\ref{cor:main}]
		Let $Y$ be the $n$-dimensional K-stable Fano variety constructed in Theorem~\ref{thm:main}, and set $v= (-K_Y)^n \in \QQ_{>0}$ to be its anticanonical volume. By construction, the point $[Y]$ belongs to the open substack $\stackstab{n,v} \subseteq \stack{n,v}$.
		
		By the local structure of the K-moduli stack, a formal neighbourhood of $[Y]$ in $ \stack{n,v}$ is isomorphic to the formal quotient stack $[\Spf(R) / \Aut(Y)]$, where $R$ is the hull of $\DefqG{Y}$.
		By Theorem~\ref{thm:main}, we know that $\Aut(Y)$ is trivial and $R$ is non-reduced. So we get that $\stack{n,v}$ is non-reduced at $[Y]$, and since $\stackstab{n,v} \subseteq \stack{n,v}$ is an open embedding, $\stackstab{n,v}$ is also non-reduced at $[Y]$.
		Finally, by the properties of good moduli spaces, a formal neighbourhood of $[Y] \in \modspace{n,v}$ is isomorphic to $\Spf(R^{\Aut(Y)})$. Again, since $\Aut(Y)$ is trivial, $R^{\Aut(Y)} = R$ is non-reduced. Because $\spacestab{n,v} \subseteq \modspace{n,v}$ is an open embedding, we conclude that $\spacestab{n,v}$ is also non-reduced at $[Y]$.
	\end{proof}
	
	\appendix
	
	\section{Points in a neighbourhood on an algebraic stack} \label{sec:stack}
	
	We introduce stacky germs to encode stabiliser actions on formal miniversal deformation spaces of points on algebraic stacks (see Definition~\ref{def:stacky_germ} and Example~\ref{exa:stacky_germ_stack}).
	We establish their invariance under étale morphisms inducing isomorphisms on stabiliser groups (see Lemma~\ref{lem:etale_induce_same_stacky_germ}).
	For algebraic stacks with affine diagonal and good moduli spaces, we use étale slices to give a criterion for realising prescribed stacky germs at closed points (see Theorem~\ref{thm:artin_stack}).
	
	We continue to denote by $\CC$ an arbitrary algebraically closed field of characteristic zero.
	Recall that a linear algebraic group over $\CC$ is reductive if and only if it is linearly reductive.
	
	\begin{definition}\label{def:stacky_germ}
		A \emph{stacky germ} is the datum of
		\begin{itemize}
			\item  a complete local noetherian $\CC$-algebra $R$ with maximal ideal $\frakm_R$ and with residue field $\CC$,
			\item an affine group scheme $G \to \Spec \CC$ of finite type,
			\item an action of $G$ on the formal spectrum of $R$ over $\CC$, i.e.\ a
			sequence of compatible $G$-actions over $\CC$ on the artinian truncations $\Spec R / \frakm_R^{n+1}$ of $R$ for all $n \in \ZZ_{\geq 0}$. The compatibility means that, for every $n \in \ZZ_{\geq 0}$, the closed immersion $\Spec R / \frakm_R^{n+1} \into \Spec R / \frakm_R^{n+2}$, induced by the natural surjection $R / \frakm_R^{n+2} \onto R / \frakm_R^{n+1}$, is $G$-equivariant. 
		\end{itemize}
		Such a stacky germ is denoted by $(G \curvearrowright R)$.
		
		Two stacky germs $(G \curvearrowright R)$ and $(G' \curvearrowright R')$ are called \emph{isomorphic} if there exist an isomorphism $G \simeq G'$ of group schemes over $\CC$ and a $\CC$-algebra isomorphism $R \simeq R'$ such that they induce equivariant isomorphisms on the artinian truncations.
	\end{definition}
	
	\begin{remark}
		In general, 
		the action on the formal spectrum of $R$ need not be induced by an `honest' action of $G$ on the affine scheme $\Spec R$.
		For example, the formal action of $\Gm$ on $\CC \pow{t}$ with weight $1$ does not come from an action because $\CC \pow{t}$ is not the direct sum of the eigenspaces $\CC t^i$. Indeed,	
		it is well known that giving an action of $\Gm$ on a $\CC$-algebra $R$ is equivalent to giving a $\ZZ$-grading on $R$.
		
		Nevertheless, if $G$ is finite then a formal action of $G$ on $R$ is always an action.
	\end{remark}
	
	\begin{remark}
		For an algebraic stack $\cX$, we denote by $| \cX |$ the underlying topological space, as defined in \cite[Tag~04XE]{stacks-project}.
		If $K$ is a field and $x \in \cX(K)$ is a $K$-point of $\cX$, then we often denote by $|x|$ its class in $|\cX|$;
		we say that $x \in \cX(K)$ is closed in $\cX$ if its class $|x|$ is closed in $| \cX |$.
		For a morphism $f \colon \cX \to \cY$ of algebraic stacks, we denote by $|f| \colon | \cX | \to | \cY |$ the induced continuous function.
	\end{remark}
	
	\begin{example} \label{exa:stacky_germ_stack}
		Let $\cM$ be an algebraic stack of finite type over $\CC$ with affine diagonal and admitting a good moduli space.
		Let $p \in \cM(\CC)$ be a closed $\CC$-point of $\cM$.
		Let $R$ be a miniversal ring of $\cM$ at $p$, i.e.\ a hull of the deformation category associated to $\cM$ at $p$ \cite[Tag~0DQT]{stacks-project}.
		Let $G = \Aut_\cM(p)$ be the stabiliser group of $p$, which is (linearly) reductive by \cite[Proposition~12.14]{alper_good}.
		By \cite[Theorem~4.16 and Remark~4.17]{AHR20} one may choose a formal miniversal deformation space $\Spf R$ of $p$ endowed with a formal action of $G$.
		We can consider the stacky germ $(G \curvearrowright R)$.
		Its isomorphism class is well defined by \cite[Theorem~4.19]{AHR20}.
		The stacky germ $(G \curvearrowright R)$ is called the \emph{stacky germ of $\cM$ at $p$} and is denoted by $(\cM, p)$. It should be thought of as the \'etale local structure of $\cM$ at the point $p$ (Lemma~\ref{lem:etale_induce_same_stacky_germ})
		and describes the formal completion $[ \Spf R \ / G]$ along the residual gerbe of a point on an algebraic stack.

		In particular, if $\cM$ is a scheme then $G$ is the trivial group and $R$ is the completion  $\widehat{\cO_{\cM,p}}$ of the local ring $\cO_{\cM,p}$.
	\end{example}
	
	\begin{example} \label{ex:stacky_germ_quotient_stack}
		Let $S$ be an affine scheme of finite type over $\CC$ endowed with an action of the linearly reductive group scheme $G$.
		Let $s \in S(\CC)$ be a closed point fixed by $G$.
		Let $\bar{s} \in [S/G](\CC)$ be the point induced by $s$ on the quotient stack $[S/G]$.
		Then the stacky germ $([S/G], \bar{s})$ of $[S/G]$ at $\bar{s}$ is $(G \curvearrowright \widehat{\cO_{S,s}})$, where the formal action of $G$ on $\widehat{\cO_{S,s}}$ is induced by the algebraic action of $G$ on $S$.
	\end{example}

	\begin{lemma} \label{lem:etale_induce_same_stacky_germ}
		Let $F \colon \cX \to \cY$ be a morphism of algebraic stacks of finite
		type over $\CC$, with affine diagonal and admitting good moduli spaces.
		Let $p \in \cX(\CC)$ and set $q=F(p)$.
		Assume that $p$ is a closed point in $\cX$, that $q$ is a closed point in $\cY$, that $F$ is \'etale at $p$ and that $F$ induces an isomorphism from the stabiliser group of $p$ to the stabiliser group of $q$.
		
		Then there is an isomorphism of stacky germs between
		$
		(\cX,p)$ and $(\cY,q)
		$.
	\end{lemma}
	
	\begin{proof}
		Set
		$
		G := \Aut_{\cX}(p)$.
		By hypothesis the morphism $F$ induces an isomorphism
		$
		\theta \colon G \xrightarrow{\sim}  \Aut_{\cY}(q)
		$.
		Since $p$ and $q$ are closed points and $\cX$ and $\cY$ admit good
		moduli spaces, the group schemes $G$ and $\Aut_{\cY}(q)$ are linearly
		reductive by \cite[Proposition~12.14]{alper_good}.
		
		If we replace $\cX$ by an open substack of $\cX$ containing $p$, we may assume that $F$ is étale.
		
		By \cite[Theorem~1.1]{AHR20} there exist an affine scheme
		$
		W
		$
		endowed with an action of $G$, a $G$-fixed point
		$w\in W(\CC)$, and an \'etale morphism
		$\Phi \colon [W/G] \to \cX$
		such that $\Phi(\bar{w})=p$, where $\bar{w} \in [W/G](\CC)$ is the image of $w$ under the natural $G$-torsor $W \to [W / G]$, and such that the induced homomorphism
		$
		G_w=\Aut_{[W/G]}(\bar{w}) \to G
		$
		is an isomorphism.
		
		Consider the composite
		$
		\Psi := F\circ \Phi\colon [W/G] \to \cY
		$.
		Since both $F$ and $\Phi$ are \'etale, the morphism $\Psi$ is \'etale.
		Moreover,
		$
		\Psi(\bar{w})=F(p)=q
		$.
		The homomorphism on stabiliser groups induced by $\Psi$ at $\bar{w}$ is the
		composite
		\[
		G_w
		\xrightarrow{\sim}
		G
		\xrightarrow{\theta}
		\Aut_{\cY}(q),
		\]
		which is an isomorphism.	
		We have consequently obtained a diagram
		\[
		\xymatrix{
			& [W/G] \ar[ld]_{\Phi} \ar[rd]^{\Psi} & \\
			\cX && \cY
		}
		\]
		of \'etale morphisms such that $\Phi(\bar{w})=p$, $\Psi(\bar{w})=q$, and both $\Phi$ and
		$\Psi$ induce isomorphisms on stabiliser groups at $\bar{w}$.
		
		The preceding diagram verifies condition~(3) of
		\cite[Theorem~4.19]{AHR20}. Hence there exist an
		isomorphism
		$
		G= \Aut_{\cX}(p) \simeq \Aut_{\cY}(q)
		$
		and an equivariant isomorphism
		between formal miniversal deformation spaces
		of $\cX$ and $\cY$ at $p$ and $q$, respectively.
		Therefore we have an isomorphism between the stacky germs $(\cX,p)$ and $(\cY, q)$.
	\end{proof}

	The main result of this appendix is the following theorem which gives information on the existence of closed points in an algebraic stack $\cM$ which are near a closed point at which we know the étale local structure of $\cM$.

	\begin{theorem} \label{thm:artin_stack}
		Let $\cM$ be an algebraic stack of finite type over $\CC$, with affine diagonal and admitting a good moduli space. Let $p \in \cM(\CC)$ be a closed point.
		
		Let $G$ be an algebraic group over $\CC$ acting on an affine scheme $S$ of finite type over $\CC$ and let $s \in S(\CC)$ be a point fixed by $G$.
		Assume that there is an isomorphism of stacky germs
		\begin{equation} \label{eq:isom_stacky_germs}
			(\cM,p) \simeq (G \curvearrowright \widehat{\cO_{S,s}})
		\end{equation}
		where the formal action on $\widehat{\cO_{S,s}}$ is induced by the algebraic action on $S$.

		Let $G'$ be a (linearly) reductive closed subgroup of $G$ and let $(G' \curvearrowright R')$ be a stacky germ.
		Let $S' \into S$ be a closed subscheme which is $G'$-invariant.
		Assume that every open neighbourhood of $s$ in $S$ contains at least one point $s' \in S'(\CC)$ such that:
		\begin{enumerate}
			\item  the orbit $G \cdot s'$ is closed in $S$,
			\item $G'$ is the stabiliser of $s'$,
			\item the natural morphism $G \times^{G'} S' \to S$ is \'etale at the point $[(e,s')]$,
			\item there exists an isomorphism of stacky germs
			\[
			(G' \curvearrowright \widehat{\cO_{S',s'}}) \simeq (G' \curvearrowright R')
			\]
			where the formal action on $\widehat{\cO_{S',s'}}$ is induced by the algebraic action on $S'$.
		\end{enumerate}
		
		Then every open neighbourhood of $p$ in $\cM$ contains a point $p' \in \cM(\CC)$ such that $p'$ is closed in $\cM$ and there is an isomorphism of stacky germs
		\[
		(\cM,p') \simeq (G' \curvearrowright R').
		\]
		In particular, $R'$ is a miniversal ring of $\cM$ at $p'$ and there is a group scheme isomorphism $G' \simeq \Aut_{\cM}(p')$.
	\end{theorem}

	\begin{example} \label{exa:Gm_agisce_pesi_1-1}
		Let $\cM$ be an algebraic stack of finite type over $\CC$, with affine diagonal and admitting a good moduli space. If $p \in \cM(\CC)$ is a closed point such that the stacky germ $(\cM,p)$ is isomorphic to $(\Gm \curvearrowright \CC \pow{t_1,t_2})$ where the formal action is given by the weights $1,-1$, then we can deduce that in every open neighbourhood of $p$ in $\cM$ there exists a closed point $p' \in \cM(\CC)$ such that the automorphism group of $p'$ in $\cM$ is trivial and the miniversal ring of $\cM$ at $p'$ is $\CC \pow{t}$.
		
		Indeed, in Theorem~\ref{thm:artin_stack} we can use the affine plane $S = \Spec \CC [t_1,t_2]$ with the action of $\Gm$ given by the weights $1,-1$, the origin as $s$ and the closed subscheme $S' = \Spec \CC[t_1,t_2]/(t_1-t_2)$ of $S$.
		Every closed point of $S'$ different from the origin lies on a closed orbit and has trivial stabiliser. Moreover, $S'$ serves as a Luna \'etale slice.
	\end{example}
	
	The rest of this section is mainly devoted to the proof of Theorem~\ref{thm:artin_stack}.
	We need some auxiliary lemmas first.
	
	\begin{lemma} \label{lem:fibra_chiusa}
		Let $\cX$ be an algebraic stack over $\CC$, let $X$ be an algebraic space of finite type over $\CC$,
		and let $\pi \colon \cX \to X$ be a morphism over $\CC$.
		Assume that either $X$ is a scheme or $X$ is a quasi-separated algebraic space.
		Let $p \in \cX(\CC)$ and let $x = \pi(p) \in X(\CC)$.
		Let $\cX_x$ be the stack-theoretic fibre over $x$ with respect to $\pi$.
		Denote by $p_x$ the point in $\cX_x(\CC)$ induced by $p$.
		
		Then $|p|$ is closed in $\vert \cX \vert$ if and only if $|p_x|$ is closed in $\vert \cX_x \vert$.
	\end{lemma}
	
	\begin{proof}
		Consider the cartesian diagram
		\begin{equation}\label{eq:diagrr}
			\xymatrix{
				\cX_x \ar[r]^j \ar[d] & \cX \ar[d]^\pi \\
				\Spec \CC \ar[r]^{i} & X
			}
		\end{equation}
		where $i$ is the morphism induced by $x$.
		
		We now prove that $i$ is a closed immersion.
		If $X$ is a scheme, this is well known because $X$ is of finite type over $\CC$.
		If $X$ is a quasi-separated algebraic space, then $X$ is a decent algebraic space by \cite[Tag~03I7]{stacks-project}, hence $i$ is a closed immersion by \cite[Lemma~08AL]{stacks-project}.
		
		Since closed immersions are stable by base change, $j$ is a closed immersion.
		Therefore $|j| \colon |\cX_x| \to |\cX|$ is a closed immersion of topological spaces. We now prove that the image of $|j|$ is the topological fibre $|\pi|^{-1}(|x|)$.
		
		By applying the functor $| \cdot |$ to \eqref{eq:diagrr} we obtain a commutative square of topological spaces
		\begin{equation*}
			\xymatrix{
				|\cX_x| \ar[r]^{|j|} \ar[d] & |\cX| \ar[d]^{|\pi|} \\
				|\Spec \CC| \ar[r]^{|i|} & |X|
			}
		\end{equation*}
		where the horizontal arrows are closed immersions of topological spaces. It follows that $|j| (| \cX_x |) \subseteq | \pi |^{-1}(|x|)$.
		We now need to prove the containment $|j| (|\cX_x|) \supseteq |\pi|^{-1}(|x|)$.
		
		Let $z \in |\pi|^{-1}(|x|)$.
		Then choose an algebraically closed field extension $L \supseteq \CC$ and an object $q \in \cX(L)$ whose class in $|\cX|$ is $z$.
		Thus $q$ determines a morphism $q \colon \Spec L \to \cX$. The class in $|X|$ of the composite $\pi \circ q \colon \Spec L \to X$ is $|x|$.
		The morphism $\pi \circ q$ factors through $i$.
		This gives an $L$-point of $\cX_x$, whose class in $| \cX_x |$ is mapped to $z=|q|$ by $|j|$.
		
		We have thus proved that $|j|$ is a closed immersion from $|\cX_x|$ onto the topological fibre $| \pi |^{-1}(|x|)$.
		The required statement on $|p|$ follows from an obvious argument with the subspace topology.
	\end{proof}
	
	\begin{remark}
		If $X$ is allowed to be a non-quasi-separated algebraic space, then Lemma~\ref{lem:fibra_chiusa} becomes false.
		A counterexample is obtained as follows.
		Consider the quotient $X = \AA^1_\CC / \ZZ$, where the infinite cyclic group $\ZZ$ acts on $\AA^1_\CC$ by translation.
		The $\ZZ$-orbit of the origin in $\AA^1_\CC$ is not closed, hence this induces a non-closed point $x \in X(\CC)$. Set $\pi \colon \cX \to X$ to be the identity and $p \in \cX(\CC)$ such that $\pi(p) = x$.
		Hence the point $p$ is not closed in $\cX$, but the point $p_x$ is closed in the fibre $\cX_x \simeq \Spec \CC$.
	\end{remark}

	\begin{lemma} \label{lem:chiusura_punti_morfismo_da_spazi_algebrici}
		Let 
		\begin{equation} \label{eq:quadrato_cartesiano_stack}
			\xymatrix{
				\cX \ar[d]_\pi \ar[r]^F & \cY \ar[d]^{\pi'} \\
				X \ar[r]^f & Y
			}
		\end{equation}
		be a cartesian diagram of algebraic stacks over $\CC$.
		Assume that $X$ and $Y$ are algebraic spaces of finite type over $\CC$.
		Assume that $X$ and $Y$ are either schemes or quasi-separated algebraic spaces.
		Let $p \in \cX(\CC)$ and $q = F(p) \in \cY(\CC)$.
		
		Then $|p|$ is closed in $\vert \cX \vert$ if and only if $|q|$ is closed in $\vert \cY \vert$.
	\end{lemma}
	
	\begin{proof}
		Let $x := \pi(p) \in X(\CC)$ and $y := \pi' (q) = f(x) \in Y(\CC)$.
		Let $\cX_x$ (resp.\ $\cY_y$) be the stack-theoretic fibre of $\pi$ (resp. $\pi'$) over $x$ (resp.\ $y$).
		Let $p_x \in \cX_x(\CC)$ (resp.\ $q_y \in \cY_y(\CC)$) be the point induced by $p$ (resp.\ $q$).
		By Lemma~\ref{lem:fibra_chiusa} it is enough to show that $|p_x|$ is closed in $|\cX_x|$ if and only if $|q_y|$ is closed in $|\cY_y|$.
		
		Since $x$ and $y$ have the same residue field $\CC$ and the square in \eqref{eq:quadrato_cartesiano_stack} is cartesian,
		we have an isomorphism $\cX_x \simeq \cY_y$
		which maps $p_x$ to $q_y$.
	\end{proof}
	
	\begin{lemma} \label{lemma:artin_local_structure_stacks_with_good_moduli}
		Let $\cX$ and $\cY$ be two algebraic stacks of finite type over $\CC$, with affine diagonal and admitting good moduli spaces.
		Let $p \in \cX(\CC)$ and $q \in \cY(\CC)$ be closed points.
		Set $G = \Aut_\cX(p)$. 
		
		If the stacky germs $(\cX,p)$ and $(\cY,q)$ are isomorphic,
		then there exist an affine scheme $W$ of finite type over $\CC$ equipped with an action of $G$, a point $w \in W(\CC)$ fixed by $G$ and a commutative diagram
		\begin{equation*}
			\xymatrix{
				\cX \ar[d] & [W/G] \ar[d] \ar[r] \ar[l] & \cY \ar[d] \\
				X & W/G \ar[r] \ar[l] & Y
			}
		\end{equation*}
		where:
		\begin{itemize}
			\item $W/G = \Spec \Gamma(W,\cO_W)^G$ is an affine scheme of finite type over $\CC$;
			\item the two squares are cartesian;
			\item the vertical morphisms are good moduli spaces;
			\item the horizontal morphisms are étale and affine and induce isomorphisms of stabiliser groups;
			\item if $\bar{w} \in [W/G](\CC)$ is the image of $w$ via the quotient morphism $W \to [W/G]$, then 
			the images of $\bar{w}$ via $[W/G] \to \cX$ and $[W/G] \to \cY$ are $p$ and $q$, respectively.
		\end{itemize}
	\end{lemma}
	
	\begin{proof}
		Since we are in characteristic zero, $G$ is smooth.
		Since $p$ is closed in $\cX$ and $\cX$ admits a good moduli space,
		by \cite[Proposition~12.14]{alper_good} $G$ is a linearly reductive algebraic group over $\CC$.
		The stacky germ isomorphism between $(\cX,p)$ and $(\cY,q)$
		is precisely condition~(1) of \cite[Theorem~4.19]{AHR20}.
		Therefore, by condition~(3) in \cite[Theorem~4.19]{AHR20} there exist an affine scheme $W$ of finite type over $\CC$ together with an action of $G$, a point $w \in W(\CC)$ fixed by $G$, and a diagram of \'etale morphisms
		\begin{equation} \label{eq:tetto_artin}
			\cX \overset{\Phi}\longleftarrow [W/G] \overset{\Psi}\longrightarrow \cY
		\end{equation}
		such that $\Phi(\bar{w})=p$ and $\Psi(\bar{w})=q$, and both $\Phi$ and $\Psi$ induce isomorphisms of stabiliser groups at $\bar{w}$.
		Here $\bar{w} \in [W/G](\CC)$ denotes the image of $w \in W(\CC)$ under the quotient morphism $W \to [W/G]$.
		
		Since $w$ is fixed by $G$, the point $\bar{w}$ is closed in $[W/G]$.
		Since $W$ is affine and $G$ is linearly reductive, the stack $[W/G]$ is cohomologically affine.
		Moreover, the diagonals of the three stacks $[W/G]$, $\cX$ and $\cY$ over $\CC$ are affine.
		By applying \cite[Proposition~3.2]{AHR20} successively to
		$\Phi$ and $\Psi$, we can find a cohomologically affine open neighbourhood $\cU$ of $\bar{w}$ in $[W/G]$ such that $\Phi \vert_\cU$ and $\Psi \vert_\cU$ are affine.
		By \cite[Tag~04YK]{stacks-project} such a $\cU$ must be of the type $[W'/G]$, where $W'$ is a $G$-invariant open neighbourhood of $w$ in $W$. Since $\cU$ is cohomologically affine, by \cite[Example~12.9]{alper_good} we have that $W'$ is affine.
		Therefore, replacing $W$ by $W'$ (and consequently replacing $[W/G]$ by $\cU = [W'/G]$, and $\Phi$ and $\Psi$ by their restrictions to $\cU = [W'/G]$), we may assume that both $\Phi$ and $\Psi$ are affine.
		In particular, $\Phi$ and $\Psi$ are separated and representable.
		
		By the universal property of good moduli spaces, the diagram~\eqref{eq:tetto_artin} extends to a commutative diagram
		\begin{equation} \label{eq:artin_common_etale}
			\xymatrix{
				\cX \ar[d] & [W/G] \ar[d]^\pi \ar[r]^\Psi \ar[l]_\Phi & \cY \ar[d] \\
				X & W/G \ar[r]^\psi \ar[l]_\phi & Y
			}
		\end{equation}
		where the vertical maps are  good moduli spaces. Here $W/G = \Spec \Gamma(W,\cO_W)^G$ is an affine scheme of finite type over $\CC$, by Hilbert's finite generation theorem on the invariant subalgebra.

		The morphism $\Phi$ is étale, separated and representable and induces an isomorphism of stabiliser groups at the closed point $\bar{w}$.
		By Luna's fundamental lemma \cite[Proposition~4.13]{AHR20} (or \cite[Theorem~3.14]{AHR26}), we can find an open neighbourhood $U$ of $\pi(\bar{w})$ in $W/G$ such that the restriction $\phi \vert_U$ is étale and $\pi^{-1}(U) \simeq U \times_X \cX$.
		By restricting to a principal open subset of the affine scheme $W/G$ which is contained in $U$ and contains $\pi(\bar{w})$, we can assume that $U$ is affine and the open immersion $U \into W/G$ is affine.
		Let $W'$ be the inverse image of $U$ along the quotient morphism $W \to W/G$.
		By replacing $W$ by $W'$ (and consequently replacing $[W/G]$ by $[W'/G]$, 
		replacing $W/G$ by $U=W'/G$, 
		replacing $\Phi$, $\Psi$, $\phi$, $\psi$ by the appropriate restrictions,
		and replacing $\pi$ by its base change to $U$), we can assume that the left square in \eqref{eq:artin_common_etale} is cartesian. (Notice that the base change of $\pi$ to $U$ is a good moduli space by \cite[Proposition~4.7(i)]{alper_good}.)
		
		By applying the same method to $\Psi$ and $\psi$, we can assume that also the right square in \eqref{eq:artin_common_etale} is cartesian.
		We have proved that both $\Phi$ and $\Psi$ are strongly \'etale according to \cite[Definition~3.13]{AHR26}.
		
		Since $X$, $W/G$, $Y$ are algebraic spaces, $\phi$ and $\psi$ trivially induce isomorphisms on stabiliser groups.
		By \cite[Lemma~0DUB]{stacks-project} the property of inducing isomorphisms on stabiliser groups is stable base change, hence it is satisfied by $\Phi$ and $\Psi$.
		This concludes the proof of Lemma~\ref{lemma:artin_local_structure_stacks_with_good_moduli}.
	\end{proof}
	
	\begin{lemma} \label{lem:local_structure_M_[S/G]}
		Let $\cM$ be an algebraic stack of finite type over $\CC$, with affine diagonal and admitting a good moduli space $M$.
		Let $p \in \cM(\CC)$ be a closed point.
		Let $G$ be an algebraic group over $\CC$ acting on an affine scheme $S$ of finite type over $\CC$ and let $s \in S(\CC)$ be a point fixed by $G$.
		Assume that there is an isomorphism of stacky germs
		\begin{equation} \label{eq:isom_stacky_germs_lemma}
			(\cM,p) \simeq (G \curvearrowright \widehat{\cO_{S,s}})
		\end{equation}
		where the formal action on $\widehat{\cO_{S,s}}$ is induced by the algebraic action on $S$.
		
		Then there exist an affine scheme $W$ of finite type over $\CC$ equipped with an action of $G$, a point $w \in W(\CC)$ fixed by $G$ and a commutative diagram
		\begin{equation}
			\label{eq:diagramma_3_quadrati_lemma}
			\xymatrix{
				& W \ar[d]_\rho \ar[r]^f & S \ar[d]^\sigma \\
				\cM \ar[d] & [W/G] \ar[d] \ar[r]^\Psi \ar[l]_\Phi & [S/G] \ar[d] \\
				M & W/G \ar[r] \ar[l] & S/G
			}
		\end{equation}
		where:
		\begin{itemize}
			\item $W/G = \Spec \Gamma(W,\cO_W)^G$ and $S/G = \Spec \Gamma(S,\cO_S)^G$ are affine schemes of finite type over $\CC$;
			\item the three squares are cartesian;
			
			\item the three vertical morphisms in the two bottom squares are good moduli spaces;
			
			\item the horizontal morphisms are étale and affine and induce isomorphisms on stabiliser groups;

			\item $\rho$ and $\sigma$ are the quotient morphisms with respect to the action of $G$;
			
			\item $f$ is $G$-equivariant and induces the morphism $\Psi$;
			
			\item $f(w) = s$ and $\Phi(\rho(w))=p$.
		\end{itemize}
	\end{lemma}
	
	\begin{proof}
		The stacky germ isomorphism \eqref{eq:isom_stacky_germs_lemma} implies that $G$ is isomorphic to $\Aut_\cM(p)$, which is linearly reductive because 	
		$p$ is a closed point of an algebraic stack admitting a good moduli space, by \cite[Proposition~12.14]{alper_good}.
		Hence $S/G = \Spec \Gamma(S,\cO_S)^G$ is the good moduli space of $[S/G]$ and is an affine scheme of finite type over $\CC$.
		
		Consider the quotient morphism $\sigma \colon S \to [S/G]$.
		Let $\bar{s} \in [S/G](\CC)$ be the image of $s \in S(\CC)$ via $\sigma$.
		By Example~\ref{ex:stacky_germ_quotient_stack} the stacky germ
		of $[S/G]$ at $\bar{s}$ is
		$(G\curvearrowright\widehat{\cO_{S,s}})$.
		Thus the isomorphism \eqref{eq:isom_stacky_germs_lemma} implies that the stacky germs
		$(\cM,p)$ and $([S/G], \bar{s})$ are isomorphic.
		Therefore, by Lemma~\ref{lemma:artin_local_structure_stacks_with_good_moduli},
		we can find an affine scheme $Z$ of finite type over $\CC$ equipped with an action of $G$,
		a point $z \in Z(\CC)$ fixed by $G$,
		and a commutative diagram
		\begin{equation*}
			\xymatrix{
				\cM \ar[d] & [Z/G] \ar[d] \ar[r]^{\Psi'}\ar[l]_{\Phi'} & [S/G] \ar[d] \\
				M & Z/G \ar[r]^{\psi'} \ar[l]_{\phi'} & S/G
			}
		\end{equation*}
		where the squares are cartesian, 
		the vertical maps are good moduli spaces,
		the horizontal maps are étale and affine,
		and
		if $\bar{z} \in [Z/G](\CC)$ denotes the image of $z$ via the quotient morphism $Z \to [Z/G]$ then $\Phi'(\bar{z})=p$ and $\Psi'(\bar{z}) = \bar{s}$.
		
		Set $W := S \times_{[S/G]} [Z/G]$,
		where the fibred product is taken with respect to the quotient morphism $\sigma \colon S \to [S/G]$ and $\Psi' \colon [Z/G] \to [S/G]$.
		We have the cartesian square
		\begin{equation*}	
			\xymatrix{
				W \ar[r]^f \ar[d]_\tau & S \ar[d]^\sigma \\
				[Z/G] \ar[r]^{\Psi'} & [S/G] 
			}
		\end{equation*}
		where $f$ is affine and étale, because $\Psi'$ is.
		Since $S$ is an affine scheme of finite type over $\CC$, we deduce that $W$ is an affine scheme of finite type over
		$\CC$.
		
		Since $\sigma \colon S \to [S/G]$ is a $G$-torsor, also its base change $\tau \colon W \to [Z/G]$ is a $G$-torsor.
		Therefore
		$W$ has a $G$-action for which $f$ is
		$G$-equivariant, and there is an equivalence
		\[
		\epsi \colon[W/G]\xrightarrow{\sim}[Z/G].
		\]
		Let $\rho \colon W \to [W/G]$ be the quotient morphism. 
		The
		equivalence $\epsi$ comes with an isomorphism
		$\epsi \circ \rho \simeq \tau$.
		Moreover, the morphism $\Psi \colon [W/G] \to [S/G]$ induced
		by the $G$-equivariant morphism $f \colon W \to S$ satisfies
		$\Psi \simeq \Psi' \circ\epsi$.
		Furthermore, the affine schemes $Z/G$ and $W/G$ are isomorphic as they are the good moduli spaces of the corresponding quotient stacks.
		With the equivalence $\epsi$ it is easy to construct the diagram~\eqref{eq:diagramma_3_quadrati_lemma}:
		the morphisms in the middle row are $\Phi = \Phi' \circ \epsi$,
		$\Psi = \Psi' \circ \epsi$ and the morphisms in the bottom row are the compositions of $\phi'$, $\psi'$ with the isomorphism $W/G \simeq Z/G$.
		
		We choose $w \in W(\CC)$ corresponding to $(s, \bar{z}) \in  S(\CC) \times [Z/G](\CC)$.
	\end{proof}

	\begin{proof}[Proof of Theorem~\ref{thm:artin_stack}]
		By Lemma~\ref{lem:local_structure_M_[S/G]}
		we find an affine scheme $W$, a point $w \in W(\CC)$ and a diagram like \eqref{eq:diagramma_3_quadrati_lemma}.
		We omit to write down all the properties listed in the statement of Lemma~\ref{lem:local_structure_M_[S/G]}.
		
		Let $\cU \subseteq \cM$ be an arbitrary open substack containing $p$.
		We must show that $\cU$ contains a point $p' \in \cM(\CC)$ which is closed in $\cM$ and such that the stacky germ $(\cM,p')$ is isomorphic to $(G' \curvearrowright R')$.
		
		Since $f$ is étale, $f$ is open. Therefore
		\[
		f\bigl((\Phi\circ\rho)^{-1}(\cU)\bigr) = f \bigl(  \rho^{-1} (  \Phi^{-1}(\cU)  )   \bigr)
		\]
		is an open neighbourhood of $s$ in $S$.
		By hypothesis there exists $s' \in S'(\CC)$ such that
		\[
		s' \in f\bigl((\Phi\circ\rho)^{-1}(\cU)\bigr)
		\]
		and $s'$ satisfies the conditions (1)-(4).
		In particular, there exists $w' \in W(\CC)$ such that $f(w')=s'$
		and
		\[
		w' \in (\Phi\circ\rho)^{-1}(\cU).
		\]
		Set $p' = \Phi (\rho(w')) \in \cU(\CC) \subseteq \cM(\CC)$.
		We need to show that $p'$ satisfies the required conditions.
		
		Let us start by proving that $p'$ is closed in $\cM$.
		By condition (1) the orbit $G \cdot s'$ is closed in $S$, therefore $\sigma(s')$ is closed in $[S/G]$.
		By applying Lemma~\ref{lem:chiusura_punti_morfismo_da_spazi_algebrici} to $\Psi$ (which comes from base change from a morphism between affine schemes of finite type over $\CC$) we obtain that $\rho(w')$ is closed in $[W/G]$.
		If we could apply Lemma~\ref{lem:chiusura_punti_morfismo_da_spazi_algebrici} to $\Phi$ (which comes from base change from a morphism between algebraic spaces of finite type over $\CC$) we would obtain that $p'$ is closed in $\cM$.
		We need to show that the good moduli space $M$ of $\cM$ is quasi-separated.
		Let $\mu \colon \cM \to M$ denote the good moduli space morphism.
		Since $\cM$ is of finite type over $\CC$, by \cite[Theorem~4.16(i),(xi)]{alper_good}, $\mu$ is surjective and quasi-compact and $M$ is of finite type over $\CC$. Let $\Delta_\cM$ (resp.\ $\Delta_M$) denote the diagonal of $\cM$ (resp.\ $M$) over $\CC$. Since $\mu$ is quasi-compact, we get that $\mu \times \mu$ is quasi-compact.
		Since $\Delta_\cM$ is affine by hypothesis, we have that $\Delta_\cM$ is quasi-compact and hence
		the composition $\Delta_M \circ \mu = (\mu \times \mu) \circ \Delta_\cM$ is quasi-compact.
		Since $\mu$ is surjective, by \cite[Lemma~050X]{stacks-project} we have that $\Delta_M$ is quasi-compact, i.e.\ $M$ is quasi-separated.

		Since $\Phi$ and $\Psi$ come from base change from morphisms between algebraic spaces, by \cite[Lemma~0DUB]{stacks-project} we have that they induce isomorphisms on stabiliser groups.
		From condition~(2) we deduce that the stabiliser group of $\sigma(s')$ is $G'$ and hence the stabiliser group of $p'$ is $G'$.
		
		Since $\Phi$ and $\Psi$ are étale and induce an isomorphism on stabiliser groups at each point of $[W/G]$, by Lemma~\ref{lem:etale_induce_same_stacky_germ} we obtain the following isomorphisms of stacky germs:
		\begin{equation} \label{eq:isom_stacky_germ1}
			(\cM,p')
			\overset{\Phi}\simeq
			\left( [W/G], \rho(w')  \right)
			\overset{\Psi}\simeq
			\left( [S/G], \sigma(s')  \right)
		\end{equation}
		Therefore it remains to show that the stacky germ $( [S/G], \sigma(s') )$ is isomorphic to $(G' \curvearrowright R')$.

		Consider the morphism
		\[
		\lambda \colon[S'/G'] \to [S/G]
		\]
		induced by the closed immersion $G'\into G$ and the $G'$-equivariant closed
		immersion $S'\into S$. There is a cartesian diagram
		\[
		\xymatrix{
			G\times^{G'}S' \ar[r]^\al \ar[d]
			& S\ar[d]^\sigma\\
			[S'/G']\ar[r]^\lambda
			& [S/G]
		}
		\]
		where $\al$ is the natural morphism given by $[(g,u)] \mapsto g \cdot u$.
		Let $\sigma' \colon S' \to [S'/G']$ be the quotient morphism.
		
		Let $E\subseteq G\times^{G'}S'$ be the étale locus of $\al$.
		Since $\al$ is $G$-equivariant, $E$ is $G$-invariant.
		Condition~(3) implies that $[(e,s')]\in E$; hence $E$ contains
		the preimage of $\sigma'(s')$ with respect to $G \times^{G'} S' \to [S'/G']$.
		Since $\al$ is the base change of $\lambda$ along the smooth surjective morphism $\sigma$ and
		since étaleness	descends under smooth surjective base change, it follows that
		$\lambda$ is étale at $\sigma'(s')$.
		
		It is clear that $\lambda$ induces an isomorphism on stabiliser groups at $\sigma'(s')$.
		Moreover $\sigma'(s')$ is closed in $[S'/G']$ because $s'$ is $G'$-fixed and $\sigma(s')$ is closed in $[S/G]$ because the orbit $G \cdot s'$ is closed in $S$.
		Therefore by Lemma~\ref{lem:etale_induce_same_stacky_germ} applied to $\lambda$
		we obtain the following isomorphisms of stacky germs
		\begin{equation} \label{eq:isom_stacky_germ2}
			(G' \curvearrowright R')
			\overset{(4)}{\simeq}
			(G' \curvearrowright \widehat{\cO_{S',s'}})
			\simeq
			\left( [S'/G'], \sigma'(s') \right)
			\overset{\lambda}{\simeq}
			\left( [S/G], \sigma(s') \right),
		\end{equation}
		where the isomorphism in the middle holds by Example~\ref{ex:stacky_germ_quotient_stack}.
		By combining the isomorphisms \eqref{eq:isom_stacky_germ1} and \eqref{eq:isom_stacky_germ2} we conclude that the stacky germ $(\cM,p')$ is isomorphic to $(G' \curvearrowright R')$.
		This concludes the proof of Theorem~\ref{thm:artin_stack}.
	\end{proof}
	
	With a similar strategy we can also prove the existence of points which are not closed:
	
	\begin{proposition} \label{prop:ultima_prop}
		Let $\cM$ be an algebraic stack of finite type over $\CC$, with affine diagonal and admitting a good moduli space. Let $p \in \cM(\CC)$ be a closed point.
		
		Let $G$ be an algebraic group over $\CC$ acting on an affine scheme $S$ of finite type over $\CC$ and let $s \in S(\CC)$ be a point which is fixed by $G$.
		Assume that there is an isomorphism of stacky germs
		\begin{equation*}
			(\cM,p) \simeq (G \curvearrowright \widehat{\cO_{S,s}})
		\end{equation*}
		where the formal action on $\widehat{\cO_{S,s}}$ is induced by the algebraic action on $S$.
		
		Assume that every open neighbourhood of $s$ in $S$ contains at least one point $s' \in S(\CC)$ such that the orbit $G \cdot s'$ is not closed in $S$.
		
		Then every open neighbourhood of $p$ in $\cM$ contains a point $p' \in \cM(\CC)$ which is not closed in $\cM$.
	\end{proposition}
	
	\begin{proof}
		The proof is very similar to the proof of Theorem~\ref{thm:artin_stack}, but shorter.
		We can find an affine scheme $W$, a point $w \in W(\CC)$ and a diagram like \eqref{eq:diagramma_3_quadrati_lemma} satisfying all the properties listed in the statement of Lemma~\ref{lem:local_structure_M_[S/G]}.
		
		Let $\cU \subseteq \cM$ be an arbitrary open substack containing $p$.
		Since $f$ is étale, $f$ is open. Therefore $f\bigl((\Phi\circ\rho)^{-1}(\cU)\bigr)$
		is an open neighbourhood of $s$ in $S$.
		By hypothesis there exists a point $s' \in S(\CC)$ such that
		\[
		s' \in f\bigl((\Phi\circ\rho)^{-1}(\cU)\bigr)
		\]
		and the orbit $G \cdot s'$ is not closed in $S$.
		In particular, there exists $w' \in W(\CC)$ such that $f(w')=s'$
		and
		\[
		w' \in (\Phi\circ\rho)^{-1}(\cU).
		\]
		Set $p' = \Phi (\rho(w')) \in \cU(\CC) \subseteq \cM(\CC)$.
		
		Since the orbit $G \cdot s'$ is not closed in $S$, the point $\sigma(s')$ is not closed in $[S/G]$.
		By Lemma~\ref{lem:chiusura_punti_morfismo_da_spazi_algebrici} $\rho(w')$ is not closed in $[W/G]$.
		By Lemma~\ref{lem:chiusura_punti_morfismo_da_spazi_algebrici} $p'$ is not closed in $\cM$.
		Here we have used that the good moduli space of $\cM$ is quasi-separated, as proved in the proof of Theorem~\ref{thm:artin_stack}.
	\end{proof}
	
	\begin{example}
		We continue Example~\ref{exa:Gm_agisce_pesi_1-1}.
		Let $\cM$ be an algebraic stack of finite type over $\CC$, with affine diagonal and admitting a good moduli space. If $p \in \cM(\CC)$ is a closed point such that the stacky germ $(\cM,p)$ is isomorphic to $(\Gm \curvearrowright \CC \pow{t_1,t_2})$ where the formal action is given by the weights $1,-1$, then we can deduce that in every open neighbourhood of $p$ in $\cM$ there exists a non-closed point $p' \in \cM(\CC)$.
		
		Indeed, in Proposition~\ref{prop:ultima_prop} we can use the affine plane $S = \Spec \CC [t_1,t_2]$ with the action of $\Gm$ given by the weights $1,-1$, and the origin as $s$:
		points on the two axes different from the origin have non-closed orbit.
	\end{example}
	
	\bibliography{Biblio}
	
\end{document}